%    \newcount\chapno
 \magnification=\magstep1 
\newcount\sectno
\newcount\subsectno
\newcount\parno
\newcount\equationno
\newif\ifsubsections
\subsectionsfalse

\def\sectnum{\the\sectno} 
\def\subsectnum{\sectnum\ifsubsections .\the\subsectno\fi} 
\def\parnum{\subsectnum .{\the\parno}}
\def\eqnum{\subsectnum .\the\equationno}

\def\abstract#1\endabstract
{
{\abstractfont
    \baselineskip=9pt
    \leftskip=4pc  \rightskip=4pc
    \bigskip
    \noindent
     ABSTRACT.\ #1
\medskip} 
}
 
\def\thanks[#1]#2\endthanks{\footnote{$^#1$}{\footnotefont\kern-6pt #2}}

\newcount\minutes
\newcount\scratch

\def\timestamp{%
\scratch=\time
\divide\scratch by 60
\edef\hours{\the\scratch}
\multiply\scratch by 60
\minutes=\time
\advance\minutes by -\scratch
\the \month/\the\day$\,$---$\,$\hours:\null
\ifnum\minutes< 10 0\fi
\the\minutes}

\def\today{\ifcase\month\or
January\or February\or March\or April\or May\or June\or
July\or August\or September\or October\or November\or December\fi
\space\number\day,\number\year}

\outer\def\newsection #1.\par{\vskip1.5pc plus.75pc \penalty-250
     \subsectno=0
     \parno=0
     \equationno=0
     \advance\sectno by1
     \leftline{\smalltitlefont \sectnum.\hskip 1pc  #1}
                \nobreak \vskip.75pc\noindent}

\outer\def\newsectiontwoline #1/#2/.\par{\vskip1.5pc plus.75pc \penalty-250
     \subsectno=0
     \parno=0
     \equationno=0
     \advance\sectno by1
     \leftline{\smalltitlefont \sectnum.\hskip 1pc  #1}
     \leftline{\smalltitlefont \hskip 22pt #2}
                \nobreak \vskip.75pc\noindent}

\outer\def\newsubsection #1.\par{\vskip1pc plus.5pc\penalty-250
     \parno=0
     \equationno=0
     \advance\subsectno by1
     \leftline{{\bf \subsectnum}\hskip 1pc  #1.}
                \nobreak \vskip.5pc\noindent}

\def\newpar #1.{\advance \parno by1
     \par
 \medbreak \noindent 
      {\bf \parnum. #1.} \hskip 6pt}

\long\def \newclaim #1. #2\par {\advance \parno by1
    \medbreak \noindent 
     {\bf \parnum \hskip 6 pt #1.\hskip 6pt} {\sl #2} \par \medbreak}

\def\eq $$#1$${\global \advance \equationno by1 $$#1\eqno(\eqnum)$$}

\def\rmarginsert[#1]{\hglue 0pt\vadjust
{\null\vskip -\baselineskip\rightline{\abstractfont\rlap{\hfil\  #1}}}}

\def\lmarginsert[#1]{\hglue 0pt\vadjust
{\null\vskip -\baselineskip\leftline{\abstractfont\llap{#1\ \hfill}}}}

\newif\ifproofmode
\proofmodefalse

\def\refpar[#1]#2.{\advance \parno by1
     \par
 \medbreak \noindent 
      {\bf \parnum \hskip 6 pt #2.\hskip 6pt}%
\expandafter\edef\csname ref#1\endcsname
{\parnum}\ifproofmode\rmarginsert[\string\ref#1]\fi}

\long\def \refclaim[#1]#2. #3\par {\advance \parno by1
    \medbreak \noindent 
{\bf \parnum \hskip 6 pt #2.\hskip 6pt}%
\expandafter\edef\csname ref#1\endcsname
{\parnum}\ifproofmode\rmarginsert[\string\ref#1]\fi
{\sl #3} \par \medbreak}

\def\refer[#1]{%
\expandafter\xdef\csname ref#1\endcsname
{\parnum}\ifproofmode\rmarginsert[\string\ref#1]\fi}

\def\refereq[#1]$$#2$$ {%
\eq$$#2$$%
\expandafter\xdef\csname ref#1\endcsname{(\eqnum)}%
\ifproofmode\rmarginsert[\string\ref#1]\fi
}

\def\refeq{\refereq}

\def \Definition #1\\ {\vskip 1pc \noindent 
      {\bf #1. Definition. \hskip 6pt}\vskip 1pc}

\def\proof{{PROOF.} \enspace}

\def\qedmark{\hbox{\vrule height 4pt width 3pt}}
\def\qedskip{\vrule height 4pt width 0pt depth 1pc}
\def\qed{\penalty 1000\quad\penalty 1000{\qedmark\qedskip}}

\def \a {\alpha}
\def \b {\beta}

\def \d {\delta}
\def \D {\triangle}

\def \g {\gamma}

\def \K {\nabla}
\def \l {\lambda}

\def \n {\,\vert\,}
\def \N {\,\Vert\,}
\def \o {\theta}
\mathchardef\p="011E    %straight phi

\def \s {\sigma}

\def \W {\Omega}

     %Second Fundamental Form

\def\Gtwo{{\mathop{{{\mbi G\/}}\kern-.5pt_{{}_2}}}}
\def\Ffour{{\mathop{{{\mbi F\/}}\kern-2.5pt_{{}_4}}}}
\def\Esix{{\mathop{{{\mbi E\/}}\kern-.5pt_{{}_6}}}}
\def\Eseven{{\mathop{{{\mbi E\/}}\kern-.5pt_{{}_7}}}}
\def\Eeight{{\mathop{{{\mbi E\/}}\kern-.5pt_{{}_8}}}}

%% Lie Algebra Symbols

\def\ca{{\liefont A}}

\def\cb{{\liefont B}}

\def\cc{{\liefont C}}

\def\cf{{\liefont F}}

\def\cg{{\liefont G}}

\def\ch{{\liefont H}}

\def\ck{{\liefont K}}

\def\cl{{\liefont L}}

\def\cm{{\liefont M}}

\def\cn{{\liefont N}}

\def\co{{\liefont O}}

\def\cp{{\liefont P}}

\def\cs{{\liefont S}}

\def\ct{{\liefont T}}

\def\cu{{\liefont U}}

\def\sdp{
\mathop{\hbox{$\raise 1pt\hbox{$\scriptscriptstyle |$}\kern-2.5pt\times$}}
        }

\def\dag{{\raise 1 pt\hbox{{$\scriptscriptstyle \dagger$}}}}

\def\*{\raise 1.5pt \hbox{*}}

\def\Im{\mathop{\rm Im}\nolimits}

\def\exp{\mathop{\rm exp}\nolimits}

\def\tr{\mathop{\rm tr}\nolimits}

\def\diag{\mathop{\rm diag}\nolimits}

\def\Ad{\mathop{\rm Ad}\nolimits}
\def\ad{\mathop{\rm ad}\nolimits}

\def\diag {\mathop{\rm diag}\nolimits}

\def\id {\mathop{\rm id}\nolimits}

 % for top of page
                                                          
\def\ifundefined#1{\expandafter\ifx\csname#1\endcsname\relax}

\def\cross{\times}
\def\longerrightarrow{-\kern-5pt\longrightarrow}

\def\star{\lower 1pt\hbox{*}}
\def \nulset {
\raise 1pt\hbox{
\hskip -3pt$\not$\kern -0.2pt \raise .7pt\hbox{${\scriptstyle\bigcirc}$}}}
\def \norm|#1|{\Vert#1\Vert}

\def \interior(#1){#1\kern -6pt \raise 7.5pt 
      \hbox{$\scriptstyle \circ$}{}\hskip 2pt}

\def\twist_#1{\kern -.15em\cross\kern -.30 em{}_{{}_{#1}}\kern .07 em}

\font\cmr=cmr10 at 10pt
\font\cmrviii=cmr8
\font\cmrvi=cmr6

\font\cmrXIV=cmr12 at 14 pt
\font\cmrXX=cmr17 at 20 pt
\font\cmrXXIV=cmr17 at 24 pt
\font\cmbxXII=cmbx12
\font\cmbxsl=cmbxsl10
\font\cmbxslviii=cmbxsl10 at 8pt
\font\cmbxslv=cmbxsl10 at 5pt

      \font\tenrm=cmr10 at 10.3 pt
      \font\sevenrm=cmr7 at 7.21 pt
      \font\fiverm=cmr5 at 5.15 pt
      \font\teni=cmmi10 at 10.3 pt
      \font\seveni=cmmi7 at 7.21 pt
      \font\fivei=cmmi5 at 5.15 pt     
      \font\tensy=cmsy10 at 10.3 pt
      \font\sevensy=cmsy7 at 7.21 pt   
      \font\fivesy=cmsy5 at 5.15 pt    
      \font\tenex=cmex10 at 10.3 pt
      \font\tenbf=cmbx10 at 10.3 pt
      \font\sevenbf=cmbx7 at 7.21 pt    
      \font\fivebf=cmbx5 at 5.15 pt  

\def\UseComputerModern   %%%  This is the default plain TeX choice
{
\textfont0=\tenrm \scriptfont0=\sevenrm \scriptscriptfont0=\fiverm
\def\rm{\fam0\tenrm}
\textfont1=\teni \scriptfont1=\seveni \scriptscriptfont1=\fivei
\def\mit{\fam1} \def\oldstyle{\fam1\teni}
\textfont2=\tensy \scriptfont2=\sevensy \scriptscriptfont2=\fivesy
\def\cal{\fam2}
\textfont3=\tenex \scriptfont3=\tenex \scriptscriptfont3=\tenex
\def\it{\fam\itfam\tenit} % \it is family 4
\textfont\itfam=\tenit
\def\sl{\fam\slfam\tensl} % \sl is family 5
\textfont\slfam=\tensl
\def\bf{\fam\bffam\tenbf} % \bf is family 6
\textfont\bffam=\tenbf \scriptfont\bffam=\sevenbf
\scriptscriptfont\bffam=\fivebf
\def\tt{\fam\ttfam\tentt} % \tt is family 7
\textfont\ttfam=\tentt
\def\abstractfont{\cmrviii}
\def\footnotefont{\cmrviii}
\def\tinyfont{\cmrvi}
\def\smalltitlefont{\cmbxXII}  %added 2/11/94
\def\titlefont{\cmrXIV}
\def\bigtitlefont{\cmrXX}
\def\verybigtitlefont{\cmrXXIV} 
\textfont9=\cmbxsl \scriptfont9=\cmbxslviii \scriptscriptfont9=\cmbxslv
\def\mbi{\fam9}
\cmr
}  

\def\liefont{\cal}

\def \bs {\bigskip}
\def \ms {\medskip}
\def \ss {\smallskip}

\def \ni {\noindent}

\def\enditem{\item{}\par\vskip-\baselineskip\noindent}
\def\ei{\enditem}

   \baselineskip=14 true pt 
   \hsize 37 true pc \hoffset= 22 true pt
   \voffset= 0 true pt
   \vsize  54 true pc

\def\id{\mathop{\rm id}\nolimits}

\def\ni{\noindent}
\def\bs{\bigskip}
\def\ms{\medskip}

\def\ss{\smallskip}
\def\ct{{\cal T}}

   \baselineskip=14 true pt 
   \hsize 35 true pc \hoffset= 25 true pt
   \vsize  52 true pc
\UseComputerModern
\subsectionsfalse
%\font\smalltitlefont=TimesB at 11 pt
\font\smalltitlefont=cmbx10 at 11 pt
%\proofmodetrue
%\margindate
%\rightline{\today}\ss
%

%\centerline{\mytitlefont Title of paper}
%\bs

\def\Bibliography
{

\font\TRten=cmr10 at 10 true pt
\font\TIten=cmti10 at 10 true pt
\font\TBten=cmbx10 at 10 true pt

\def\ourindent{\hfil\vskip-\baselineskip}

 \frenchspacing
 \parindent=0pt

 %fonts for various fields in a bibliographic reference
 \def \keyfnt{\TRten}
 \def \authornamefnt{\TRten}
 \def \booktitlefnt{\TIten}
 \def \articletitlefnt{\TRten}
 \def \journalnamefnt{\TIten}
 \def \volumefnt{\TBten}
 \def \publishernamefnt{\TRten}
 \def \pagesfnt{\TRten}
 \def \yearfnt{\TRten}
 \def \commentfnt{\TRten}

 \def \bookitem //##1//##2//##3//##4//##5//##6//##7//##8//
      %#1=key,#2=author,#3=title,#4=publisher,#6=year,#8=comment
      { \goodbreak{\par\hskip-40pt{\keyfnt [##1]}\ourindent{\authornamefnt ##2,}}
              {\booktitlefnt ##3.\/}\thinspace
              {\publishernamefnt ##4,}
              {\yearfnt ##6.}
              {\commentfnt ##8}
       }

\def \b{\bookitem}

\def \articleitem //##1//##2//##3//##4//##5//##6//##7//##8//
      %#1=key,#2=author,#3=title,#4=journal,#5=volume,#6=year,#7=pages,#8=comment
      { \goodbreak{\par\hskip-40pt{\keyfnt [##1]}\ourindent{\authornamefnt ##2,}}
              {\articletitlefnt ##3},
              {\journalnamefnt ##4\/}
              {\volumefnt ##5}
              {\hbox{\yearfnt(\hskip -1pt ##6)}},
              {\pagesfnt ##7.}
              {\commentfnt ##8}
       }
\def \a{\articleitem}

\def \preprintitem //##1//##2//##3//##4//##5//##6//##7//##8//
      %#1=key,#2=author,#3=title,#8=comment
      { \goodbreak{\par\hskip-40pt{\keyfnt [##1]}\ourindent{\authornamefnt ##2,}}
              {\articletitlefnt ##3},
              {\commentfnt ##8}
       }
\def \p{\preprintitem}

   \vskip 1in
   \centerline{References}
   \vskip .5in
}

% end of bibliographic macros

\def\bu{\bullet}

%\input cldefs
%\today
%\proofmodetrue
\def\bu{\bullet}
\def\p{\partial}
\font\bfs=cmb10 at 9 pt

\def\bu{\bullet}

\font\cmbxXIV=cmbx12 at 14pt

\font\cmbxXIV=cmbx12 at 14pt

{\cmbxXIV 
\centerline {Poisson Actions and Scattering Theory} 
\centerline {for Integrable
Systems}
}
\bs\bs
\centerline { Chuu-Lian Terng\footnote{$^1$} {Research supported in part by 
NSF Grant DMS 9626130 }
and  Karen Uhlenbeck\footnote{$^2$}{Research supported in part by  Sid
Richardson Regents' Chair Funds, University of Texas system}}

\font\bfs=cmb10 at 10 pt
\font\cmrviii=cmr8

\bs\bs
\centerline {\bfs Abstract\/}
\bs\bs
{\cmrviii  
Conservation laws, heirarchies, scattering theory and B\"acklund
transformations are known to be the building blocks of integrable partial
differential equations. We identify these as facets of a theory of
 Poisson group actions, and apply the theory to the
ZS-AKNS nxn heirarchy (which includes the non-linear Schr\"odinger
equation, modified KdV, and the n-wave equation). 
We first find a simple model Poisson group action that contains flows for systems
with a Lax pair whose terms all decay on $R$.  B\"acklund transformations
and flows arise from subgroups of this single Poisson group.   For the ZS-AKNS nxn
heirarchy defined by a constant $a\in u(n)$, the simple model is no longer correct.
The $a$ determines two separate Poisson structures. The flows come from the
Poisson action of the centralizer $H_a$ of $a$ in the dual Poisson group (due to the
behavior of  $e^{a\l x}$ at infinity).  When
$a$ has distinct eigenvalues, $H_a$ is abelian and it acts symplectically.  The phase
space of these flows is the space $S_a$ of left cosets of the centralizer of $a$ in
$D_-$, where $D_-$ is a certain loop group. The group $D_-$ contains both a Poisson
subgroup corresponding to the continuous scattering data, and a rational loop
group corresponding to the discrete scattering data.  The $H_a$-action is the right
dressing action on $S_a$.  B\"acklund transformations arise from the
action of the simple rational loops on $S_a$ by right multiplication. 
Various geometric equations arise from appropriate choice of $a$ and 
restrictions of the phase space and flows.  In particuar, we discuss 
applications to the sine-Gordon equation, harmonic maps, Schr\"odinger
flows on symmetric spaces, Darboux orthogonal coordinates, and
isometric immersions of one space-form in another. 
 }

\vfil\eject

\centerline {\bf Table of Contents\/}
\ms

\+\kern .6in &\kern .3in & \cr
\+ & 1. & Introduction\cr
\+ & 2. & Poisson Actions\cr
\+ & 3. & Negative flows in the decay case  \cr
\+ & 4.  & Poisson structure for negative flows (decay case) \cr
\+ & 5.  &  Positive flows in the asymptotically constant case\cr
\+ & 6.  & Action of the rational loop group \cr
\+ & 7. & Scattering data and Birkhoff decompositions  \cr
\+ & 8 & Poisson structure for the positive flows (asymptotic case) \cr
\+ & 9. & Symplectic structures for the restricted case \cr
\+ & 10. &  B\"acklund transformations for the $j$-th flow\cr
\+ &11.  & Geometric Non-linear Schr\"odinger equation  \cr
\+ &12.  & First flows and flat metrics \cr

\bs

\newsection Introducton.\par

       Soliton theory is an enticingly elegant
part of modern mathematics. It has a multitude of interpretations in geometry, 
analysis and algebra. The main goal of this paper is to relate loop groups actions,
scattering theory, and B\"acklund transformations within the same narrative, via
Poisson actions. Our work is motivated by Beals and Coifman's  rigorous and
beautiful treatment of scattering and inverse scattering theory of the first order
systems ([BC 1, 2, 3]).  An expository version of our main
result on scattering theory is contained in lecture notes by Richard Palais [Pa].
In retrospect,, we also find that many of our results in the $su(2)$ case are
contained in the book by Faddeev and Takhtajian ([FT]).  Throughout the paper, 
the matrix non-linear Schr\"odinger equation
is used as a motivating example. In the final section we discuss a number of
applications in geometry, including Darboux orthogonal coordinates and isometric
immerisons of space forms in space forms.  We begin the paper with a survey.

\ms
\ni $\bu$ \hskip 6pt {\bfs Finite dimensional mechanics\/}
\ms

To give some perspective, we start with a short review of finite dimensional Hamiltonian
systems, complete integrability, symplectic actions, Poisson actions, and moment maps.  A
more detailed review of Poisson actions is given in section 2.  A {\it symplectic structure\/}
on a $2n$-dimensional manifold
$M$ is a closed, non-degenerate two form $w$ on $M$. Since $w$ is non-degenerate, it induces
an isomorphism $J:T^*M\to TM$. A Hamiltonian on $M$ is a smooth funtion $f:M\to R$.  The
Hamiltonian vector field
$X_f$ corresponding to
$f$ is the symplectic dual of $df$, i.e., 
$$X_f=J(df), \qquad {\rm or\/}\qquad i_{X_f}w= df.$$  It follows from this definition
that $X_f$ is symplectic, i.e.,
$L_{X_f} w=0$, or equivalently the one parameter subgroup generated by $X_f$
preserves $w$. 

A Poisson structure on $M$ is a Lie bracket $\{\, , \,\}$ on $C^\infty(M,R)$, which satisfies the
 Leibnitz rule $$\{fg,h\}=f\{g,h\} + g\{f,h\}.$$ A symplectic form $w$ induces a natural Poisson
structure on
$M$  by 
$$\{f,g\}=w(X_f,X_g)= df(X_g).$$  Then the map from $C^\infty(M,R)$ to the Lie algebra of
vector fields on $M$ defined by
$f\mapsto X_f$ is a Lie algebra homomorphism, i.e., $$[X_f,X_g]= X_{\{f,g\}}.$$  Two
Hamiltonians $f,g$ are said to be {\it in involution\/} if $\{f,g\}=0$. In this case the
corresponding Hamiltonian flows commute, and $g$ is a {\it conservation law\/} for the
Hamiltonian system 
\refeq[de]$${dx\over dt}=X_f(x(t)),$$ i.e., $g$ is constant on the integral curves of $X_f$. 

The Hamiltonian system \refde{} on $M^{2n}$ is called {\it completely integrable\/} if
there exists $n$ conservation laws $f_1=f, f_2, \cdots, f_n$ that are in involution and $df_1,
\cdots, df_n$ are linearly independent. For example, the Hamiltonian systems given by the
Kowalevsky top, the Toda system and the geodesic flow on an ellipsoid are completely
integrable.

Suppose $f$ is completely integrable and  $f=f_1, \cdots, f_n$ are in involution. Then $X_{f_1},
\cdots, X_{f_n}$ generate an action of $R^n$ on $M$.  If the map
$\mu=(f_1, \cdots, f_n):M\to R^n$ is proper, then $X_{f_1}, \cdots, X_{f_n}$  
generate an action of the $n$-torus $T^n$ on each $R^n$-orbit of $M$. Let $\o_1,
\cdots, \o_n$ denote the angular coordinates on the torus orbit. Then $(f_1,\cdots,
f_n, \o_1,
\cdots, \o_n)$ is a coordinate system on $M$ and the system \refde{} is linearized
in these coordinates. These are the {\it action-angle coordinates\/} in Liouville's
Theorem (for detail see [AbM], [Ar]). 

The notion of complete integrability can be extended to the notion of symplectic 
action of a Lie group. An action of $G$ on $(M,w)$ is {\it symplectic\/} if the action
preserves $w$.  A symplectic action of a Lie group $G$ on
$M$ is {\it Hamiltonian\/} if there exists a map $\mu:M\to \cg^*$ such that the
infinitesimal vector field on $M$ corresponding to $\xi\in \cg$ is the  Hamiltonian
vector field of the function
$f_\xi$ defined by $f_{\xi}(x)=\mu(x)(\xi)$. Such $\mu$ is called a  {\it moment
map\/}. When
$G$ is abelian, the flows generated by the action commute.  In particular,  the study
of completely integrable systems on
$M^{2n}$ is the same as the study of Hamiltonian actions of $R^n$ on $M^{2n}$.  When $G$ is
non-abelian,
 the flows generated by $\eta$ in the centralizer of $\xi$ $\cg_\xi=\{\eta\in \cg\n
[\xi,\eta]=0\}$ commute with the flow generated by $\xi$. In other words, $f_\eta$ is a
conservation law for the flow generated by $\xi$.  But the flows generated by $\eta_1,
\eta_2\in \cg_\xi$ in general do not commute.  

\ms
\ni $\bu$ \hskip 6pt {\bfs Poisson groups\/}
\ms

  Given two Poisson manifolds $(M_1, \{\,,\,\}_1)$ and $(M_2, \{\,,\,\}_2)$, the {\it product
Poisson  structure\/} on $M_1\times M_2$ is defined by 
$$\{f,g\}(x,y)=\{f(\cdot, y), g(\cdot, y)\}_1 (x) + \{f(x,\cdot), g(x,\cdot)\}(y).$$ A map
$\phi:M_1\to M_2$ is Poisson if
$\phi$ preserves the Poisson structure, i.e., $$\{f,g\}_2\circ \phi= \{f\circ \phi,
g\circ\phi\}_1.$$   A Poisson group is a Poisson manifold $(G,\{\,,\,\})$ such that 
$G$ is a Lie group and the multiplication map $m:G\times G\to
G$ is a Poisson map, where $G\times G$ is equipped with the product Poisson  
structure. The modern study of Poisson groups was initiated by Drinfe\'ld in
[Dr] and there are several good articles by Lu and Weinstein [LW] and
Semenov-Tian-Shanksy 
 [Se].  Given a Poisson group
$G$, there is a canonical construction of a dual Poisson group
$G^*$ (cf. [LW]). The simplest Poisson group is a Lie group $G$ with the trivial
Poisson structure, and its dual Poisson group is the dual
$\cg^*$ of Lie algebra $\cg$ with the standard Lie Poisson structure and  viewed as
an abelian Lie group.    In general, Poisson groups are best understood as part of a 
 Manin triple. A Manin triple is a triple
$(\cg, \cg_+, \cg_-)$, where $\cg$ is a Lie algebra with a non-degenerate
bi-linear form $< , >$, 
$\cg_+, \cg_-$ are Lie subalgebras of $\cg$, $\cg=\cg_+ + \cg_-$ as direct sum  of
vector spaces, and $<\cg_+,\cg_+>=<\cg_-,\cg_->=0$.  Then the corresponding Lie
group $G_+$ is Poisson and $G_-$ is its Poisson dual.  The triple $(G,G_+, G_-)$ is
called a {\it double group\/} in the literature. In this paper, we will call this
triple a {\it Manin triple group\/} to avoid confusion with the completely different
concept of a double loop group.  For example,
$(SL(n), SU(n), B_n)$ is a Manin triple group, where
$B_n$ is the subgroup of upper triangular matrices in $SL(n,C)$ with real diagonal entries and
$<x,y>=\Im(\tr(xy))$ is the non-degenerate bi-linear form. 
For the trivial Poisson structure
on $G$, the Manin triple group is $(G\sdp_{\ad} \cg^\ast, G, \cg^\ast)$, where 
$G\sdp_{\ad} \cg^\ast$ is the semi-direct product.

\ms

\ni $\bu$ \hskip 6pt {\bfs Adler-Kostant-Symes Theorem\/}
\ms

Let $\cg$ be a Lie algebra equipped with an ad-invariant, non-degenerate bilinear
form $< , >$. Suppose $\ck$ and $\cn$ are  Lie subalgebras of $\cg$ such that $\cg$ is
the direct sum of $\ck$ and $\cn$ as vector space. Then the space $\ck^\perp$
perpendicular to $\ck$ in $\cg$ with respect to $< , >$ can be identified as the dual
$\cn^*$ of $\cn$. Let $M\subset
\ck^\perp$ be a coadjoint $N$-orbit equipped with the
standard co-adjoint orbit symplectic structure. The Adler-Kostant-Symes
theorem ([AdM], [Kos]) states that if $f$ and $g$ are Ad-invariant function from $\cg$
to
$R$ then
$f\n M$ and $g\n M$ are commuting Hamiltonians. For example, let $<x,y>=\tr(xy)$
and
$sl(n,R)=\ck+ \cn$, where $\ck=so(n)$ and  $\cn$ is the subalgebra of real, trace zero, upper
triangular matrices. Then $\ck^\perp$ is the space of real, symmetric, trace zero matrices,
and the coadjoint $N$-orbit $M$ at $x_0=\sum_{i=1}^{n-1} (e_{i,i+1}+e_{i+1,i})$ is the set of
all tridiagonal matrices $z=\sum_{i=1}^n x_ie_{ii} + \sum_{i=1}^{n-1} y_i (e_{i,
i+1}+e_{i+1,i})$ such that all $y_i>0$ and $\sum_i x_i=0$. Note that $f_k(x)=\tr(x^k)$ is
Ad-invariant function on $sl(n,R)$. So the Hamiltonians $H_2=f_2\n M$, $\cdots$, $H_n=f_n\n
M$ are commuting, and the  Hamiltonian system on
$M$ corresponding to $H_2$ is the Toda lattice. 

 Adler and van Moerbeke [AdM] have shown that many finite dimensional
completely integrable systems can be obtained by applying the
Adler-Kostant-Symes theorem to suitable Lie algebras. For more examples, see
also the paper by Reyman [R].

\ms
\ni $\bu$ \hskip 6pt {\bfs Poisson actions and dressing actions\/}
\ms

An action of a Poisson group $G$ on a Poisson manifold $M$ is Poisson if the action
$G\times M\to M$ is a Poisson map. When $G$ is equipped with the trivial Poisson structure,
a $G$-action on a symplectic manifold is Poisson if and only if it is symplectic. The coadjoint
action of $G$ on $\cg^*$ is Poisson in this trivial structure. In general, if $(G,G_+,
G_-)$ is a Manin triple group such that the multiplication map from $G_-\times G_+$
to
$G$ is an isomorphism, then the action of $G_+$ on $G_-$ defined by
$g_+\ast g_-=\tilde g_+$, where
$\tilde g_+$ is obtained from the factorization 
$$g_+g_-=\tilde g_-\tilde g_+\in G_-\times G_+,$$ is Poisson.  This action of $G_+$ on
$G_-$ is  called the {\it dressing action\/}. To construct a global dressing action,
every element in $G$ must be factored as a product $g_+g_-\in
G_+\times G_-$.  For example, in reference to the example in the previous
paragraph, the factorization of
$g\in GL(n)$ as
$g_+g_-\in U(n)\times B_n$ can be obtained by applying the Gram-Schmidt process
to the columns of $g$. In general, this factorization cannot be carried out in the
entire group $G$. 

\ms
\ni $\bu$ \hskip 6pt {\bfs Birkhoff decompositions theorems\/}
\ms

We remark here that all of the definitions and results in symplectic and Poisson geometry
mentioned above make sense in infinite dimensions. 

Two typical examples of infinite dimensional Manin triple groups
 are given by:
\ss
\item {(1)} $G=$ the loop group of smooth maps from $S^1$ to $GL(n,C)$, $G_+$ is the subgroup
of $g\in G$ such that $g$ is the boundary value of a holomorphic map in the disk $\n\l\n<1$,
and $G_-$  is the subgroup of
$g\in G$ such that $g$ is the boundary value of a holomorphic maps in $\n\l\n>1$ and $g(1)=I$.
\item {(2)} $G$ and $G_+$ are the same as in example (1), and $G_-$  is the subgroup of
$g\in G$ such that $g(S^1)\subset U(n)$ and $g(1)=I$. 
\ss

\ni The two Birkhoff decomposition theorems, which are carefully explained by
Pressley and Segal ([PrS]) state that the multiplication map from
$G_+\times G_-$ to $G$ is injective onto an open dense subset of $G$ in example
(1) and is a diffeomorphism in example (2).  We also need a third analytic
theorem on how the decomposition depends on a parameter $x\in (x_0, \infty)$
(Theorem 7.14). 

\ms
\ni $\bu$ \hskip 6pt {\bfs Soliton equations and inverse scattering\/}
\ms

Infinite dimensional completely integrable systems are defined
in terms of the existence of action-angle coordinates. 
All  interesting infinite dimensional completely integrable Hamiltonian systems seem to
be generalizations of the ``classical'' soliton equations. To set the stage,  we give a
brief, biased history of some of the work of these equations that is directly
related to our paper. It is impossible to give a full history here and we have
omitted many major developments. 

Solitons were first observed by J. Scott Russel in 1834 while riding on horseback following
the bow-wave of a barge along a narrow canal. In 1895, Korteweg and de Vries [KdV]
derived the equation
$$q_t=- q_{xxx} - 6 qq_x \eqno(KdV)$$ to model the wave propagation in a shallow channel of
water, and obtained solitary wave solutions, i.e., $q(x,t)=f(x-ct)$ and $f$ decays at
$\pm\infty$. The modern theory of soliton equations started with the famous 
numerical computation of the interaction of solitary waves of
  the KdV equation  by
Zabusky and Kruskal ([ZK]) in 1965.   In 1967,   Gardner,
Green, Kruskal, and Miura [GGKM] used a method called inverse scattering of the 
one-dimensional linear Schr\"odinger operator to solve the Cauchy problem for
rapidly decaying initial data for the KdV equation.  In 1968, Lax ([La]) introduced 
the concept of Lax-pair for KdV and wrote KdV as the condition for a pair of
commuting linear operators.  Zakharov and Faddeev ([ZF] 1971) gave a Hamiltonian
formulation of  KdV, and proved that KdV is completely integrable by finding
action-angle variables. Zakharov and Shabat ([ZS] 1972) found a Lax pair of
$2\times 2$ first order differential operators for the non-linear  Schr\"odinger
equation:
$$q_t= {i\over 2}(q_{xx} + 2 \n q\n^2 q) \eqno(NLS)$$ and solved the Cauchy problem via a
similar inverse scattering.  Wadati ([Wa] 1973) used the same inverse scattering
transform for the Modified KdV  equation
$$q_t= -q_{xxx} +  6q^2 q_x.\eqno(mKdV)$$  Ablowitz, Kaup, Newell
and Segar ([AKNS1] 1973) again used this same inverse scattering
transform for the sine-Gordon equation 
$$q_{xt}=\sin q, \eqno(SGE),$$ and  also observed ([AKNS2]) that all these equations have a Lax
pair of $2\times 2$ linear operators.  In 1973, Zakharov and Manakov ( [ZMa1], [ZMa2])
``solved'' the
$3$-wave equation 
$$(u_{ij})_t={b_i-b_j\over a_i-a_j} (u_{ij})_x + \sum_{k\not= i, j}
\left({b_k-b_j\over a_k-a_j} - {b_i-b_k\over a_i-a_k}\right) u_{ik} u_{kj}, \qquad
i\not=j.$$ using a Lax pair of $3\times 3$ first order linear operators. In 1976, 
Gelfand and Dikii [GDi] found an evolution equation on the space of
$n$-th order differential operators on the line with a Lax pair (this generalizes the
KdV equation). In 1978,  Zakharov and Mikhailov studied
 harmonic maps from $R^{1,1}$ to Lie group $G$ using a Lax pair
of $\cg$-valued first order linear operators ([ZMi1], [ZMi2]). 

The scattering theory of the $n\times n$ linear system was studied by Shabat [Sh], and Beals
and Coifman [BC1], [BC2]. In a series of
papers, Beals and Coifman studied the scattering and inverse scattering theory of the first
order $n\times n$ linear operator:
\refeq[di]$$\cases{{\p\psi\over \p x} = (a\l + u)\psi,&\cr \lim_{x\to -\infty} e^{-a\l
x}\psi(x,\l)= I,&\cr e^{a\l x} \psi(x,\l)\, {\rm is\, bounded\, in\,} x \,{\rm as\,} \l\to
\infty.&\cr}$$ Here $a=\diag(a_1, \cdots, a_n)$ is a constant diagonal matrix with
distinct eigenvalues
$a_1,
\cdots, a_n$, and
$u$ lies in the space $\cs(R,gl_\ast(n))$ of Schwartz maps from $R$ to the space $gl_\ast(n)$ of
all
$y\in gl(n)$ with zero diagonal entries. The ``scattering data'',
$S$, of $u$ is defined in terms of the singularity of
$e^{-a\l x}\psi(x,\l)$ in $\l$. Assume $b\in gl(n)$ is a diagonal matrix. Beals and
Coifman defined an evolution equation, the $j$-th flow associated to $a,b$ on
$\cs(R,gl_\ast(n))$, such that if
$u(x,t)$ is a solution of this equation then the scattering data $S(\cdot, t)$ of $u(\cdot, t)$ is a
solution of the following linear equation:
$${\p S\over \p t}= [S, \l^j b],$$ i.e., $S(\l,t)=e^{-b\l^j t}S(\l,0)e^{b\l^j t}$. Then by the inverse
scattering transform they solved the Cauchy problem for the $j$-th flow equation globally.
When
$n=2$ and
$u\in su(2)$, the second flow defined by $a=b=\diag(i,-i)$ is the 
non-linear Schr\"odinger equation, and when
$n=3$ and $u\in su(3)$, the first flow defined by $a, b$ with $b\not=a$ is the
$3$-wave equation.  They prove that 
\refeq[dm]$$w={\rm Re\/}\, \int_{-\infty}^\infty
\tr(-\ad(a)^{-1}(v_1)v_2)dx$$ is a symplectic structure on
$\cs(R,gl_\ast(n))$ and all the $j$-th flows are commuting Hamiltonian flows.  In
1991,  Beals and Sattinger ([BS]) proved that the $j$-th flow equation is completely
integrable by finding  action-angle variables.  A good survey on recent results is
contained in the article by Beals, Deift and Zhou [BDZ]. 

\ms
\ni $\bu$ \hskip 6pt {\bfs Soliton equations and geometry\/}
\ms

Even earlier than they appeared in applied problems, soliton equations occurred in classical
differential geometry. It was known in the mid 19th century that solutions of the
sine-Gordon equation correspond to surfaces in
$R^3$ with constant Gaussian curvature $-1$. The B\"acklund transformations (cf. [Da1]) of
surfaces in $R^3$ generate families of new surfaces with $-1$ curvature from a given one, and
hence give a method of generating new solutions of the  sine-Gordon equation from a given
one by solving two compatible ordinary differential equations. Inspired by this classical
result, B\"acklund transformations have been constructed for a large class of the equations
already mentioned (cf. [Mi], [SZ1, 2], [GZ], [TU1]).  Many more equations in
differential geometry possesses B\"acklund transformations.  For example,
equations for submanifolds with constant curvature in Euclidean space ([TT], [Ten]),
Darboux orthogonal coordinate systems ([Da2]), and  harmonic maps from  $R^{1,1}$
into a Lie group ([U1]). 

Another interesting soliton equation made its appearance in differential geometry at the
beginning of the twentieth century.  Da Rios, a student of Levi-Civita, studied the free motion
of a thin vortex tube in a liquid medium in his master degree thesis ([dR]). He
modeled this motion using the evolution of curves in $R^3$ (the {\it vortex filament
equation\/} or the {\it smoke ring equation\/}):
\refeq[dl]$$\g_t=\g_x\times \g_{xx},$$ i.e., $\g$ evolves along the direction of the
binormal with curvature as speed.  The corresponding evolution of the geometric 
quantity
$q=k\exp(\int \tau \, dx)$ satisfies the non-linear Schr\"odinger equation, where
$k(\cdot, t)$ and $\tau(\cdot, t)$  are the curvature and torsion of the curve
$\g(\cdot, t)$. This is the Hasimoto transformation of the vortex filament equation
to the non-linear Schr\"odinger equation.  For an interesting historical account of
the multiple rediscoveries of the non-linear Schr\"odinger equation for vortex
tubes see an article by Ricca [Ri].

Recently, techniques developed in soliton theory have also been used successfully in several
geometric problems whose differential equations are elliptic. For example, the studies of
harmonic maps from $S^2$ to a compact Lie group by Uhlenbeck ([U1]), harmonic maps from a
torus to $S^3$ by Hitchen ([Hi1]), into a symmetric space by Burstall, Ferus, Pedit and Pinkall
([BFPP]), constant mean curvature tori in $S^3$ by Pinkall and Sterling ([PiS]), constant mean
curvature tori in $3$-dimensional space forms by Bobenko ([Bo]), and minimal tori in $S^4$ by
Ferus, Pedit, Sterling and Pinkall ([FPPS]). For a detailed and beautiful account of these
developments see the survey book by Guest [Gu].

\ms
\ni $\bu$ \hskip 6pt {\bfs Main goal of the paper\/}
\ms

The main point of our paper is to show how  B\"acklund
transformations, scattering theory, and the hierarchy of flows can be obtained in a
uniform and natural way  from ``dressing actions'' of suitable infinite dimensional
Manin triple groups.  

\ms
\ni  $\bu$ \hskip 6pt {\bfs Decay case\/}
\ms

Most of the soliton equations considered in our paper are evolution equations on the
space $\cs_{1,a}$ of  $gl(n)$-valued connections with one complex parameter:
$${d\over dx}+a\l + u.$$   Here $a\in u(n)$ is a fixed diagonal element and $u\in \cs(R,
\cu_a^\perp)$, where
$\cu_a^\perp$ is the orthogonal complement of the centralizer $\cu_a=\{y\in u(n)\n
[a,y]=0\}$ of $a$.

To help explain the basic Poisson group action for soliton equations, we study a simpler case
first: flows on the space $\cs_+$ of $gl(n)$-valued connections of the form 
$${d\over dx} + A(x,\l),$$
where $A(x,\l)=\sum_{j=0}^k \a_j(x)\l^j$ for some $k$ and $\a_j\in
\cs(R,u(n))$ for all $0\leq j\leq k$. Hence $A(x,\l)$ is rapidly
decaying in $x$ for each $\l\in C$.   Consider now the infinite dimensional Manin triple group
$(G,G_+,G_-)$, where
$G$ is the group of holomorphic maps $g$ from $\co\setminus \{\infty\}$ satisfying the
condition $g(\bar\l)^*g(\l)=I$,
$G_+$ is the subgroup of $g\in G$ that are holomorphic in $C$ and $G_-$ is the
subgroup of $g\in G$ that is holomorphic in $\co$, where $\co$ is a small
neighborhood of
$\infty$ in the Riemann sphere $C\cup \{\infty\}$.  Since $\cs_+$ can be identified as
a subspace of the dual of the Lie algebra $C(R,\cg_-)$, and is invariant under the
coadjoint action, 
$\cs_+$ is a Poisson manifold with the standard Lie-Poisson structure. Here we use
$\l$ to denote the loop variable for $G$ and $x$ for the variable $x\in R$.   The trivialization $F$
of $D = {d\over dx} +A(x,\l)$ in $\cs_+$ is the solution of 
\refeq[dg]$$\cases{F^{-1}F_x=A(x,\l),&\cr \lim_{x\to -\infty}F(x,\l)=I,&\cr}$$ and the
{\it monodromy\/} of $D$ is  $$F_\infty(\l)= \lim_{x\to \infty}F(x,\l).$$ Since
$A(x,\l)$ is decaying in $x$, it follows that the linear system \refdg{} has a unique global
solution. This identifies
$\cs_+$ as a subset of $C(R,G_+)$. Moreover,  the monodromy $F_\infty$ exists and  is an
element in $G_+$. The group
$G_-$ acts on $C(R,G_+)$ by the pointwise dressing action of $G_-$ on $G_+$, hence it induces
an action $\ast$ of $G_-$ on $\cs_+$.  In fact, given $g\in G_-$, factor
$$g(\cdot)F(x,\cdot)=\tilde F(x,\cdot)\tilde g(x,\cdot)\in G_+\times G_-,$$ then
$g\ast D= {d\over dx} + \tilde F^{-1}\tilde F_x$.   The fundamental theorem for the
decay case appears in section 4. We show that the action of
$G_-$ on $\cs_+$ is Poisson with the monodromy on the line as the
moment map. We call the flows generated by the action of $G_-$ on $\cs_+$ the 
{\it negative flows\/}. 

\ms
\ni $\bu$ \hskip 6pt {\bfs Rational loop group action\/}
\ms

For a fixed constant $a\in u(n)$, the phase space $\cs_{1,a}$ is a coadjoint orbit of
$C(R,G_-)$, and the $w$ defined by formula \refdm{}, is the  Kostant-Kirillov
symplectic form. Since $a\l+u$ does not decay in $x$,  the monodromy of the
connection
${d\over d x} + a\l + u(x)\in \cs_{1,a}$ on the line is not defined.  The ``action'' of
$G_-$ on $\cs_{1,a}$ can still be defined formally by the dressing action. In fact,
first we identify $A$ with its trivialization
$E\in C(R,G_+)$ of $A$ normalized at
$x=0$, i.e., $E$ is the solution of
$$\cases{E^{-1}E_x= a\l + u,&\cr E(0,\l)=I.&\cr}$$ Let $\tilde E(x,\cdot)$ denote
the dressing action of $g$ on $E(x,\cdot)$ for each $x\in R$. Then $$\tilde E^{-1}\tilde
E_x=a\l +\tilde u$$ for some smooth $\tilde u$. In general, $\tilde u$ does not decay
at $\pm\infty$ for $g\in G_-$. So the action of $G_-$ on $\cs_{1,a}$ does not exist. But $\tilde u$
does belong to the Schwartz class for
$g\in G_-^m$, the subgroup of rational maps in $G_-$.
Hence the subgroup $G_-^m$ does act on $\cs_{1,a}$.  Moreover, if $g\in G_-^m$ is a linear
fractional transformation then $\tilde u$ can be obtained by solving an ordinary differential
equation or by an algebraic formula in terms of $g$ and $E$.  These results are proved in
section 6. 

\ms
\ni $\bu$ \hskip 6pt {\bfs Homogeneous structure of scattering data\/}
\ms

We are motivated by results from scattering theory to choose  the group $D_-$  of
meromorphic maps $f$ from
$C\setminus R$ to $GL(n,C)$, which satisfy the following conditions:
\item {(i)} $f(\bar \l)^*f(\l)=I$,
\item {(ii)} $f$ has a smooth extension to the closure $\bar C_\pm$, 
\item {(iii)} $f$ has an asymptotic expansion at $\infty$,
\item {(iv)} $f_\pm(r)=\lim_{s\searrow} f(r\pm is)$ such that
$f_+=v_+h_+$ factors with $v_+$ unitary and $h_+$ upper triangular and $h_+-I$ 
in the Schwartz class.

\ms
Note that $G_-^m$ is a subgroup of $D_-$. But $D_-$ is not a subgroup of $G_-$ because we do
not assume $f$ is holomorphic at $\l=\infty$.  Let $D_-^c$ denote the subgroup of
$f\in D_-$ such that $f$ is holomorphic in $C\setminus R$. In section 7, we prove that $D_-$ is
diffeomorphic to $G_-^m\times D_-^c$ by translating the Birkhoff decomposition theorems for
maps from the unit circle to maps from the real line using a linear fractional transformation. 
We identify the space
$\cs_{1,a}$ as the homogeneous space $D_-/H_-$ of left cosets of $H_-$ in
$D_-$ (scattering cosets), where $H_-$ is the subgroup of $f\in D_-$ that commutes with $a$. 
 In fact, given $f\in D_-$, we use the Birkhoff decomposition 
\refeq[dk]$$f(\l)^{-1}e^{a\l x}=E(\l,x)M(\l,x)^{-1},$$ with $E\in C(R,G_+)$ and
$M\in C(R,D_-)$. Then  $E^{-1}E_x$ is of the form
$a\l+u$ for some decay map $u$, and the map sending the left coset $H_-f$ to ${d\over
dx}+E^{-1}E_x$  gives the identification of $D_-/H_-$ and $\cs_{1,a}$. Moreover, the right
action of $D_-$ on
$D_-/H_-$  induces an action of $D_-$ on $\cs_{1,a}$, which extends the action $\ast$ of the
$G_-^m$ on $\cs_{1,a}$ defined in section 6.   

\ms
\ni $\bu$ \hskip 6pt {\bfs Poisson structure of positive flows\/}
\ms

Let $H_+$ denote the subgroup of  $G_+$ generated by $\{e^p\n p$ is a polynomial 
in $\l$ which commutes with $a\}$. In section 8, we construct an action of  $H_+$ on
$D_-/H_-$ by the ``dressing action'' of $H_+$ on $D_-$. Hence it induces an action of 
$H_+$ on
$\cs_{1,a}$. If
$a$ is regular (i.e., $a$ has distinct eigenvalues) then $H_+$ is abelian, the action  of
$H_+$ on
$\cs_{1,a}$ is Hamiltonian, and the flows generated by $H_+$ are the commuting 
hierarchy of the $j$-th flows. If $a$ is singular, then
$H_+$ is a non-abelian Poisson group and $H_+$ contains a distinguished infinite
dimensional abelian subgroup generated by polynomials in $a$. Although the action
of $H_+$ is not symplectic, we prove the action of $H_+$ on $\cs_{1,a}$ is Poisson by 
constructing a moment map. Here we need to prove the difficult result that
$M_{\pm\infty}\in H_-$, where
$M_{\pm\infty}(\l)=\lim_{x\to
\pm\infty} M(x,\l)$. Then
$M_{-\infty}^{-1}M_\infty$ is a moment map for the $H_+$-action. 
 We also  show that the pull back of
the symplectic form $w$ to the space of continuous scattering cosets 
$D_-^c/(H_-\cap D_-^c)$ is non-degenerate. We believe the restriction of $w$ to
each algebraic component of the space of discrete scattering cosets
$G_-^m/(H_-\cap G_-^m)$ is also non-degenerate, and we  prove this in one case. 

\ms
\ni $\bu$ \hskip 6pt {\bfs B\"acklund transformations\/}
\ms

Since $G_-^m$ acts on the phase space $\cs_{1,a}$, it induces an action of $G_-^m$ on the space
of solutions of the $j$-th flow. In general, if $G$ acts on $M$, the induced action of
$G$ on the space of solutions of a dynamical system is not easy to write down. In section 9, we
prove that  the induced action of
$G_-^m$ on the space of solutions of the $j$-th flow on $\cs_{1,a}$ can be 
constructed again by dressing action as done in section 6. In fact, if $a\l + u$ is a
solution of the $j$-th flow and 
$g\in G_-^m$ is a linear fractional transformation, then
$g\sharp (a\l +u)$ can be obtained by solving two compatible  ordinary
differential equations.  The action of such $g$ gives the classical B\"acklund
transformation for the sine-Gordon equation. The orbit of the rational  negative
loop group $G_-^m$ through the vacuum (trivial) solution can be computed
explicitly, and is the space of pure solitons.  Using the action of $G_-^m$, we are
also able to construct periodic (breather) solutions for the harmonic map equation
and the $j$-th flow equation with $j\geq 2$. 

\ms
\ni $\bu$ \hskip 6pt {\bfs Geometric Non-linear Schr\"odinger equation\/}
\ms

In section 10, we apply soliton theory to the Schr\"odinger flow on $Gr(k,C^n)$. 
Suppose
$(M,g,J)$ is a complex Hermitian manifold. The geometric non-linear
Schr\"odinger equation (GNLS) is the following evolution equation of curves on
$M$:
$$J\phi_t = \D \phi = \K_{\phi_x} \phi_x,\eqno(GNLS) $$ where $\K$ is the
Levi-Civita connection of the metric $g$.  When
$M=S^2$, this equation is equivalent via the Hasimoto transformation to the
non-linear Schr\"odinger equation. In fact, if $\g$ evolves according  to the vortex
filament equation
\refdl{} and
$x$ is the arc length parameter, then $\g_x$ satisfies the geometric non-linear
Schr\"odinger equation on $S^2$.  When 
$M$ is the complex  Grassmannian manifold $Gr(k,C^n)$, the GNLS gives the  matrix
non-linear Schr\"odinger equation (MNLS) studied by Fordy and Kulish [FK] for maps
$q$ from
$R^2$ to the space $\cm_{k\times (n-k)}$ of $k\times (n-k)$ complex matrices:
$$q_t = {i\over 2} (q_{xx} + 2qq^*q),\eqno(MNLS)$$ where $q^*=\bar q^t$. The 
MNLS is the second flow on $\cs_{1,a}$ defined by 
$$a=\pmatrix{iI_k&0\cr 0& -iI_{n-k}\cr}.$$  This flow has a Lax pair:
$$\eqalign{&\bigg[{\p\over \p x} + \pmatrix{i&0\cr 0&-i\cr}\l + \pmatrix{0&q\cr
-q^*&0\cr}, \cr
&\qquad {\p \over \p t} + \pmatrix{i&0\cr 0&-i}\l^2 + \pmatrix{0&q\cr -q^*&0\cr}\l
+{1\over 2i}\pmatrix{qq^*& -q_x\cr -q_x& -q^*q\cr}\bigg]=0.\cr}$$
By applying soliton theory to the MNLS, we can solve the
 Cauchy problem globally with decay initial data, and obtain a Poisson action of $H_+$ on
$\cs_{1,a}$ such that the flow generated by
$a\l^2$ is the MNLS. The flow generated by
$b\l^j$ with
$b\in\cu_a$ commutes with MNLS. If $n=2$ then $H_+$ is abelian and the action is symplectic. If
$n>2$ then $H_+$ is non-abelian and the flows generated by $b\l^k$ with $b\in \cu_a$
commute with the MNLS. But the flows generated by  
$b_1\l^j$ and $b_2\l^s$ with $b_1,b_2\in
\cu_a$ do not commute if $[b_1,b_2]\not=0$. Not all these flows are described by differential
equations. The flow generated by $b\l^k$, $b\not=a$ and $k\not=1$, are mixed
integral-differential flows. 

\ms
\ni $\bu$ \hskip 6pt {\bfs Restriction of the phase space by an automorphism\/}
\ms

The phase space of the modified KdV (mKdV) equation is the following subspace of
$\cs_{1,a}$:
$$\cs'_{1,a}=\left\{{d\over dx} + \pmatrix{i&0\cr 0&-i\cr} \l + \pmatrix{0& q\cr -q&0\cr}
\biggl| \quad
q\in \cs(R,R)\right\}.$$  The third flow defined by $b=a=\diag(i,-i)$ leaves
$\cs'_{1,a}$ invariant and is the modified KdV flow.  While all the even flows
vanish  on $\cs_{1,a}'$, all odd flows leave $\cs'_{1,a}$ invariant. 
This is a special case of  restrictions given by finite order automorphisms. To explain
this in a more general context, we let $\cu$ be a semi-simple Lie algebra (not necessary a
subalgebra of $su(n)$), and let
$<,>$ denote the Killing form.  Given
$a\in \cu$, let $\cs_{1,a}$ denote the space of all connections of the form ${d\over dx}+a\l +
u$, where $u$ is a Schwartz class map from $R$ to the orthogonal complement $\cu_a^\perp$ of
the centralizer $\cu_a$ of $a$ in $\cu$.  Then $\ad(a)$ maps $\cu_a^\perp$
isomorphically onto $\cu_a^\perp$. Hence $$w(v_1,v_2)={\rm
Re\/}\,\int_{-\infty}^\infty <-\ad(a)^{-1}(v_1),v_2>dx$$ still defines a symplectic
structure on
$\cs_{1,a}$. 

Suppose $\s$ is an order
$k$ Lie algebra automorphism of $\cu$ such that there is an eigendecomposition of $\s$  $$\cu=
\cu_0+
\cdots + \cu_{k-1},$$ where
$\cu_j$ is the eigenspace with eigenvalues $ e^{2(j-1)\pi i /k}$ with $1\leq j\leq k$.
Assume
$a\in \cu_1$, and consider the following subspace of $\cs_{1,a}$:
$$\cs^\s_{1,a}(\cu)=\left\{{d\over dx} + a\l + u\,\,\biggl|\quad u\in\cu_0\cap
\cu_a^\perp.\right\}$$ Note that when $\cu=su(2)$, $\s(x)=\bar x$, and $a=\diag(i,-i)$, we have
$\cs^\s_{1,a}=\cs_{1,a}'$.

It was shown by the first author [Te2] that
there exist a sequence of symplectic structures $w_r$ such that $w_{-1}=w$ and
all positive flows are Hamiltonian with respect to $w_r$.   In section 11, we study the
restriction of the sequence $w_r$ of symplectic forms  and the hierarchy of flows
to the subspace $\cs^\s_{1,a}$. We generalize results proved in [Te2] when $\s$ is of order $2$
and a result for the generalized modified KdV equation proved by Kupershmidt and
Wilson [KW] when
$\cu=gl(n,C)$, $\s$ is the order $n$ automorphism defined by the conjugation of the
operator $c\in GL(n)$ that permutes the standard basis of $C^n$ cyclically, and
$a=\diag(1,\a, \cdot,
\a^{n-1})$  with $\a=\exp(2\pi i/n)$. In fact, we prove:
\item {(i)} If $j\not\equiv 1$ (mod $k$), then the $j$-th flow vanishes on
$\cs_{1,a}^\s$, and if
$j\equiv 1$ (mod $k$) then the $j$-th flow leaves $\cs_{1,a}^\s$ invariant.
\item {(ii)} The restriction of $w_r$ on $\cs_{1,a}^\s$ is zero if $r\not\equiv 0$ (mod
$k$), and is non-degenerate if $r\equiv 0$ (mod $k$). 
\item {(iii)} Let $J_{rk}$ denote the Poisson structure corresponding to $w_{rk}$, and
$F_{jk+1}$ the Hamiltonian for the $(jk+1)$-th flow with respect to $J_0$. Then the $(k+1)$-th
flow satisfies the Lenard relation $$u_t=J_0(\K F_{k+1}) = J_k(\K F_1).$$

We should point out that when $\cu\subset su(n)$, $\s$ must have order $2$. So the order $k$
automorphisms occur in a more general context, in situations for which the
scattering theory is considerably more difficult than the case we have discussed.
This leads us to the question of other algebraic situations. 

\ms
\ni $\bu$ \hskip 6pt {\bfs Other semi-simple Lie algebras\/}
\ms

In this paper we have proved that all rational factorizations can be carried out, and all the
formal scattering coset data yield actual geometric flows when $\cu=su(n)$. It follows that
any problem for a Lie algebra $\cu\subset su(n)$ becomes purely an algebraic subproblem.
However, many interesting equations in differential geometry arise as flows on a twisted
space
$\cs_{1,a}^\s(\cu)$, where $\cu\not\subset su(n)$.  We believe that some form of the discrete
factorization theory and construction of scattering coset can be carried out
for many real semi-simple Lie algebras. However, one normally expects a certain
number of the factorization theorems to fail off a ``big cell''. Even more
complications arise in trying to handle systems which lie properly in the full
complex group. For example, the Gelfand-Dikii hierarchy for a $k$-th order
differential operators is linked to a restriction by an order $k$ automorphism
($k$-twist) in the full $gl(k,C)$, and the scattering theory is along rays in the
directions of $k$-th roots of unity. Our formal observations about twists apply, and
can help understand pure soliton solutions, but do not address the scattering
theory difficulties. 

\ms
\ni $\bu$ \hskip 6pt {\bfs First flows and flat metrics\/}
\ms

A symmetric space $U/K$ is formed by a splitting of the Lie algebra $\cu=\ck+\cp$, where
$$[\ck, \ck]\subset\ck, \quad [\ck,\cp]\subset \cp, \quad [\cp, \cp]\subset \ck.$$
The rank of a symmetric space is the maximal number of linearly independent commuting
elements in $\cp$, i.e., the dimension of a maximal abelian subalgebra $\ct$ in $\cp$.
Choose a basis $b_1, \cdots, b_k$ of $\ct$. Then for each element
$[f]$ of the scattering coset, from our point of view (at least formally, rigorously if
$\cu\subset su(n)$), there are
$k$ commuting first flows in variables we call $x_1, \cdots, x_k$. This yields a flat connection
$${\p\over \p x_i} + b_i\l + u_i$$ of $k$ variables for each scattering coset $[f]$.  
For example,  Darboux orthogonal coordinates in $R^n$ ([Da2]),  isometric
immersions of $R^n$ into $R^{2n}$ with flat normal bundle and maximal rank ([Te2]),
equations of hydrodynamic types ([DN1], [DN2], [Dub1], [Ts]) and Frobenius manifolds
([Dub2], [Hi2]) are of this type.  In the appendix, we apply some of the soliton theory
to these examples. 

\ms
The authors would like to thank Mark Adler, Percy Deift, Gang Tian and Pierre Van
Moerbeke for many helpful discussions. We are grateful to  Dick Palais and
Gudlaugur Thorbergsson for reading a draft of this paper.  

\bs

\newsection Review of Poisson Actions.\par

In this section, we review basic definitions and theorems on Poisson Lie groups and 
Poisson actions. Two good introductions for this material are articles by Lu and
Weinstein [LW] and Semenov-Tian-Shansky [Se1].

  A {\it Poisson structure\/} on a smooth manifold $M$ is a smooth section $\pi$ of
$L(T^*M,TM)$ such that the bilinear map $$\{,\}:C^\infty(M,R)\times C^\infty(M,R)\to
C^\infty(M,R)$$ defined by
$\{f,g\}=dg(\pi(df))$ is a Lie bracket and satisfies the condition
$$\{fg,h\}=f\{g,h\}+g\{f,h\}, \quad {\rm for\, all\,} f,g,h\in C^\infty(M,R).$$ 
We will refer to either $\{,\}$ or $\pi$ as the Poisson structure on $M$. The
section $\pi$ can also be viewed as a section of $(T^*M\otimes T^*M)^*$ or a
section of $TM\otimes TM$, which will still be denoted by $\pi$. Symplectic
manifolds are well-known examples of Poisson manifolds. 

   Let $(M,\{,\}_M)$ and $(N,\{,\}_N)$ be two Poisson manifolds. A smooth map $\phi:M\to N$
is called a {\it Poisson map\/} if $\{f_1\circ \phi,f_2\circ\phi\}_M =\{f_1,f_2\}_N\circ
\phi$. {\it The product Poisson structure\/} on $M\times N$ is defined by 
$$\{f,g\}(x,y)=\{f(*,y),g(*,y)\}_M(x) +\{f(x,*),g(x,*)\}_N(y).$$
A submanifold $N$ of $M$ is a {\it Poisson submanifold\/} if there exists a Poisson
structure on $N$ such that the inclusion map $i:(N,\{,\}_N)\to (M,\{,\}_M)$ is Poisson.

The dual $\cg^*$ of a Lie algebra has a natural Lie-Poisson structure  by 
$$\pi_\ell(x,y)=\ell([x,y]), \qquad \ell\in \cg^*, x, y\in \cg=(\cg^*)^*,$$ with
coadjoint orbits as its symplectic leaves. If $\cg$ has a non-degenerate
ad-invariant form
$(,)$, then by identifying
$\cg^*$ with $\cg$ via $(,)$, the Lie-Poisson structure on $\cg$ is 
$\pi_x(y,z)=(x,[y,z])$ for all $x,y,z\in\cg$. 

\refpar[Aa] Definition. A {\it Poisson group\/} is a Lie group $G$ together with a Poisson
structure $\pi$ such that the multiplication map $m:G\times G\to G$ is a Poisson map, where
$G\times G$ is equipped with the product Poisson structure. 

\ss
Note that $\pi(e)=0$ when $\pi$ is viewed as a map from $G\to
TG\times TG$. Moreover, the dual of $d\pi_e$ is a map from $\cg^*\times \cg^* \to \cg^*$, 
which defines a Lie bracket on $\cg^*$. The corresponding simply connected Lie group
$G^*$ has a natural Poisson structure $\pi^*$ such that the dual of $d(\pi^*)_e$ is
the Lie bracket on $\cg$. We will call $(G^*,\pi^*)$ the {\it dual Poisson group\/}
of $(G,\pi)$. This pair often fits into a larger group and we call the collection of three groups a
Manin triple group. We first explain the Manin triple at the level of Lie algebras.

\refpar[Ab] Definition. A {\it Manin triple\/} is a collection of three Lie algebras $(\cg, \cg_+,
\cg_-)$ and  an ad-invariant non-degenerate bilinear form $< , >$ on $\cg$ with the
properties:
\item {(1)} $\cg_+, \cg_-$ are subalgebras of $\cg$ and $\cg=\cg_+ +\cg_-$ as direct
sum of vector spaces,
\item {(2)} $\cg_+, \cg_-$ are isotropic, i.e., $<
\cg_+,\cg_+>=<\cg_-,\cg_->=0$. 

\ss

Let $(\cg,\cg_+,\cg_-)$ be a Manin triple
with respect to $< ,  >$.  Then $\cg_+\simeq \cg_-^*$ and the infinitesimal vector field
corresponding to $x_-\in \cg_-$ for the coadjoint action of $G_-$ on $\cg_+$ is 
$$v_{x_-}(y_+)= [x_-,y_+]_+.$$ 
The Lie Poisson structure on $\cg_+$ is 
$$(\pi_+)_{x_+}(y_-)= [x_+,y_-]_+.$$
If there are corresponding Lie groups $(G,G_+,G_-)$
we call this a {\it Manin triple group\/}. If $(G,G_+,G_-)$ is a Manin triple group, then  $G_+$
and $G_-$ have natural Poisson group structures. To describe the Poisson
structures on $G_+$ and $G_-$, we first set up some notation: Given $x_\pm\in
\cg_\pm$, let 
$\ell_{x_\pm},
\tau_{x_\pm}$ denote the
 $1$-forms on $G_\mp$ defined by   
$$\eqalign{\ell_{x_-}(y_+g_+)= <x_-, y_+>, &\qquad 
\tau_{x_-}(g_+y_+)=<x_-, y_+>,\cr \ell_{x_+}(y_-g_-)= <x_+, y_->, &\qquad
\tau_{x_+}(g_-y_-)= <x_+, y_->.\cr}$$ Then the Poisson structures on $G_{\pm}$
are given explicitly:
$$\eqalign{(\pi_+)_{g_+}(\ell_{x_-}, \ell_{y_-}) &= <(g_+^{-1}x_-g_+)_+,
g_+^{-1}y_-g_+>\cr
(\pi_-)_{g_-} (\tau_{x_+}, \tau_{y_+}) &=
<(g_-x_+g_-^{-1})_-, g_-y_+g_-^{-1}>.\cr}$$ This is equivalent to
$$(\pi_+)_{g_+}(\ell_{x_-})= g_+(g_+^{-1}x_-g_+)_+, \qquad 
(\pi_-)_{g_-}(\tau_{x_+}) = (g_-x_+g_-^{-1})_-g_-,$$
 where $g_{\pm}\in G_{\pm}$, $x_\pm\in \cg_\pm$ and
$y_\pm$ denotes the projection of $y\in\cg$ onto $\cg_\pm$ with respect to the 
decomposition
$\cg=\cg_++\cg_-$. Here we identify $\cg_-$ as
$\cg_+^*$, $\cg_+$ as $\cg_-^*$ via $< , >$, and use the matrix convention
$gx=(\ell_g)_*(x)$, $gxg^{-1}=\Ad(g)(x)$, and so forth. Since $$\ell_{x_+}(g_-)=
\tau_{(g_-^{-1}x_+g_-)_+}(g_-),$$ we have 
$$\eqalign{&(\pi_-)_{g_-}(\ell_{x_+},\ell_{y_+})\cr &= 
<(g_-(g_-^{-1}x_+g_-)_+ g_-^{-1})_-, \,\,g_-(g_-^{-1}y_+g_-)_+g_-^{-1}>\cr 
&=<(g_-(g_-^{-1}x_+g_ - -  (g_-^{-1}x_+g_-)_-)g_-^{-1})_-, \,\,
g_-(g_-^{-1}y_+g_-)_+g_-^{-1}>\cr 
&= -<g_-(g_-^{-1}x_+g_-)_-g_-^{-1}, \,\, g_-(g_-^{-1}y_+g_-)_+g_-^{-1}>\cr
&= -<(g_-^{-1}x_+g_-)_-, \,\, (g_-^{-1}y_+g_-)_+> \cr &= -<(g_-^{-1}x_+g_-)_-, \,\,
g_-^{-1}y_+g_->.\cr}$$ 
Hence $(G_+,\pi_+)$ is the dual Poisson group of $(G_-,\pi_-)$.  Conversely, if $K$ is
a Poisson group and $K^*$ is its dual Poisson group, then there exist   an
Ad-invariant form $<,>$ and a Lie bracket on $\cg=\ck+\ck^*$ such that
$(\cg,\ck,\ck^*)$ is a Manin triple.  Hence there is a bijective
correspondence between the Manin triples and simply connected Poisson groups. 
The Manin triple group
$(G,G_+,G_-)$ is called a {\it double group\/} in the literature. In some cases, 
multiplication in $G$ can not be globally defined. In this case, we call $(G,G_+,G_-)$ a
{\it local Manin triple group\/}. 

\ms

\refpar[fg] Example. Let $G=SL(n,C)$, $G_+=SU(n)$, $G_-$ the subgroup of upper triangular
matrices with real diagonal entries, and $<x,y>=\Im(\tr(xy))$ the non-degenerate bi-invariant
form on $\cg$.   Then
$(G,G_+, G_-)$ is a Manin triple group, and the multiplication map $G_+\times G_-\to G$ and
$G_-\times G_+\to G$ are isomorphisms. The decomposition of $g\in SL(n,C)$ as $g=g_+g_-\in
G_+\times G_-$ and
$g=h_-h_+\in G_-\times G_+$ are obtained by applying the Gram-Schmidt process to the
columns and rows of $g$ respectively. 

\ms 

\refpar[Af] Examples.
The type of Poisson groups we need in this paper are generally credited to Cherednik 
([Ch]). Let $\W_+$ and $\W_-$ be two domains of  $S^2=C\cup
\{\infty\}$  such that $S^2=\W_+\cup \W_-$ and both $\W_+$ and
$\W_-$ are invariant under complex conjugation.  Let $\co= \W_+\cap \W_-$.   A
map $g:\co\to SL(n,C)$ is called $su(n)${\it-holomorphic\/} if
$g$ is holomorphic and satisfies the reality condition 
$g(\bar\l)^*g(\l)=I$ for all $\l\in \co$. Let 
$$G=\{g:\co\to GL(n,C)\n g\,\, {\rm is\,\,} su(n) {\rm -\,\, holomorphic\/}\}.$$ Now
we fix a normalization point
$\l_o\in C\cup\{\infty\}$.
 If $\l_o\in \W_+$, define
$$\eqalign{
G_+&=\{g\in G\n g\,{\rm extends\,\,} su(n){\rm - holomorphically\,\, to \/}\,\,
\W_+\,\, g(\l_o)=I\},\cr  G_-&=\{g\in G\n g  \,{\rm extends\,}\, su(n){\rm -
holomorphically\,\, to
\/}\,\,\W_-\}.\cr }$$
Similarly, if $\l_o\in \W_-$, we define 
$$\eqalign{G_+&=\{g\in G\n g\,{\rm extends\,\,} su(n){\rm - holomorphically\,\, to
\/}\,\,
\W_+\},\cr  G_-&=\{g\in G\n g\,{\rm extends\,\,} su(n){\rm - holomorphically\,\, to \/}\,\,
\W_-\,\, g(\l_o)=I\}.\cr }$$
The normalization point $\l_o$ determines an Ad-invariant bilinear form
$< , >$ on $\cg=\cg_++\cg_-$ such that $(\cg,\cg_+,\cg_-)$ is a Manin triple. In fact,
$$<u,v> =\cases{{1\over 2\pi i}\oint_{\g} {\tr(u(\l)v(\l))\over
(\l-\l_o)^2}d\l, & if $\l_o\in C$,\cr
{1\over 2\pi i}\oint_{\g} \tr(u(\l)v(\l))d\l, & if $\l_o=\infty$,\cr}$$
where $\g= \partial \co$. Note that if $u(\l)=\sum_k u_k(\l-\l_o)^k$ and 
$ v=\sum_k v_k(\l-\l_o)^k$, then
$$< u,v >_{\l_o}=\cases{\sum_k \tr(u_kv_{-k+1}), & if $\l_o\in C$, \cr \sum_k
\tr(u_kv_{-k-1}), & if $\l_o=\infty$.\cr}$$ 
The main examples we use in this paper are: 
\item {(i)} $\W_+= C$, $\W_-=\co_\infty$ a  neighborhood of $\infty$,
\item {(ii)} $\W_+= C\setminus \{0\}$, $\W_-=\co_0\cup\co_\infty$.

\ni It follows from the Birkhoff Decomposition Theorem (cf. p. 120 Theorem 8.1.2
in the book by Pressley and Segal ([PrS])  that the multiplication map
$G_+\times G_-\to G$ for example (i) is injective and maps onto an open dense
subset of $G$. McIntosh shows that the multiplication map for example (ii)  is a
diffeomorphism [Mc].

\ss

Now suppose $(G,G_+,G_-)$ is a Manin triple group, and the multiplication  map
$G_+\times G_-\to G$ is a diffeomorphism. Then given
$g_\pm\in G_\pm$, we decompose $$g_+g_-=f_-f_+\in G_-G_+, \qquad g_-g_+=h_+h_-\in
G_+G_-.$$  Define 
$$g_+\#g_-=f_-,\qquad g_-\#g_+= h_+.$$
Then $\#$ defines the dressing  action of $G_+$ on $G_-$ on the left, and the
dressing action of
$G_-$ on
$G_+$ on the left respectively.   Let
$x_-\in
\cg_-$,  and $\tilde x_-$ denote the infinitesimal vector field
of the action of $G_-$ on $G_+$. Then
$$\tilde x_-(g_+)= g_+(g_+^{-1}x_-g_+)_+, \qquad \tilde x_+(g_-)= g_-(g_-^{-1}x_+g_-)_-.$$
There are clearly also corresponding dressing action of $G_-$ 
on $G_+$ and $G_+$ on $G_-$ on the right.   

Since the image of the multiplication map is an open dense subset for Example \refAf{} (i)
and the whole group $G$ for Example \refAf{} (ii), the dressing actions for the
corresponding Manin triple groups are local and global respectively. However, the
Lie algebra actions are  defined for all elements in both cases. 

\refpar[Ac] Definition. An action of a Poisson group $G$ on a Poisson manifold $P$
is {\it Poisson\/} if the action $G\times P\to P$ is a Poisson map. 

\ss
It is clear that if the
$G$-action on $P$ is Poisson, $M$ is a Poisson
submanifold of $P$, and $M$ is invariant under $G$, then the $G$-action on $M$ is also
Poisson. Here one must be careful as the requirement that $M\subset P$ is Poisson
is quite restrictive. 
\ss

A {\it symplectic structure\/} on $P$ is a Poisson structure $\pi$ such that
$\pi_x:TP^*_x\to TP_x$ is injective for all $x\in P$. This definition agrees with the standard
one when $P$ is finite dimensional, and is the definition of a weak symplectic structure
defined in the lecture notes of Chernoff and Marsden [CM] when $P$ is of infinite
dimension. For simplicity of notation, we still call such structure a symplectic
structure.   A
$G$-action on
$P$ is called symplectic if $g_*(\pi)=\pi$ for all
$g\in G$.  If $G$ is equipped with the trivial Poisson structure ($\pi_G=0$), then an
action of $G$ on a symplectic manifold
$P$ is Poisson if and only if it is symplectic. However, in general these two notions
of actions are different on symplectic manifolds.  

A {\it moment map\/} of a symplectic action of $G$ on a symplectic manifold $P$ is a
$G$-equivariant map $\mu:P\to \cg^*$ such that $\pi_P(df_\xi)$ is the infinitesimal vector field
$\tilde \xi$ associated to
$\xi$, where $f_\xi$ is the function on $P$ defined by $f_\xi(x)=\mu(x)(\xi)$. When
 the action is Poisson, we can not expect to define a Poisson map $\mu:P\to
\cg^*$. The following theorem gives a natural generalization of
moment map for Poisson actions. 

\refclaim[Ad] Theorem ([Lu]). Suppose the Poisson group $(G,\pi)$ acts on the
Poisson manifold $(P,\pi_P)$, and there exists a
$G$-equivariant Poisson map $$m:(P,\pi_P)\to (G^*,\pi^*)$$ such that 
$$\pi_P(((dm)m^{-1})(\xi))= \tilde \xi, \quad \forall\,\, \xi\in\cg,$$ where $\tilde\xi$ is the
infinitesimal vector field on $P$ associated to $\xi$ and $(G^*,\pi^*)$ is the dual Poisson
group of $(G,\pi)$. Then the action of
$(G,\pi)$ on $(P,\pi_P)$ is Poisson.

\refpar[fi] Definition. A {\it moment map\/} for a Poisson action of a Poisson group $G$ on a
Poisson manifold $P$ is a map $m:P\to G^*$ which satisfies the assumptions in the above
theorem. 

\refpar[Ae] Example.  Suppose $(G,G_+,G_-)$ is a Manin triple group, and the
multiplication maps $G_+\times G_-\to G$ and $G_-\times G_+\to G$ are
diffeomorphisms. Then the dressing action of
$(G_-,\pi_-)$ on $(G_+,\pi_+)$ is Poisson and the identity map 
$\id:G_+\to G_-^*=G_+$ is a moment map. To see this, note first that
the identity map is Poisson and equivariant. So by Theorem \refAd{} it suffices to check
$$(\pi_+)_{g_+}(dg_+g_+^{-1}(x_-))=(\pi_+)_{g_+}(\ell_{x_-})=g_+(g_+^{-1}x_-g_+)_+=\tilde
x_-(g_+).$$ Similarly, the dressing action of $(G_+,\pi_+)$ on
$(G_-, \pi_-)$ is Poisson.

\bs

\newsection Negative flows in the decay case.\par

	    Our starting point is the Manin triple $(\cg,\cg_+,\cg_-,< , >)$ of Cherednik
type (Example \refAf{} (i)) with  $\W_+= C$, $\W_-=\co_\infty$ and 
$$<u,v> = {1\over 2\pi i}\oint_{ \g_\infty} \tr(u(\l),v(\l)) d\l,$$ where
$\g_\infty=\partial \co_\infty$ is a contour around $\infty$. 
  The basic geometric object is a $\cg_+$-valued connection on the real line $R$ of the form
$$D= {d\over dx} +A(x,\l)=  {d\over dx}+ \a_k(x) \l^k + \a_{k-1}(x)\l^{k-1} + \cdots + \a_0(x).$$
From the analytic point of view there are three distinct theories which have
very different algebraic structures:
\ss

\item {(1)} Asymptotically constant cases --- the leading term $\a_k$ is a constant $a\in
gl(n,C)$ and 
$\a_j(x)$ decays in $x$ for $0\leq j<k$. 

\item {(2)} Decay case --- $\a_j(x)$ decays in $x$ for all $0\leq j\leq k$.

\item {(3)} Periodic case --- $\a_j(x)$ is periodic in $x$ for all $j$. 

\ss\ni
  Most of the classical scattering theory deals with the asymptotically constant
case, which is the case we discuss in most of the paper. For the periodic case we refer the
readers to papers by Krichever [Kr1], [Kr2].  We start with the decay case, as a warm-up
for the asymptotically constant case. 

Fix an element $\a_k\in L^1(R)$.  An important example
would be
$\a_k(x)= \rho(x) a$ for $\rho\in L^1(R)$ and $a\in gl(n,C)$. If $\rho= dy/dx$, then
we can rewrite the connection in $y$ as 
$${dy\over dx} \left({d\over dy} + a\l^k\right)= {d\over dx} + \rho(x) a\l^k.$$
Hence the decay case is in reality the case of a ``finite interval''. However, we
use the parametrization of the infinite interval to demonstrate structural
relationships with the asymptotically constant case. 

    Let $C(R,G_\pm)$ be a linear subspace of maps from $R$ to $G_\pm$, that has a
formal Lie group structure with Lie algebra $C(R,\cg_\pm)$, where $C(R)$ consists
of functions which decay at least as fast as those in $L^1(R)$. Identify a map
$A\in C(R, \cg_+)$ with an element in $C(R,\cg_-)^*$ via the pairing
$$<< A,T>> = \int_{-\infty}^\infty < A,T >\, dx.$$
(Note that if $Adx$ is thought as a one form, then the above formulation is coordinate
invariant). 
 Let $\cs$ be a subset of $C(R,\cg_+)$ that is invariant
under the coadjoint action of $C(R,G_-)$. The infinitesimal vector fields for the
coadjoint action of $C(R,G_-)$ on $\cs$ are
$$v_T(A)(x)= [A(x),T_-(x)]_+,$$
where $+$ indicates the orthogonal projection from $\cg$ onto $\cg_+$. The Poisson
structure on $\cs$ is given by
$$\pi_A(T_-)=-[A,T_-]_+,$$
and gives rise to a symplectic structure on the coadjoint orbits of $C(R,G_-)$ on
$\cs$. 

The coadjoint orbit of $\a(x)\l^k$ under $C(R,G_-)$ is clearly contained in the set of
polynomials of degree $k$ of the form
$$A= \a(x)\l^k+ \a_{k-1}(x) \l^{k-1} + \cdots + \a_1(x)\l + \a_0(x)$$
with the condition that $\a_{k-1}(x)$ is of the form $[\a(x), v(x)]$ for some $v$. For many
choice of $\a$, this will be the only constraint. The
vector field $v_T$ for $T=\sum_{k=1}^\infty T_j(x)\l^{-j}$ is
$$v_T(A)=[A,T]_+ = \sum_{j=0}^{k-1}\left(\sum_{i=j+1}^{k}
[\a_i(x),T_{i-j}(x)]\right)\l^j,$$
where $\a_k=\a$. 

\ss

The negative flows in the decay case can be easily described. Let $\cp(R,\cg_+)$
denote the Lie
 algebra of maps $A:R\to\cg_+$ such that $A(x)(\l)$ is a polynomial in $\l$ and decay
in $x$. Let $\cp_k$ denote the set of all $A\in \cp(R,\cg_-)$ of degree $k$, and
$\cp_{k,\a}$ the set of all $A\in \cp(R,\cg_-)$ whose leading term is $\a\l^k$.   Then
$\cp(R,\cg_+)$, $\cp_k$ and
$\cp_{k,\a}$ are invariant under the coadjoint action of $C(R,G_-)$, and 
$$\pi_A(T_-)=-[A,T_-]_+$$ gives the Poisson structure.

\refpar[hy] Definition.   The {\it trivialization of  $A=\sum_{j=0}^k \a_j(x)\l^j$
normalized at $x=-\infty$\/} is the solution  
$F(A)\in C(R,G_+)$ of $$F^{-1}F_x=A, \qquad \lim_{x\to -\infty}
F(x,\l)=I.$$  

\ss
Given $b\in su(n)$ and $A=\sum_{j=0}^k \a_j(x)\l^j$, then $F(A)^{-1}(x)bF(A)(x)\in \cg_+$. 
Write the expansion of $F(A)^{-1}bF(A)$ at
$\l=0$ to get
\refeq[bk]$$F(A)^{-1}bF(A)= \b_0+\b_1\l+\b_2\l^2 +\cdots.$$  The
$\b_j$'s can be computed explicitly from $A$. Since 
\refeq[ff]$$(F^{-1}bF)_x +
[A, F^{-1}bF]=0,$$ we can compare coefficients of $\l^j$ in equation \refff{} to get
$$\cases{(\b_0)_x+[\a_0,\b_0]=0,&\cr (\b_j)_x+ [\a_0,\b_j] +\sum_{i=1}^{{\rm min\/}
\{j,k\}} [\a_i,\b_{j-i}]=0.&\cr}$$
The $\b_j$'s can be solved explicitly from $\a_0, \cdots, \a_k$ as follows: Let $g:R\to
GL(n,C)$ be the solution to 
$$\cases{ g^{-1}g_x=\a_0 &\cr  \lim_{x\to -\infty} g(x)=I.&\cr}$$
Then 
\refeq[bp]$$\cases{\b_0= g^{-1}bg,&\cr 
\b_j(x) = -g^{-1}(x)   \sum_{i=1}^{\min \{j,k\}}\left(\int_{-\infty}^x 
g(y)[\a_i(y),\b_{j-i}(y)]g^{-1}(y)dy\right) g(x).&\cr}$$ 
 Hence we have obtained a family of integral
equations to describe the $\b_j$'s.

The Lax pair for this system is written
\refeq[hu]$$\left[{\p\over \p x} + A,\,\, {\p \over \p t} + (F^{-1}b\l^{-m}F)_-\right]=0.$$ 
 It follows from the definition of $\b$'s that the coefficient of $\l^j$ with $j<0$ in the left
hand side of equation \refhu{} is automatically zero. Setting the coefficients of
$\l^j$ ($j\geq 0$) in equation \refhu{} to zero gives a system of equations
describing a flow on
$\cp_k$:
\refeq[ga]$$\cases{{d\a_k\over dt}=0,&\cr {d\a_j\over dt}= \sum_{i=j+1}^{\min \{k,
m+j\}} [\a_i,\b_{m+j-i}].\cr}$$
We cal this flow the $-m$-{\it flow on
$\cp_{k,\a}$ defined by $b$}.  Equation \refhu{} also gives
$$\eqalign{A_t &= ((\l^{-m} F^{-1}bF)_x)_- + [A, (\l^{-m} F^{-1}bF)_-]\cr
&=[(\l^{-m} F^{-1}bF), A]_- + [A, (\l^{-m} F^{-1}bF)_-] \cr &= [A,
(F^{-1}b\l^{-m}F)_-]_+.\cr}$$ 
So the $-m$-th flow can also be written as  
\refeq[Ah]$$A_t= [A,(F^{-1}b\l^{-m}F)_-]_+.$$
Since the vector field 
$$\xi_{b,j}(A)= [A,(F^{-1}b\l^{-m}F)_-]_+$$ is bounded in $L^1$, it is not difficult to see that the
$-m$-flow is global. We will prove these flows generates a natural Poisson group action on
$\cp_{k,\a}$ in the next section. 

\ms 

For our basic model, $k=1$, we have
$$\eqalign{&A = \a(x) \l + u(x), \cr & v_T(u)=[\a(x),T_1(x)],\cr 
&\{T,V\}_A = \int_{-\infty}^\infty \tr(\a(x)[T_1(x),V_1(x)])dx,\cr}$$
where $T=\sum_{j=1}^\infty T_j\l^{-j}$ and $V=\sum_{j=1}^\infty V_j\l^{-j}$. 
This gives our next proposition.

\refclaim[hw] Proposition.  
The $-m$-th flow on $\cp_{1,\a}$ defined by
$b$ is 
$$u_t=[\a, \b_{m-1}],$$ where $\b_j$ is defined inductively by
$$\cases{\b_0 = g^{-1}b g,&\cr 
\b_j(x)= -g^{-1}(x)\left(\int_{-\infty}^x g[\a,\b_{j-1}]g^{-1} dy\right) g(x),&\cr}$$
and 
$g$ is the solution to $g^{-1}g_x= u$ and $\lim_{x\to -\infty} g(x) = I$. 

A simple change of gauge (cf.
[Te2]) implies that the $-1$-flow describes the geometric equation for harmonic maps 
from $R^{1,1}$ into $U(n)$ in characteristic coordinates: 

\refclaim[ce] Proposition. Fix a smooth $L^1$-map $\a:R\to u(n)$ and $b\in u(n)$.
Suppose $u(x,t)$ is a solution of the $-1$-flow equation on $\cp_{1,\a}$ defined by $b$:
\refeq[dd]$$u_t=[\a, g^{-1}bg],\qquad {\rm where\,\,}
g^{-1}g_x=u,\,\, \lim_{x\to -\infty} g(x)= I.$$  Then there exists a unique solution 
$E(x,t,\l)$  for  $$\cases{E^{-1}E_x= \a\l+u, &\cr E^{-1}E_t= \l^{-1}g^{-1}bg,&\cr
E_\l(0,0)=I.&\cr}$$ Set $s(x,t)=E(x,t,-1)E(x,t,1)^{-1}$. Then $s: R^{1,1}\to U(n)$
is harmonic,  
$(s^{-1}s_x)(x,t)$ is conjugate to $\a(x)$, and $(s^{-1}s_t)(x,t)$ is conjugate to $b$
for all
$t\in R$.

Harmonic maps into a symmetric space are obtained by restriction ([Te2]). This is discussed
in section 9. Also, a more elaborate choice of Cherednik splittings
allows more complicated examples like the harmonic map equation in space-time 
(laboratory) coordinates. 
\bs

\newsection Poisson structure for negative flows (decay case).\par

The dressing action defines a local action of $G_-$ on
$\cp(R,\cg_+)$ which is Poisson and generates the
negative flows.  The notation is the same as in section 3. 

\refclaim[Be] Theorem. For $A\in \cp(R,\cg_+)$, let
$F(A):R\to G_+$ denote the trivialization of $A$ normalized at $x=-\infty$. Given $g_-\in
G_-$, let  $\tilde F(x) = g_-\sharp (F(A)(x))$, where $\sharp$ denotes the dressing action
of $G_-$ at $G_+$ for each $x\in R$. Define
$$g_-\ast A = \tilde F^{-1}\tilde F_x.$$   Then $g_-\ast A$ defines a local action of 
$G_-$ on $\cp(R,\cg_+)$. Moreover,  the infinitesimal vector field $\tilde \xi_-$ associated to
 $\xi_-\in\cg_-$ for this action is
\refeq[bj]$$\tilde \xi_-(A)= -[A, (F^{-1}(A)\xi_-F(A))_-]_+.$$

\proof  It is clear that $(g_-\ast A)$ defines a local action of $G_-$ on $\cc(R,\cg_+)$. Now
we compute the infinitesimal vector field $\tilde \xi_-$ on $\cc(R,\cg_+)$. 
Write $g_-F=\tilde F f_-$, and let $\d$ denote the tangent variation. Then 
$(\d g_-) F= \d \tilde F + F\d f_-$, which implies that 
$$F^{-1}(\d g_-)F= F^{-1}\d \tilde F + \d f_-.$$
If $\xi_-=\d g_-$, then we have 
\refeq[Bf]$$\d f_-= (F^{-1}\xi_- F)_-,\qquad F^{-1}\d \tilde F = (F^{-1}\xi_- F)_+.$$
Since $g_-* A= \tilde F^{-1}\tilde F_x$, we obtain
$$\eqalign{\tilde \xi_-(A)&= -F^{-1}(\d \tilde F) F^{-1}F_x + F^{-1} (\d \tilde F)_x\cr
&= -(F^{-1}\xi_-F)_+A + F^{-1}(F(F^{-1}\xi_-F)_+)_x\cr
&=-(F^{-1}\xi_-F)_+A  + A(F^{-1}\xi_-F)_+ + ((F^{-1}\xi_-F)_x)_+\cr
&= [A, (F^{-1}\xi_-F)_+] + [F^{-1}\xi_-F, A]_+\cr
&= -[A,(F^{-1}\xi_-F)_-]_+.\cr}$$

Since $x\mapsto A(x)(\l)$ is in $L^1(R)\cap C^\infty(R)$ and $x\mapsto F(x,\l)$
is bounded for all $\l$, we have
$\tilde \xi_-$ is tangent to $\cp(R,\cg_+)$. \qed 

\refclaim[cg] Corollary. The local action of $G_-$ on $\cp(R,\cg_+)$ leaves
$\cp_{k,\a}$ invariant, and the flow generated by $\xi_-=-b\l^{-m}$  is the
$-m$-flow on $\cp_{k,\a}$ defined by $b$. 

\refclaim[Bb] Theorem. The local action of $G_-$ on $\cp(R,\cg_+)$ is Poisson. The infinitesimal
vector field corresponding to $\xi_-$ is $\tilde \xi_-(A)=
-[A,(F^{-1}\xi_-F)_-]_+$, where $F$ is the trivialization of $A$ normalized at
$x=-\infty$.  In fact,  the map
$\phi:\cp(R,\cg_+)\to G_+=G_-^*$ defined by
$\phi(A)=
\lim_{x\to\infty}  F(A)(x)$  is a moment map for this action. 

To prove the theorem,  we first need a lemma:

\refclaim[Bc] Lemma. $d\phi_A(B)=\bigl(\int_{-\infty}^\infty
F(A)BF(A)^{-1}dx\bigr)\phi(A)$.

\proof Let $F$ denote $F(A)$, and $\d F= dF_A(B)$. Taking the
variation of the equation $F^{-1}F_x=A$, we get $(F^{-1}\d F)_x +[A,F^{-1}\d F]=B$. This
implies that
$$F^{-1}\d F =F(A)^{-1}\left(\int_{-\infty}^x
F(A)(y,\l)B(y,\l)F^{-1}(A)(y,\l)dy\right)F(A).$$ 
Then the lemma follows from taking the limit as $x\to \infty$. \qed

\refpar[fj] Proof of Theorem \refBb{}.
It suffices to prove that $\phi$ satisfies the
assumption in Theorem \refAd{}. First we prove that $\phi$ is $G_-$-equivariant. Taking the
limit of $g_-F=\tilde F f_-$ as $x\to -\infty$, we get 
$$\lim_{x\to -\infty}\tilde F(\l,x)= I, \qquad \lim_{x\to -\infty}
f_-(\l,x)=g_-(\l).$$ So $F(g_-*A)=\tilde F$ and 
$$\phi(g_-*A)=\lim_{x\to \infty}\tilde F= g_-\phi(A)(\lim_{x\to \infty}
f_-)^{-1}=g_-\#
\phi(A).$$ This proves that $\phi$ is $G_-$-equivariant.

Given $\xi_-\in\cg_-$ and $B\in \cp(R, \cg_+)$, using Lemma \refBc{} we get
$$\eqalign{<< d\phi_A(B)(\phi(A))^{-1},\xi_->>&= << F(A)BF(A)^{-1}, \xi_- >>\cr 
& = << B, F(A)^{-1}\xi_-F(A) >>\cr& = << B, (F^{-1}\xi_-F)_- >>.\cr}$$ So
$(\Pi_+)_A(d\phi_A\phi(A)^{-1},\xi_-)= \tilde \xi_-(A)$.

It remains to prove that $\phi$ is a Poisson map. Given $\xi_-,\eta_-\in \cg_-$, let
$g_+= \phi(A)$, and $\ell_i$ the linear functional on $T(G_+)_{g_+}$ defined by
$$\ell_1(z_+g_+)=<\xi_-,z_+>, \qquad \ell_2(z_+g_+)=<\eta_-,z_+>.$$
It follows form Lemma \refBc{} that 
$$\ell_1\circ d\phi_A(B)=< B, (F^{-1}\xi_-F)_- >, \qquad 
\ell_2\circ d\phi_A(B)=< B, (F^{-1}\eta_-F)_->.$$
But $(F^{-1}\xi_-F)_x+[A,F^{-1}\xi_-F]=0$. So we get
$$\eqalign{&\Pi_+(\ell_1\circ d\phi_A,\ell_2\circ d\phi_A)\cr &= -<< [A,
(F^{-1}\xi_-F)_-], (F^{-1}\eta_-F)_- >>\cr &= -<< [A,
(F^{-1}\xi_-F)-(F^{-1}\xi_-F)_+],(F^{-1}\eta_-F)_->> \cr 
&= << (F^{-1}\xi_-F)_x+[A,(F^{-1}\xi_-F)_+], (F^{-1}\eta_-F)_->>\cr
&= << (F^{-1}\xi_-F)_x,(F^{-1}\eta_-F)_->> +  << [A, (F^{-1}\xi_-F)_+],
(F^{-1}\eta_-F)_->>\cr
&= < F^{-1}\xi_-F,(F^{-1}\eta_-F)_-> \n_{x=-\infty}^{x=\infty} -<< (F^{-1}\xi_-F),
((F^{-1}\eta_-F)_-)_x>>\cr &\qquad+ << [A, (F^{-1}\xi_-F)_+],
(F^{-1}\eta_-F)_->>.\cr}$$ The first term is equal to
$$\eqalign{& < g_+^{-1}\xi_-g_+, (g_+^{-1}\eta_-g_+)_-> - <\xi_-,\eta_-> =
< g_+^{-1}\xi_-g_+, (g_+^{-1}\eta_-g_+)_- >\cr &= (\pi_+)_{g_+}(\xi_-g_+,\eta_-g_+)
=(\pi_+)_{g_+}(\ell_1,\ell_2),\cr}$$ where $<\xi_-,\eta_->=0$ because $\cg_-$ is
isotropic with respect to $(,)$. The second term  is $$\eqalign{<<
(F^{-1}\xi_-F),((F^{-1}\eta_-F)_-)_x >>  &= <<
(F^{-1}\xi_-F)_+,(F^{-1}\eta_-F)_x>>\cr &=
<< (F^{-1}\xi_-F)_+,-[A,F^{-1}\eta_-F]>>\cr &=<< [A,
(F^{-1}\xi_-F)_+],F^{-1}\eta_-F>>
\cr &=<< [A, (F^{-1}\xi_-F)_+],(F^{-1}\eta_-F)_->>,\cr}$$ which cancels the third
term. This proves that $\phi$ is Poisson. Since $\phi$ satisfies all assumptions of
Theorem
\refBb{}, the action of $G_-$ on $\cl_+$ is Poisson and $\phi$ is a moment map.  \qed

\bs

\newsection Positive flows in the asymptotically constant case.\par

In this section, we will use the same Manin triple as in section 3, and describe flows in the
asymptotically constant case. We restrict our discussion to the simplest cases.

  Fix $a\in su(n)$, and set
$$\eqalign{U_a&=\{g\in SU(n)\n ga=ag\},\cr
\cu_a&=\{y\in su(n)\n [a,y]=0\}, \cr
\cu_a^\perp&=\{z\in su(n)\n <z,\cu_a>=0\}.\cr}$$  Given a vector  space $V$, we let
$\cs(R,V)$ denote the space of all maps from $R$ to $V$ that are in the Schwartz class. 
Let $\cs_{1,a}$ denote the
space of all maps $A:R\to \cg_+$ such that $A(x)(\l)= a\l+u(x)$ with $u\in
\cs(R,\cu_a^\perp)$. The basic symplectic structure on $\cs_{1,a}$ is
similar to what we have described already for the decay case. 
However, the structure of the natural flows is different because we may not
normalize at $x=-\infty$. Integration as described in the negative flows will tend
to destroy the decay condition. The  $-1$-flow does in some sense exist:
$\cs_{1,a}$
\refeq[gj]$$\cases{u_t= [a, g^{-1}bg], &\cr g_x=g u, &\cr\lim_{x\to -\infty} g=
I.&\cr}$$ However, the right-hand boundary at $\infty$ will not be under control and
the symplectic structure does not make coherent sense. 

	Rather than identify $A$ with the trivialization $F$
normalized at $x=-\infty$, we use two different trivializations. For the purposes of
constructing B\"acklund transformations, we identify $A$ with the trivialization $E$
normalized at $x=0$, i.e., 
$$E^{-1}E_x= a\l + u,\qquad E(0,\l)=I.$$  When we describe the Poisson structure of the positive
flows we use $M(x,\l)$, where
$$(e^{a\l x}M)^{-1}(e^{a\l x}M)_x= A, \qquad \lim_{x\to -\infty} M(x,\l)= I.$$
Since both $E$ and $e^{a\l x} M$ solve the same linear equation, there exists $f(\l)$ such
that 
$$f(\l)E(x,\l) =  e^{a\l x} M(x,\l).$$ Note that 
$f(\l)=M(0,\l)$ contains all the spectral information. The general condition 
is that $f$ is not holomorphic at
$\l=\infty$, but that both $f(\l)$ and
$M(x,\l)$ have asymptotic expansion at $\l=\infty$. This is known to be the case
in scattering theory, and we need our theory to mesh with this analysis. 

The positive flows for the asymptotically constant case are defined in a similar
fashion as the negative flows for the decay case with the restriction that the
generators commute with $a$. The hierarchy of flows is now mixed ordinary
differential and integral equations. Let
$A= a\l+ u$ with $u\in
\cs(R,\cu_a^\perp)$, and
$M$ as above.   Fix
$b\in u(n)$ such that $[a,b]=0$. Then  $M^{-1}bM$ has an asymptotic expansion at $\l=
\infty$ (cf. [BC1,2]): 
$$M^{-1}bM \,\, \sim\,\,  Q_{b,0}+Q_{b,1}\l^{-1} + Q_{b,2}\l^{-2} + \cdots , \qquad
Q_{b,0}=b.$$ Since  $$M^{-1}bM=E^{-1}f^{-1} e^{a\l x}b e^{-a\l x}fE = E^{-1}f^{-1}bf E,$$
we get
$(M^{-1}bM)_x+[a\l+u,M^{-1}bM]=0$. So we have
\refeq[db]$$(Q_{b,i})_x + [u,Q_{b,i}] + [a, Q_{b,i+1}]=0.$$
This defines $Q_{b,i}$'s recursively.  

An element $a\in u(n)$ is {\it regular\/} if $a$ has distinct
eigenvalues. Otherwise,
$a$ is {\it singular\/}. If 
$a$ is regular, then it is known that 
$Q_{b,i}$'s are polynomial differential operators in $u$ (cf [Sa]). But when $a$ is
singular, the $Q_{b,i}$'s are integral-differential operators in $u$. To be more
precise, we decompose 
$$Q_{b,j}=P_{b,j} + T_{b,j}\in \cu_a + \cu_a^\perp.$$
Using equation \refdb{}, $P$'s and $T$'s can be solved recursively. In fact,
\refeq[bq]$$\cases{P_{b,0}= 0,&\cr T_{b,0}=b,&\cr P_{b,j+1}=
-\ad(a)^{-1}((P_{b,j})_x + [u,Q_{b,j}]^\perp),&\cr
T_{b,j+1}=-\int_{-\infty}^x[u,P_{b,j+1}]^{\ct} dy,&\cr}$$ where $v^\perp$ and
$v^{\ct}$  denote the projection onto $\cu_a^\perp$ and
$\cu_a$ respectively, and 
 $-\ad(a)$ maps $\cu_a^\perp$ isomorphically to $\cu_a^\perp$.
It follows from induction and formula \refbq{} that the $T_{b,j}$ are bounded and
the $P_{b,j}$ are in the Schwartz class.    

Consider the Lax pair
\refeq[gb]$$\left[{\p\over \p x} + A,\,\, {\p\over \p t} + (M^{-1}b\l^jM)_+ \right]=0.$$
Set the coefficient of $\l^{j-k}$, $0\leq k<j$, in equation \refgb{} equals to zero
to get
\refeq[gd]$$\cases{[a,Q_{b,0}]=0,&\cr  (Q_{b,k})_x + [u,Q_{b,k}] +
[a,Q_{b,k+1}]=0,& $1\leq k < j$.\cr }$$ This defines the
$Q_{b,j}$'s.  The constant term gives 
\refeq[bi]$$u_t= (Q_{b,j})_x+ [u,Q_{b,j}] = [Q_{b,j+1}, a]$$
which is called the $j$-{\it th flow equation on $\cs_{1,a}$ defined by $b$}. Equation
\refgb{} can also be written as 
$$\eqalign{A_t &= ((M^{-1}b\l^j M)_+)_x + [A, (M^{-1}b\l^j M)_+]\cr
&= [M^{-1}b\l^j M, M^{-1}M_x]_+ + [A, (M^{-1}b\l^j M)_+]\cr &= 
[M^{-1}b\l^j M, A- M^{-1}a\l M]_+ + [A, (M^{-1}b\l^j M)_+]\cr &= [M^{-1}b\l^j M, A]_+
+ [A, (M^{-1}b\l^j M)_+]\cr &= [(M^{-1}b\l^j M)_-, A]_+ = [Q_{b,j+1},a].\cr}$$
It is clear that the following three statements are equivalent:
\item {(i)}  $[{\p\over \p x}+A, {\p\over \p t} +B]=0$,
\item {(ii)} the connection $1$-form $\o = Adx+Bdt$ is flat for all $\l$, i.e.,
$d\o=-\o\wedge \o$,
\item {(iii)} $\cases{E^{-1}E_x=A, &\cr E^{-1}E_t=B,&\cr}$ is solvable.

\ni
 So we have

\refclaim[hn] Proposition.  $A=a\l+ u$ is a solution of the $j$-th flow \refbi{} on
$\cs_{1,a}$ defined by
$b$ if and only if $$\o(x,t,\l)= (a\l+u)dx + (b\l^j + Q_{b,1}\l^{j-1} +\cdots + Q_{b,j})dt$$ is
flat on the $(x,t)$-plane for each $\l$. 

\refpar[hp] Definition. The one parameter family of connection $1$-form $\o$ defined in
Proposition \refhn{} is called the {\it flat connection associated to the solution\/}
$A$ of the
$j$-th flow.  The unique solution $E:R^2\times C\to GL(n,C)$ of 
$$\cases{E^{-1}E_x= a\l +u,&\cr
E^{-1}E_t=b\l^j + Q_{b,1}\l^{j-1}+ \cdots + Q_{b,j},&\cr E(0,0,\l)= I&\cr}$$
is called {\it the trivialization of the flat connection $\o$ normalized at the
origin\/} or {\it the trivialization of the solution $A$ at $(x,t)=(0,0)$\/}.

 When $a$ is regular, positive flows are the familiar hierarchy of commuting Hamiltonian
flows described by differential equations.  When $a$ is singular, positive flows  generate a
non-abelian Poisson group action. This will be described in section 8.  

\refpar[ho] Example. For $su(2)$ with $a=\diag(i,-i)$, $\cs_{1,a}$ is the set of $A$ of the
form $a\l+u$, where
$$u=\left(\matrix {0&f\cr -\bar f&0\cr}\right)$$ and $f:R\to C$ is in the Schwartz class. 
The first flow is the translation $u_t=u_x$, the second flow defined by $a$ is the non-linear
Schr\"odinger equation (NLS)
\refeq[lm]$$q_t={i\over 2}(q_{xx} + 2\n q\n^2 q),$$ and the positive flows are the hierarchy
of commuting flows associated to the non-linear Schr\"odinger equation.  

\ss

\refpar[hx] Example. For $a=\diag(a_1, \cdots, a_n)\in su(n)$ with  $a_1<\cdots <a_n$,
$\cs_{1,a}$ is the set of all $A=a\l + u$, where $u=(u_{ij})\in su(n)$ and $u_{ii}=0$ for all $1\leq
i\leq n$. The first flow on $\cs_{1,a}$
defined by $a$ is the translation
$$u_t=u_x.$$ The first flow on $\cs_{1,a}$ defined by $b=\diag(b_1, \cdots, b_n)$ ($b\not= a$)
is the $n$-wave equation ([ZMa1, 2]) for $u$:
$$(u_{ij})_t={b_i-b_j\over a_i-a_j} (u_{ij})_x + \sum_{k\not= i, j}
\left({b_k-b_j\over a_k-a_j} - {b_i-b_k\over a_i-a_k}\right) u_{ik} u_{kj}, \qquad
i\not=j.$$

\ms

If $a$ is singular and $[b,a]=0$, then the $j$-th flow on $\cs_{1,a}$  defined by $b$ is in general an
integro-differential equation. But the $j$-th flow on $\cs_{1,a}$ defined by $a$ is
again a differential operator:

\refclaim[bt] Proposition. $Q_{a,j}(u)$ is always a polynomial differential operator in $u$.

\proof  It is easy
to see that $Q_{a,1}=u$. We will prove this Proposition by induction. Suppose
$Q_{a,i}$ is a polynomial differential operator in $u$ for $i\leq  j$. Write
$$Q_{a,i}=P_{a,i}+T_{a,i}\,\,\in\,\, \cu_a^\perp + \cu_a$$ as before.  Using formula
\refbq{}, we see that $P_{a,j+1}$ is a polynomial differential operator in $u$. But we can not
conclude from formula \refbq{} that $T_{a,j+1}$ is a polynomial differential operator in $u$.  
Suppose $a$ has $k$ distinct eigenvalues $c_1, \cdots, c_k$. Then 
$$f(t)=(t-c_1)(t-c_2)\cdots (t-c_k)$$ is the minimal polynomial of $a$. So
$f(M^{-1}aM)=0$, which implies that the formal power series 
\refeq[bu]$$f(a+Q_{a,1}\l^{-1} + Q_{a,2} \l^{-2} +\cdots )=0.$$
Notice that $f'(a)$ is invertible and $T_{a,j+1}$ commutes with $a$.  Now
compare coefficient of $\l^{-(j+1)}$ in equation \refbu{} implies that $T_{a, j+1}$ can
written in terms of $a, Q_{a, 1}, \cdots , Q_{a,j}$. This proves that $Q_{a,j+1}$ is a
polynomial differential operator in $u$.
\qed

\ss

\refpar[ca] Example. For $u(n)$ with 
$$a=\pmatrix{iI_k&0\cr 0&-iI_{n-k}\cr},$$  
 $$\cs_{1,a}=\left\{a\l+u\,\,\bigg| \,\,u=\pmatrix{0&X\cr -X^*&0\cr},  X\in \cm_{k\times
(n-k)}\right\},$$ where $\cm_{k\times (n-k)}$ is the space of $k\times (n-k)$ complex
matrices.  Identifying
$\cs_{1,a}$ as
$\cs(R,\cm_{k\times (n-k)})$, then the bi-linear form
$$<u,v>=\int_{-\infty}^\infty \tr(uv)dx$$  on
$\cs(R,\cu_a^\perp)$ induces the following bi-linear form on $\cs(R,\cm_{k\times (n-k)})$:
$$<X,Y>= -\int_{-\infty}^\infty \tr(XY^*+X^*Y)dx.$$
The orbit symplectic structure on $\cs_{1,a}$ 
induces the following symplectic structure on $\cs(R, \cm_{k\times (n-k)})$:
$$w(X,Y)=<{i\over 2}X,Y>.$$ According to Propositions \refbt{}, the $j$-th flow
defined by
$a$ can be written down explicitly. For 
$$u=\pmatrix{0&B^*\cr -B^*&0\cr},$$ we have 
$Q_{a,0}=a$,  $Q_{a,1}= u$, 
$$Q_{a,2}=\pmatrix{{1\over 2i} BB^*& {i\over 2} B_x\cr {i\over 2} B^*_x& -{1\over
2i} B^*B\cr}.$$ The first three flows on $\cs(R,\cm_{k\times(n-k)})$ are 
$$\eqalign{B_t &= B_x\cr
B_t &= {i\over 2} (B_{xx} + 2 BB^\ast B)\cr
B_t &= -{1\over 4} B_{xxx} - {3\over 4} (B_xB^\ast B + BB^\ast B_x).\cr}$$ 
Notice that the second flow is the matrix non-linear Schr\"odinger equation
associated to
$Gr(k,C^n)$ by Fordy and Kulish [FK].  By Proposition \refhn{}, $B$ is a solution of the second
flow if and only if 
$$\left(a\l + u\right)dx +
\left(a \l^2 + u\l + Q_{a,2}\right)dt$$ is flat for all $\l$.   

\bs

\newsection Action of the rational loop group.\par

The rational loop group is used to construct the soliton data for the positive flows
discussed in section 5.  We first define a local action $\sharp$ of $G_-$ on
$C(R,\cg_+)$ via the dressing action. In general the $G_-$-action does not preserve
the space $\cs_{1,a}$ (because the Schwartz condition on $u$ for $A=a\l +u$ is not
preserved even locally).  However, we prove that the action
$\sharp$ of the subgroup
$G_-^m$ of rational maps in
$G_-$ leaves $\cs_{1,a}$ invariant. We also show that the factorization can be done explicitly.
In particular, the action $g_-\sharp A$ can be computed explicitly in terms of the
trivialization $E(A)$ of $A$ normalized at $x=0$. In fact, $g_-\sharp A$ is given by an
algebraic formula in terms of $E(A)$ and $g$. 

Let  $A\in C(R,\cg_+)$, and $E(x,\l)$ denote the trivialization of $A$ normalized at
$x=0$. Then the map $A\mapsto E$ identifies $C(R,\cg_+)$
with a subset of $C(R,G_+)$. (We write $E(x)(\l)= E(x,\l)$).

Given $f_-\in G_-$ and $A\in C(R,\cg_+)$, define 
$$f\sharp A= \tilde E^{-1} (\tilde E)_x,$$ where $\tilde E(x)= f_-\sharp E(x)$ is the
dressing action of $G_-$ on $G_+$ for each $x\in R$. In other words, we factor $$f_-
E(x)=
\tilde E(x) \tilde f_-(x)\in G_+\times G_-.$$ Clearly, this defines  a local action of
$G_-$ on $C(R,\cg_+)$.   For
$\xi_-\in \cg_-$, the corresponding infinitesimal vector field on $C(R,\cg_+)$ is
$$\tilde \xi_-(A)= -[A,(E(A)^{-1}\xi_-E(A))_-]_+.$$

\refclaim[gu] Proposition. 
Let $a$ be a fixed diagonal element in $u(n)$, and $C_{1,a}$ the space of all $A\in
C(R,\cg_+)$ such that $A(x)(\l)=a\l+u(x)$ and $u:R\to \cu_a^\perp$ is a smooth map.
Then $\xi_-\mapsto \tilde \xi_-$ defines an action of $\cg_-$ on
$C_{1,a}$. 

In general, the action of $G_-$ does not preserve the Schwartz condition for $\cs_{1,a}$. So it
does not define an action on $\cs_{1,a}$.  But  the subgroup $G_-^m$ of
rational maps does preserve the decay condition.

\refclaim[mj] Theorem. Let $G_-^m$ be the subgroup of rational maps $g\in G_-$. Then the
$\sharp$ action of $G_-^m$ on $C(R,\cg_+)$ leaves the space $\cs_{1,a}$ invariant. Moreover, 
let $g\in G_-^m$, $A\in \cs_{1,a}$, and $E$ the trivialization of $A$ normalized at $x=0$, then
\item {(i)} we can factor $gE(x)=\tilde E(x)\tilde g(x)\in G_+\times G_-^m$ and $g\sharp
A= \tilde E^{-1}(\tilde E)_x$,
\item {(ii)} $g\sharp A$ can be constructed algebraically from $E$ and $g$.

To prove this theorem, we first recall the following
result of the second author [U1]: 

\refclaim[mg] Proposition ([U1]).  Let $z\in C\setminus R$, $V$ a complex linear subspace of
$C^n$,  $\pi$ the projection of $C^n$ onto $V$, and $\pi^\perp= I-\pi$. Set 
\refeq[ra]$$g_{z,\pi}(\l)= \pi + {\l- z\over \l- \bar z} \pi^\perp. $$ Then
\item {(i)} $g_{z,\pi}\in G_-^m$,
\item {(ii)} $G^m_-$ is generated by $\{g_{z,\pi}\n z\in C\setminus R, \pi$ is a projection$\}$.
($g_{z,\pi}$ will be called a {\it simple\/} element).  

\refclaim[mh] Proposition. \item {(i)} Let $g(\l)= \Pi_{j=1}^r {\l- z_j\over \l- \bar z_j}$,
and
$A=a\l+u$. Then $g\in G_-^m$ and $g\sharp A=A$.  
\item {(ii)} Let $v_1, \ldots, v_k$ be a unitary basis of the linear subspace $V$, $\pi_j$
the projection of $C^n$ onto $Cv_j$, and $\pi$ the projection onto $V$. Then 
$$\prod_{j=1}^k g_{z,\pi_j} = \left({\l- z\over \l-\bar z}\right)^{k-1}g_{z,\pi}.$$

\proof  Statement (i) follows from the fact that $g$ commutes with $G_+$ and $G_-$.
Statement (ii) follows from a direct computation.\qed

The above two Propositions imply that to prove Theorem \refmj{} it suffices to prove
$g_{z,\pi}\sharp A\in \cs_{1,a}$, where $\pi$ is the projections onto a one dimensional
subspace.  First, we give an explicit construction of $g_{z,\pi}\sharp A$. 

\refclaim[mi] Theorem. 
Let $A=a\l + u\in \cs_{1,a}$, and $E$ the trivialization of $A$
normalized at $x=0$.  Let $z\in C\setminus R$, $V$ a complex linear subspace of $C^n$, and $\pi$
the projection onto $V$.  Set 
$$\eqalign{\tilde V(x)&= E(x,z)^*(V),\cr \tilde \pi(x) &= {\rm the\, projection\, of\,} C^n\, {\rm 
onto \,} \tilde V(x),\cr \tilde E(x,\l) &= g_{z,\pi}(\l)E(x,\l) g_{z,\tilde
\pi(x)} \, ^{-1}\cr &= \left(\pi + {\l-z\over \l-\bar z}\pi^\perp\right)E(x,\l)\left(\tilde \pi(x) +
{\l-\bar z\over\l-z} \tilde \pi(x)^\perp\right).\cr}$$  Then: 
\item {(i)} $g_{z,\pi}\sharp E =\tilde E$.
\item {(ii)} $\tilde \pi^\perp(\tilde\pi_x + (a\bar z + u)\tilde \pi) = 0$.
\item {(iii)} If $v:R\to C^n$ is a smooth map such that $v(x)\in \tilde V(x)$ for all $x\in R$,
then $v_x(x)+ (a\bar z+ u) v(x) \in \tilde V(x)$ for all $x$. 
\item {(iv)} $g_{z,\pi}\sharp A= A + (z-\bar z)[\tilde \pi, a]$.

\proof First we claim that $\tilde E(x,\l)$ is holomorphic for $\l\in C$. By definition, $\tilde E$ is
holomorphic in $\l\in C\setminus \{z,\bar z\}$ and has possible poles at $z, \bar z$ with order
one.  The residues of $\tilde E$ at these two points can be computed easily:
$$\eqalign{{\rm Res\/}(\tilde E, z) &=(z-\bar z) \pi E(x,z)\tilde \pi^\perp(x),\cr
{\rm Res\/}(\tilde E,\bar z) &=(\bar z- z) \pi^\perp E(x,\bar
z)\tilde \pi(x).\cr}$$   Since $A(x,\bar z)^* + A(x,z)=0$ and
$E(0,\l)=I$, $E(x,\bar z)^*E(x,z)=I$. This implies that  $$\tilde V(x)
=E(x,z)^*(V)= E(x,\bar z)^{-1}(V).$$ So both residues are zero, and 
the claim is proved. In particular, we have $g_{z,\pi}E(x) = \tilde E(x) g_{z,\tilde
\pi(x)}\in G_+\times G_-$. This implies (i).

By Proposition \refgu{}, $\tilde E^{-1}(\tilde E)_x = a\l + \tilde
u(x)$ for some smooth $\tilde u:R\to \cu_a^\perp$. We get from the formula for $\tilde E$ that
\refeq[mm]$$\eqalign{a\l+\tilde u = &g_{z,\tilde\pi}(a\l+u)g_{z,\tilde\pi}^{-1} -
(g_{z,\tilde\pi})_xg_{z,\tilde\pi}^{-1}\cr
&=\left(\tilde \pi + {\l-z\over \l-\bar z} \tilde \pi^\perp\right)(a\l
+u)\left(\tilde
\pi + {\l-\bar z\over
\l- z} \tilde \pi^\perp\right)\cr &\qquad - \left(\tilde \pi_x + {\l-z\over \l-\bar z}
\tilde
\pi_x^\perp\right) \left(\tilde \pi + {\l-\bar z\over \l- z} \tilde
\pi^\perp\right).\cr}$$
 Since the left hand side is holomorphic at $\l= z$, the residue
of the right hand side at
$\l= z$ is zero. This gives $(\tilde \pi(az+u)-\tilde \pi_x)\tilde \pi^\perp = 0$, which is
equivalent to (ii). 

Statement (iii) follows from (ii) since 
$$\eqalign{v_x+(a\bar z+u)v&= (\tilde \pi (v))_x + (a\bar z+u)v\cr
&= \tilde \pi_x v+ \tilde \pi v_x + (a\bar z + u) v\cr 
&=(\tilde \pi_x+(a\bar z+u)\tilde \pi)(v) + \tilde \pi (v_x)\in \tilde V(x).\cr}$$ 

To prove (iv), we multiply $g_{z,\tilde\pi}$ to both sides of equation \refmm{}
and  get
$$\eqalign{((\l-\bar z)\tilde \pi + (\l-z)\tilde \pi^\perp)(a\l+ u) - ((\l-\bar z)\tilde \pi_x +
(\l -z)\tilde \pi_x^\perp)&\cr = (a\l + \tilde u)((\l-\bar z)\tilde \pi + (\l - z) \tilde
\pi^\perp).&\cr}$$ Set $\l= z$ and $\l=\bar z$ in the above equation, we get
\refeq[gk]$$\cases{\tilde\pi (az+ u)- \tilde \pi_x = (az+\tilde u) \tilde \pi,&\cr
\tilde\pi^\perp(a\bar z + u)- \tilde\pi_x^\perp = (a\bar z + \tilde
u)\tilde\pi^\perp.&\cr}$$
Add the two equations in \refgk{} to get 
\refeq[gl]$$\tilde u= u + (z-\bar z)[\tilde\pi, a].\qed$$

\refclaim[gm] Theorem. The map $\tilde \pi$ in Theorem \refmi{} is the solution of
the following ordinary differential equation:
\refeq[fr]$$\cases{(\tilde \pi)_x + [az+u,\tilde \pi] = (\bar z-z)[\tilde\pi, a]\tilde\pi, &\cr
\tilde\pi^*=\tilde\pi, \quad \tilde\pi^2=\tilde\pi, \quad \tilde\pi(0)=\pi.&\cr}$$
Moreover, if $\tilde \pi$ is a solution of this equation then $[\tilde \pi, a]$ is in the
Schwartz class. 

\proof  Substitute equation \refgl{} into the first equation of \refgk{} to get the equation
\reffr{}. 

By Proposition \refmh{}, to prove $[\tilde \pi,a]$ is in the Schwartz class it suffices to prove
it for the case when
$V$ is of one dimensional. By Theorem \refmi{} (iii) there exist smooth maps $v:R\to C^n$ and
$\phi:R\to C$, such that
$v(x)$ spans the linear subspace $\tilde V(x)$ and 
$$v_x+(a\bar z + u) v = \phi v.$$
Set $w= \exp(-\int_{-\infty}^x \phi) v$. Then  $w(x)$ generates
$\tilde V(x)$ and 
\refeq[mk]$$w_x+(a\bar z + u) w= 0.$$ 
We may assume that  $$a=\diag(ic_1, \cdots, ic_n), \qquad c_1\leq \cdots \leq c_n.$$ Let
$\psi_j:R\to C^n$ denote the solution of
$$\cases{(\psi_j)_x + (a\bar z + u) \psi_j= 0,&\cr \lim_{x\to -\infty} e^{-ic_j\bar z x}
\psi_j(x) = e_j,&\cr}$$ where $\{e_1, \cdots, e_n\}$ is the standard basis of
$R^n$. The construction of the $\psi_j$ is a standard textbook part of the scattering theory.
Then
$\psi_1, \cdots,
\psi_n$ form a basis of the solution for equation \refmk{}. So there exist constants
$b_1,\cdots, b_n$ such that
$ w= \sum_{j=1}^n b_j \psi_j$. Let $\bar z= r+is$ with $s>0$ and choose $j$ to be the
smallest integer such that $b_j\not=0$. Then 
$$e^{-ic_j\bar zx}w=\sum_{k\leq j}^n e^{-ic_j\bar z x} b_k\psi_k 
= \sum_{k\leq j}^n  e^{i(-c_j+c_k)\bar z x} b_k (e^{-ic_k \bar z x} \psi_k)$$
Since $\lim_{x\to -\infty} e^{i(-c_j+c_k)\bar z x} = 0$ if $c_k<c_j$, we get 
$$\lim_{x\to -\infty} e^{-ic_j \bar z x}w(x) = \sum_{c_k=c_j} b_ke_k,$$ which is an
eigenvector for $a$. So $\lim_{x\to -\infty}[\tilde\pi(x),a]=0$. 
Moreover, 
$$e^{-ic_j\bar z x}w(x) = \sum_{c_k=c_j} b_ke_k + \sum_{c_j<c_k}
e^{i(-c_j+c_k)\bar z x}b_ke^{-ic_k\bar z x}\psi_k.$$ Since $\lim_{x\to -\infty}
e^{-ic_k\bar z x}\psi_k(x)= e_k$, 
$$\bigg| e^{-ic_j\bar z x}w(x)-\sum_{c_k=c_j} b_ke_k\bigg| =
e^{(-c_j+c_m)x\Im(z)}\co(1),$$ where $c_m$ is the next non-zero term. Hence
$\tilde
\pi(x)-\lim_{x\to -\infty}\tilde
\pi(x)$ decays exponentially, so $[\tilde \pi, a]$ also decays exponentially as $x\to
-\infty$.    Similarly, we can
prove that  $ [\tilde \pi(x),a]=0$ decays exponentially when $x\to \infty$.

From equation \reffr{}
$$\n \tilde \pi_x\n \leq 2\n z\n  \n [a,\tilde \pi]\n + 4\n u\n.$$
So $\tilde \pi_x$ decays like $u$. Repeated differentiation of equation \reffr{} give
the desired result. In fact, $\tilde \pi_x\in \cs(R)$ as well. 
\qed

If $U$ is a matrix whose columns form a basis of $V$, then the
projection of $C^n$ onto $V$ is $\pi= U(U^*U)^{-1}U^*$. This follows from
elementary linear algebra.  So we have:

\refclaim[gn] Corollary. Let $V$ be a $k$-dimensional linear subspace of $C^n$, and $U$ a
matrix whose columns form a basis of $V$. Then 
$g_{z,\pi}\sharp A = A + (\bar z-z)[\tilde \pi,a]$, where 
\refeq[hr]$$\tilde \pi(x)= E^*(x,z)U(U^*E(x,z)E^*(x,z)U)^{-1}U^* E(x,z).$$

\refpar[fk] Proof of Theorem \refmj{}. 

Given $g\in G_-^m$,  write
$g$ as product of simple elements
$\prod_{j=1}^k g_{z_j,\pi_j}$ (note that the factorization of $g$ into simple elements is
not unique, for example see Proposition \refmh{} (ii)). Use Theorem \refmi{} to
see that
$g\sharp A\in \cs_{1,a}$, and  
$g\sharp E$ and $g\sharp A$ are obtained by algebraic formulae from $E$ and $g$.
\qed

\bs

\newsection Scattering data and Birkhoff decomposition.\par

The asymptotically constant case is the case standardly treated in the soliton literature. We
obtain many hints of how to describe the theory, since most of what we need is
already contained in scattering theory literature ([BC1, 2], [Zh1, 2]).  The main
purpose of this section is to give a homogeneous structure for the space of
scattering data, to obtain the Inverse Scattering Transform using the standard
Birkhoff decompositions,  and to relate the action of the rational loop group
described in section 6 to the scattering data in a natural and simple way.

We first review results of Beals-Coifman ([BC1, 2]) and Zhou ([Zh2]) on scattering theory for
$n\times n$ first order linear system. Let 
$$a=\diag(ia_1,\cdots, ia_n)\in u(n), \qquad a_1\leq a_2 \leq \cdots \leq a_n,$$
and $A=a\l+u\in \cs_{1,a}$. Consider the linear system
\refeq[gx]$$\cases{\psi_x=\psi(a\l+u),&\cr \lim_{x\to -\infty}e^{-a\l x}\psi(x,\l)= I,&\cr 
m(x,\l)=e^{-a\l x}\psi(x,\l)\,\, {\rm is\, bounded\, in \,} x.\cr}$$
($m$ will be called the {\it normalized (matrix) eigenfunction\/} of $A$). 

\refclaim[gy] Theorem ([BC1, 2],[Zh2]).  Given $A=a\l+u\in \cs_{1,a}$ there exists a bounded
discrete subset $D$ of $C\setminus R$ such that the normalized eigenfunction 
$m(x,\l)= e^{-a\l x}\psi(x,\l)$ is
holomorphic in $\l\in C\setminus (R\cup D)$ and has poles at  $z\in D$.  Moreover,  there exists a
dense open subset $\cs'_{1,a}$ of
$\cs_{1,a}$ such that for  $A=a\l+u\in\cs_{1,a}'$, the normalized eigenfunction $m(x,\l)$ 
satisfies the following conditions:
\item {(i)} The subset $D$ is finite, and $m$ has only simple poles at $z\in D$,
\item {(ii)} The matrix function $m$ can be extended
smoothly to the real axis from the upper and lower half $\l$-plane, 
\item {(iii)} As a function of $\l$, $m$  has an asymptotic expansion at $\l=\infty$.

The open dense subset $\cs'_{1,a}$ contains all $\{u\in\cs_{1,a}$ such that the $L^1$-norm of $u$
is less than $1$ and  all $u$ with compact support. 

\refclaim[] Theorem ([BC1,2]). Let $m$ be the normalized eigenfunction of
$A=a\l+u\in
\cs_{1,a}$, and $b\in u(n)$ such that $[a,b]=0$. Set $Q_b= m^{-1}bm$. Then $Q_b$
has an asymptotic expansion at $\l=\infty$:
$$Q_{b,\l} \sim b + Q_{b,1}\l^{-1} + Q_{b,2}\l^{-2} + \cdots .$$  Moreover,
\item {(i)} $(Q_{b,j})_x+ [u, Q_{b,j}] = [Q_{b, j+1},a]$.
\item {(ii)} The $j$-th flow $u_t=[Q_{b,j+1},a]$ is symplectic with respect to the
symplectic structure $w(v_1,v_2)=<-\ad(a)^{-1}(v_1), v_2>$. 

\ms 
Recall
$$\eqalign{G_+&=\{g:C\to GL(n,C)\n g\,\, {\rm is\, holomorphic\/},\,\, g(\bar\l)^*g(\l)=I\},\cr
G_-&=\{g:\co_\infty\to GL(n,C)\n g \,\,{\rm is\, holomorphic\/},\,\,
g(\bar\l)^*g(\l)=I, g(\infty)=I\}.\cr}$$ 
Since $m(x,\l)$ is not holomorphic at $\l=\infty$, we must change $G_-$, and restrict $G_+$ to
have a singularity at $\l=\infty$ of the type exp(polynomial). 

We are motivated by Theorem \refgy{} to choose  a different negative
group $D_-$:

\refpar[hb] Definition.   Let $D_-$ denote the group of meromorphic maps $f$ from
$C\setminus R$ to
$GL(n,C)$ satisfying the following conditions:
\item {(i)} $f(\bar \l)^*f(\l)=I$.
\item {(ii)} $f$ has a smooth extension to the closure $\bar C_\pm$, i.e.,
$f_\pm(r)=\lim_{s\searrow 0} f(r\pm is)$ exists and is smooth  for $r\in R$, (since
$f(\bar\l)^*f(\l)=I$, we have  $f_-(r)=(f_+(r)^*)^{-1}$).
\item {(iii)} $f$ has an asymptotic expansion at $\infty$.
\item {(iv)} $f_+-I$ lies in the Schwartz class modulo unitary maps. In other words,
if we factor
$f_+=v_+h_+$ with $v_+$ unitary and $h_+$ upper triangular then $h_+-I$ is in the Schwartz
class. 

\ms

Let $m(x,\l)$ be the normalized eigenfunction for $A\in \cs_{1,a}'$, and $E$ the trivialization of
$A$ normalized at $x=0$. Since both $e^{a\l x} m(x,\l)$ and $E(x,\l)$ satisfy the ordinary
differential equation in $x$:
$$E^{-1}E_x= (e^{a\l x}m)^{-1}(e^{a\l x}m)_x = A,$$ there exists $f$ such that 
$$e^{a\l x} m(x,\l)= f(\l) E(x,\l).$$ In fact, $f(\l)=m(0,\l)$. By Theorem \refgy{}, $f\in D_-$.   

Beals and Coifman [BC1,2] defined the scattering data of $A=a\l+u\in \cs_{1,a}'$ to
be the map
$S:R\cup D\to GL(n)$: for $z\in D$, $S(z)$ is the element in $GL(n,C)$ such that
$$(I-(\l-z)^{-1}e^{azx}S(z)e^{-azx} ) m(x,\l)$$ has a removable singularity at $\l=z$, and for
$r\in R$, $S(r)= v_-^{-1}(r)v_+(r)$, where $$\lim_{s\searrow 0}
m(x,r+is)=e^{arx}v_\pm(r)e^{-arx}m(x,r).$$ They prove:
\item {(i)} The map sending $A$ to $S$ is injective. 
\item {(ii)}  If $u(x,t)$ is a solution of the $j$-th flow on $\cs_{1,a}$ defined by
$b$, $u_t=[Q_{b,j+1},a]$, and 
$S(\l, t)$ is the corresponding scattering data, then 
$$S_t(\l,t)= [S(\l,t), \l^j b].$$
In particular,   $S(\l,t)= e^{-b\l^j t}S(\l,0)e^{b\l^j t}$. 

\ni
We note that 
scattering data $S$ for $A$ is determined by $f(\l)=m(0,\l)$. In fact,   
$S(r)= f_-(r)^{-1}f_+(r)$ for $r\in R$ and $S(z)$ can be obtained from the residue of $f(\l)$
at $z\in D$. 
\ss

\refpar[fb] Remark. The rational group $G_-^m$ defined in section 6 is a subgroup of
$D_-$. 
\ss

Instead of using $S$ as the scattering data, we use the left coset $H_-f$ in $D_-/H_-$ as
the scattering data of $A$, where $H_-$ is the subgroup of $h\in D_-$ that commutes
with $a$. We will call $[f]= H_-f$ the {\it scattering coset of $A$\/}. One advantage of
using the scattering cosets is that the inverse scattering transform can be
obtained from the standard Birkhoff Decomposition Theorems. Another advantage
is that the natural action of the subgroup
$G_-^m$ of rational maps on
$D_-/H_-$ on the right by multiplication induces the action of
$G_-^m$ on $\cs_{1,a}$ defined in section 6. To explain this, we first prove a decomposition
theorem.

\refclaim[mb] Theorem. Let $D^c_-$ denote the subgroup of $v\in D_-$ such that $v$ is
holomorphic in  $C\setminus R$.   Then any $f\in D_-$ can be uniquely factored into 
$$f= gh = \tilde h \tilde g,$$
where $g,\tilde g\in G_-^m$ and $h,\tilde h\in D^c_-$.  Moreover, the multiplication map
$$G_-^m\times D_-^c\to D_-$$ is a diffeomorphism. 

This theorem is the real line version of the Birkhoff decomposition theorem, which can be seen 
by transforming  the domain $C_+$ to the unit disk and real axis to the unit circle $S^1$ by a
linear fractional transformation. To be more precise, let $LGL(n,C)$ denote the loop group of
smooth maps from $S^1$ to $GL(n,C)$, and $L^+GL(n,C)$ the group of maps $g\in
LGL(n,C)$ such that $g$ is the boundary value of a holomorphic map
$$g:\{z\n \n z\n <1\}\to GL(n,C)\}.$$ Let  $\W U(n)$ denote the based loop group of
maps $g:S^1\to U(n)$ such that $g(-1)=I$,  Recall that the standard Birkhoff
Decomposition Theorem (cf. [PrS] p. 120, Theorem 8.1.1) is:

\refclaim[gs] Birkhoff Decomposition Theorem. Any $g\in LGL(n,C)$ can be factored
uniquely as
$$g= g_+ g_- = h_- h_+,$$ where $g_+, h_+\in L^+GL(n,C)$ and $g_-, h_-\in \W U(n)$.
In other words, the multiplication map 
$$L^+GL(n,C)\times \W U(n)\to LGL(n,C)$$ is a diffeomorphism. 

A direct computation shows:

\refclaim[ma] Proposition. Given  $g:S^1\to GL(n,C)$, define 
$\Phi(g):R\to GL(n,C)$  by $\Phi(g)(r)= g( {1+ir\over 1-ir})$. Then
\item {(i)} $g$ is smooth if and only if $\Phi(g)$ is smooth and has the same asymptotic
expansions at $-\infty$ and $\infty$,
\item {(ii)} $g-I$ is infinitely flat at $z=-1$ if and only if
$\Phi(g)-I$ is in the Schwartz class,
\item {(iii)}  $g:C\to GL(n,C)$ satisfies the reality condition $g(1/\bar
z)^*g(z)= I$ if and only if $\tilde g(\l)=g({1+i\l\over 1-i\l})$ satisfies the reality
condition
$\tilde g(\bar\l)^*\tilde g(\l)=I$. 

\refclaim[] Corollary.  The group $D_-$ is isomorphic to the group of smooth loops
$g:S^1\to GL(n,C)$  that are boundary values of meromorphic maps with
finitely many poles in $\n z\n <1$ and $g^*g-I$ is infinitely flat at $z=-1$. 

As a consequence of Theorem \refgs{} and Proposition \refma{}, we have 

\refclaim[ch] Corollary. If $f:R\to GL(n,C)$ is smooth and has an asymptotic expansion at
$\l=\infty$, then $f$ can be factored
$$f=gv,$$ where $g$ is unitary and $v$ is the boundary value of a holomorphic map on $C_+$.

\refpar[fl] Proof of Theorem \refmb{}. 

It follows from Corollary \refch{} that given  $f\in D_-$, we can factor  $f_\pm$
$$f_\pm(r)=h_\pm(r)g_\pm(r),\qquad r\in R$$ where $h_\pm$ is the
boundary value of a holomorphic map $h$ on $C_\pm$ and $g_\pm$ is a smooth map from
$R$ to $U(n)$. It follows from $f_- =(f_+^*)^{-1}$ that we have  $g_+=g_-$ and
$h(\bar\l)^*h(\l)=I$.  Write 
$f= hg$. Since $f$ is meromorphic and $h$ is holomorphic in $C_+$, $g$ is
meromorphic in $C_+$. However, $g(r)^*g(r)=I$ for $r\in R$ implies that $g$
extends holomorphically across the real axis. So $g$ is meromorphic in $C$ and
bounded near infinity. This implies that $g$ is rational, i.e., $g\in G_-^m$. 
\qed

Recall that $a=\diag(ia_1, \cdots, ia_n)\in u(n)$ is a fixed diagonal matrix, and
$G_+$ is the group of holomorphic maps $g:C\to GL(n,C)$.  

\refclaim[mc] Theorem.  Let $f\in D_-$, $k$ a positive
integer, and $b\in u(n)$ such that $[a,b]=0$. Let $e_{b,k}(x)(\l)=e^{b\l^k x}$. 
Then there exists a unique
$E(x,\l)$ and $M(x,\l)$ such that 
$$f^{-1}e_{b,k}(x) = E(x,\cdot) M^{-1}(x,\cdot)\in G_+\times D_-.$$

\proof Write $f=hg$ as in Theorem \refmb{} with $h\in D^c_-$ and $g\in G_-^m$.
Write
$h=pv$, where $p$ is upper triangular and $v$ is unitary. By
definition of $D_-$, $p-I$ is in the Schwartz class when restricted to the real axis
in the $\l$-plane.  So
$e^{-1}_{b,k}(x)p^{-1}e_{b,k}(x)$ has an asymptotic expansion at $r=\pm\infty$ for
each $x$.  Write 
$e^{-1}_{b,k}(x)p^{-1}e_{b,k}(x) = \tilde v(x)\tilde h(x)$,
where $\tilde v$ is unitary and $\tilde h$ is the boundary value of a holomorphic map on 
$C_+$. Notice $f, p, v$ and $h$ do not depend on $x$, whereas the rest of the
matrix functions do depend on $x$.  So
$$h^{-1}e_{b,k}(x)=v^{-1}e_{b,k}(x)\tilde v(x) \tilde h(x) =B(x)\tilde h(x),$$ where
$B(x)=v^{-1}e_{b,k}(x)\tilde v(x)$ is unitary.  Both $\tilde h(x)$ and $h$
are holomorphic in $\l\in C_+$, and $e_{b,k}(x)$ is holomorphic in $C_+$. Hence $B(x)$ is
holomorphic in $\l\in C_+$. However, $B(x)$ is unitary hence it is  holomorphic in
$\l\in C$.  

Next we claim that we can factor  $g^{-1}B(x)=E(x)g_1^{-1}(x)$ with $E$ holomorphic
in
$C$ and $g_1\in G_-^m$.  This can be proved exactly the same way as Theorem
\refmj{}.  Then $$\eqalign{f^{-1}e_{b,k}(x) &= (hg)^{-1}e_{b,k}(x) = g^{-1} h^{-1}
e_{b,k}(x)\cr & = g^{-1}B(x)\tilde h(x) = E(x)g_1^{-1}(x)\tilde h(x) =
E(x)M^{-1}(x),\cr}$$ which finishes the proof. \qed

\refpar[hd] Definition. A matrix
$q$ is called $a$-diagonal if
$q_{jk}=0$ whenever $a_j\not=a_k$, $q$ is (strictly) upper $a$-triangular if
$q_{jk}=0$ whenever
$a_j>a_k$  (and $q_{jk}=0$
or $I$ if $a_j=a_k$), and  $q$ is (strictly) lower $a$-triangular if $q_{jk}=0$ whenever
$a_j<a_k$ (and $q_{jk}=0$ or $I$ if $a_j=a_k$).  Let $q_d$
denote the {\it $a$-diagonal projection of\/} $q$, i.e.,
$$(q_d)_{ij}=\cases{q_{ij},&if $a_i=a_j$,\cr 0,& if $a_j\not=a_k$.\cr}$$ 

\refclaim[he] Proposition. Any $f\in D_-$ can be factored
uniquely as $$f= pv = q \tilde v,$$ where $p$ is upper $a$-triangular, $q$ is lower
$a$-triangular, $v$ and $ \tilde v$ are unitary, and the $a$-diagonal projections $p_d, q_d$ are
holomorphic in $C_\pm$. 

\proof Write $g=p_0v_0$, where $p_0$ is upper $a$-triangular and $v_0$ is unitary.
Such $p_0, v_0$ are not unique because an element in $U_a=\{y\in U(n)\n ay=ya\}$ is both
$a$-triangular and unitary.   Write $p_0=p_1p_2$, where 
$p_1$ is  strictly upper $a$-triangular and $p_2$ is in $a$-diagonal.    Factor $p_2=p_3h$,
where
$p_3$ is holomorphic in $C_+$ and $h$ is unitary. Then $g=p_1p_3hv_0=pv$, where
$p=p_1p_3$ is upper $a$-triangular and $v=hv_0$ is unitary.  Since $p_1$ is strictly
upper $a$-triangular,
$p_d=p_3$ is holomorphic. 
\qed

To study how the Birkhoff factorization of Theorem \refmc{} depends on parameter $x$, we
introduce the class of Schwartz maps from $[r_0,\infty)$ to a Hilbert space.  Let $H$
be a Hilbert space, a map $\phi:[r_0,\infty)\to H$ is in $\cs([r_0,\infty),H)$ if for each
pair of integers $(m,s)$ there exists a constant $c_{m,s}$ such that 
$$\left\Vert \left({d\over dx}\right)^m \phi(x)\right\Vert \leq {c_{m,s}\over (1+\n x\n)^s}.$$
Let $H_1$ denote the Sobolev space for maps from $R^+=[0,\infty)$ to
$\cu_a^\perp$. In other words,
$u\in H_1$ if
$$\eqalign{\N u\N^2_1 &=\int_0^\infty \left(\left\Vert{du\over dr}\right\Vert^2 + \N
u\N^2\right) dr
\cr &=\int_0^\infty (y^2+1)\N \hat u\N^2 dy \quad < \infty,\cr}$$ where $\hat u$ is the
Fourier transform of $u$. 

The following is a functional analytic extension of Birkhoff decomposition.

\refclaim[ft] Theorem.  Let $I+D(x,\cdot)=(I+h(x,\cdot))V(x,\cdot)$ be the Birkhoff
decomposition, where
$h(x)(r)$ is the boundary value of a holomorphic map in  the upper half plane and
$V(x)(r)$ is unitary. If
$D\in\cs(R^+,H_1)$, then $h$ and $V-I$ are in $\cs(R^+,H_1)$.

\proof This should be regarded as an implicit function theorem. It is based on the two facts
about the Sobolev space $H_1$. The first is that $H_1$ is an algebra under multiplication
and $\exp:H_1\to H_1$ is smooth.  The second is that the linear Birkhoff
decomposition can be defined using the Fourier transform $\cf$. Let $\Psi:
L^2(R)\to L^2(R)$ denote the linear operator defined by
$$\Psi(f)(y)=\cases{f(y)+f(-y)^*, & if $y\geq 0$, \cr 0, & otherwise,\cr}$$ and let
$\pi_+:H_1\to H_1$ be the bounded linear map 
$$\eqalign{\pi_+&= \cf^{-1}\Psi \cf, \qquad {\rm i.e.\/}, \cr
\pi_+(f)(r)&= \int_0^\infty (\hat f(y)+ \hat f(-y)^*) e^{iry}dy.\cr}$$
We claim that $f=\pi_+(f) + (I-\pi_+)(f)$ is the linear Birkhoff Decomposition, or
equivalently,
$\pi_+(f)$ is the boundary value of a holomorphic map on $C_+$ and $(I-\pi_+)(f)$ is
is in $u(n)$. To see this, we note that $\pi_+(f)$ is the boundary value of the
holomorphic map $$\l\in C_+ \mapsto \int_0^\infty (\hat f(y)+\hat f(-y)^*)e^{i\l
y}dy.$$ Then $$\eqalign{(I-\pi_+)(f)&=f(r) - \int_0^\infty (\hat f(y)+\hat
f(-y)^*)e^{iry}dy\cr &=\int_{-\infty}^\infty \hat f(y)e^{iry}dy -  \int_0^\infty (\hat f(y)+\hat
f(-y)^*)e^{iry}dy \cr &=\int_{-\infty}^0 \hat f(y)e^{iry}dy -\int_0^\infty \hat
f(-y)e^{iry}dy.\cr}$$ It follows that $(I-\pi_+)(f)^*=-(I-\pi_+)(f)$. 

Due to the linearity of $\pi_+$, it is easy to see
that this extends to the parameter version in $x$. We write this as $$\pi_+: \cs(R^+,
H_1)\to
\cs(R^+,H_1).$$ Now the Birkhoff decomposition is a non-linear operator. However
we are near the identity, so it can be regarded as a perturbation of the linear
operation because the exponential map is smooth on $H_1(R)$.  

Let $Y:\cs(R^+,H_1)\to \cs(R^+,H_1)$ be the map defined by
$$Y(f)=e^{\pi_+(f)}e^{(I-\pi_+)(f)}.$$
Given $D$, we wish to find $\tilde D$ such that 
$$(I+D)=\exp(\pi_+(\tilde D))\exp((I-\pi_+)\tilde D) = Y(\tilde D).$$ Since $dY_0=
I$, for $x$ sufficiently large
$$\N\pi_+(\tilde D)\N_1 \leq \N \tilde D\N_1 \leq O_s,$$ where 
$O_s(x)=  c \N D\N_1\leq {cc_s\over
(1+\n x\n)^s}$. 
The estimate on derivatives in $x$ is more difficult. Let
$I+h= \exp(\pi_+(\tilde D))$. Then 
$$(I+h)^{-1}D_x V^{-1} = (I+h)^{-1}h_x + V_x V^{-1}.$$
On the right, the first term is holomorphic in the upper half plane, the second term is
unitary. Hence 
$$(I+h)^{-1}h_x = \pi_+((I+h)^{-1}D_x V^{-1})=
\pi_+((I+h)^{-1}D_x(I+D)^{-1}(I+h)).$$
Or $$h_x= (I+h) \pi_+((I+h)^{-1}D_x(I+D)^{-1}(I+h)).$$ Certainly, $D_x(I+D)^{-1}\in
\cs(R^+,H_1)$, $ \pi_+$ is linear, and $H_1$ is an algebra.  Using the Leibnitz rule
repeatedly, we can obtain
$$\bigg\| \bigg({\p\over \p x}\bigg)^m h(x)\bigg\|_1 \leq C_m(h) \, \max_{j\leq l}\bigg\|
\bigg({\p^m D\over \p x^m}(x)\bigg)\bigg\|_1,$$ where $C_m(h)=C(\N h\N, \cdots, \N (\p/\p
x)^{m-1} h\N)$. Estimates in the Schwartz topology follows by induction on $m$. \qed

\refpar[fu] Remark. The awkwardness of this proof reminds one that the classical use of the
Schwartz space is probably not as natural for the analysis as various choices of Hilbert or
Banach spaces in $x$ would be. The above proof would then be a straight forward
use of the usual implicit function theorem (rather than a reproof). Notice that in 
fact
$H_1(R)$ could be replaced by any $H_k(R)$, $k>{1\over 2}$. 

Recall that $e_{a,1}(x)(\l)=e^{a\l x}$. 

\refclaim[hc] Theorem. In the Birkhoff factorization of Theorem \refmc{}
$$f^{-1}e_{a,1}(x) = E(x)M^{-1}(x)\in G_+\times D_-.$$ We have in addition the
following properties
\itemitem {(i)}  $E^{-1}E_x= A$, where $A(x,\l)= a\l + u(x)$ for some $u\in
\cs(R,\cu_a^\perp)$,
\itemitem {(ii)} $M_{\pm \infty} \in H_-$, where $M_{\pm \infty}(\l)=\lim_{x\to
\pm\infty} M(x,\l)$,
\itemitem {(iii)} if $\l$ is not a pole of $f$ then $\tilde M_\pm(\cdot, \l)= M(\cdot, \l)-
M_{\pm\infty}(\l)$ is in the Schwartz class.\ei 

To prove this theorem, we need the following Lemma:

\refclaim[fs] Lemma. Given $q\in \cs(R)$ and $\b<0$ a constant, and set
$$\eqalign{\D(x)(r)&=\int_{-\infty}^0 q(y+\b x)e^{iry}dy,\cr
\xi(x)(r) &= \int_0^\infty q(y+\b x) e^{iry} dy.\cr}$$ Then 
\item {(i)} $\D\in \cs(R^+,H_1)$,
\item {(ii)} for any integer $m\geq 0$ there exists a constant $c_m$ such that  
$$\left|{\p^m
\xi\over \p x^m}(x,r) \right|
\leq c_m.$$ 

\proof Write $\D(x,r)=\D(x)(r)$ and $\xi(x,r)=\xi(x)(r)$. Then $$\bigg({\p\over \p
x}\bigg)^m\D(x,r)=\b^m\int_{-\infty}^0 {\p^m q\over
\p y^m}(y+\b x)e^{iry} dy.$$ Since $q\in \cs(R^-)$, 
$$\left|\bigg({\p \over \p y}\bigg)^m q(y)\right|\leq {c_{m,s}\over (1+\n y\n)^s}.$$
But by the Plancherel Theorem
$$\eqalign{\bigg\| \bigg({\p\over \p x}\bigg)^m \D(x,\cdot)\bigg\|_1^2 
&=\b^{2m}\int_{-\infty}^0 (1+y^2)\big({\p^m q\over \p y^m}(y+\b x) \big)^2 dy\cr 
&\leq \b^{2m}\int_{-\infty}^0 {(1+\n y\n^2)c_{m,s}^2\over (1+ \n y\n + \n \b x\n)^{2s}} dy\cr
&\leq {\b^{2m} c_{m,s}^2\over (1+ \n \b x\n)^{2s-3}}\cr}.$$
An adjustment of the constants completes the proof of (i). A straight forward computation
gives (ii). \qed

\refpar[rb] Proof of Theorem \refhc{}. 

 Take the variation with respect to $f\in D_-$ in the formulas
$$\eqalign{&f^{-1}(\l)e^{a\l x}= E(x,\l)M^{-1}(x,\l),\cr
&A(x,\l)=E^{-1}(x,\l)E_x(x,\l).\cr}$$  We get
$$\eqalign{-E^{-1}f^{-1}\d f E &=  E^{-1}\d E - M^{-1}\d M,\cr \d A &= [A, (E^{-1}f^{-1}\d f
E)_-]_+.\cr}$$ For $f=I$, we have $A=a\l$. So $\d A(x)$ is independent of $\l$ and lies in
$\cu_a^\perp$ for all $x$. This implies that 
$E^{-1}E_x=a\l + u$ for some $u:R\to\cu_a^\perp$. The fact that $u\in \cs(R, \cu_a^\perp)$
follows directly from (ii) and (iii).  Write 
$$u= M^{-1}M_x + M^{-1}[a\l, M].$$
Now $M(\cdot, \l)-M_{\infty}\in \cs(R^+)$ and $[M_\infty, a]=0$ imply that $u\n R^+\in
\cs(R^+)$. The corresponding argument gives $u\n R^-\in \cs(R^-)$. 

We first prove the theorem for  $f\in D^c_-$. Use Proposition \refhe{} to write 
$$f=p_d(I+p)v,$$ where $p_d$ is $a$-diagonal and holomorphic in $C_+$, $p$ is strictly
upper $a$-triangular, and $v$ is unitary.  We will be looking at $x\to \infty$. Examine the
formula for 
$\l\in C_+$:
\refeq[hf]$$(e^{-a\l x} p(\l) e^{a\l x})_{jk} = \cases{0, & if $a_j\geq a_k$,\cr
p_{jk}(\l)e^{-i(a_j-a_k)\l x},& if $a_j<a_k$.\cr}$$
Here $p_{jk}\n R$ lies in the Schwartz space if $a_j<a_k$. 

Use inverse Fourier transform to write  $$p_{jk}(r)=\int_{-\infty}^\infty \hat
p_{jk}(y)e^{iry}dy.$$ So $$p_{jk}(r) e^{-i(a_j-a_k)rx} =\int_{-\infty}^\infty \hat p_{jk}(y+
(a_j-a_k)x)e^{iry}dy.$$ The piece $\int_0^\infty \hat p_{jk}(y+(a_j-a_k)x)e^{iry}
dy$ is the boundary value of a holomorphic map in $C_+$, which can be written
$$\xi_{jk}(x,\l)= \int_0^\infty \hat p_{jk}(y+(a_j-a_k)x)e^{i\l y} dy.$$
So 
$$\eqalign{p_{jk}(r) e^{-i(a_j-a_k)rx} &= \xi_{jk}(r,x) + \D_{jk}(x,r), \qquad {\rm where\,}\cr
\D_{jk}(x,r) &=\int_{-\infty}^0 \hat p_{jk}(y+(a_j-a_k)x) e^{iry} dy.\cr}$$
It follows from Lemma \reffs{} that $\D\in \cs(R^+,H_1)$ and $\N {\p^m \xi\over \p
x^m(x,r)}\N \leq c_m$. 

Now write $$\eqalign{e_a^{-1}p_d(I+p)e_a &= p_de_a^{-1}(I+p)e_a =p_d(I +\xi + \D) \cr
&=p_d((I+\xi)(I+\D)-\xi \D) \cr &= p_d(I+\xi)(I+\D -(I+\xi)^{-1}\xi \D).\cr}$$
We claim that $D=\D-(I+\xi)^{-1}\xi\D\in \cs(R^+,H_1)$. Note that 
$$(I+\xi)^{-1}= I-\xi + \xi^2 -\xi^3 + \cdots +\xi^n$$ is a finite series since $\xi$ is
strictly upper $a$-triangular. The rules of multiplication of $\cs(R^+,H_1)$ 
by a smooth bounded function give the result that $D\in \cs(R^+,H_1)$. 

Let $(I+D)=(I+h)V$ be the Birkhoff decomposition. By Theorem \refft{}, $h$
and $V-I$ are in $\cs(R^+,H_1)$. 
So  $$\eqalign{e_a^{-1}f &= e_a^{-1}p_d(I+p)v =
e_a^{-1}p_d(I+p)e_a e_a^{-1}v\cr &=p_d(I+\xi)(I+D)e_a^{-1}v\cr
&=p_d(I+\xi)(I+h)V e_a^{-1}v.\cr}$$    By definition $M=p_d(I+\xi)(I+h)$, and 
$$M-p_d(I+\xi) = p_d(I+\xi)h.$$ Since $h\in \cs(R^+,H_1)$, and we have
uniform estimates on all derivatives of $p_d(I+p)$, $M-p_d(I+\xi) \in
\cs(R^+, H_1)$.   The same argument, in which a factorization $f=q\tilde v$ for $q$
lower $a$-triangular and
$\tilde v$ unitary, proves Schwartz space decay as $x\to -\infty$. 

To complete the proof, given $f\in D_-$, write $f=hg\in D_-^c\times G_-^m$. Write
$$h^{-1}e_{a,1}(x)= E_0(x)M_0^{-1}(x)\in G_+\times D_-.$$ By Theorem \refmi{}, we
factor 
$g^{-1}E_0(x)=E\tilde(x) g(x)\in G_+\times G_-^m$. Then
$$\eqalign{f^{-1}e_{a,1}(x) &= g^{-1}h^{-1}e_{a,1}(x) = g^{-1}E_0(x)M_0^{-1}(x)
=E(x)\tilde g(x) M_0^{-1} (x)\cr& = E(x)M(x).\cr}$$ By Theorem \refgm{} $\tilde g$
satisfies condition (ii) and (iii). But we just proved that $M_0$ satisfies (ii) and (iii),
so is
$M=\tilde g M_0^{-1}$. \qed

Note the convergence is actually uniform in the argument in Theorem \refhc{}. So
we have

\refclaim[nb] Theorem. As in Theorem \refhc{} let $f\in
D_-$,  $f^{-1}e_{a,1}(x)=E(x)M(x)^{-1}\in G_+\times D_-$, and $f_+=\lim_{s\searrow
0} f(r+is)$.  Factor $f_+=Pv = Q
\tilde v$, where $v, \tilde v$ are  unitary,  $P$ is upper $a$-triangular, $Q$ is lower
$a$-triangular,  and $P_d, Q_d$ is holomorphic in $C_+$. Then 
$$\eqalign{\lim_{x\to \infty} e^{ar x} M_+(x,r)e^{-arx}&= P(r),\cr
\lim_{x\to -\infty} e^{ar x} M_+(x,r)e^{-arx}&= Q(r).\cr}$$

\refclaim[go] Theorem.  Let $\Psi:D_-\to \cs_{1,a}$ be the map defined by $\Psi(f) =
E^{-1}E_x$, where $E$ is obtained from $f$ as in Theorem \refhc{}. Let $H_-$ denote the
subgroup of $f\in D_-$ such that $fa=af$.  Then
\item {(i)} $\cs^0_{1,a}= \Psi(D_-)$ is an open and dense subset of $\cs_{1,a}$,
\item {(ii)} $\Psi(f)=\Psi(g)$ if and only if there exist $h\in H_-$ such that $g=hf$,
\item {(iii)} $\cs^0_{1,a}$ is isomorphic to the homogeneous space $D_-/H_-$ of left cosets
of $H_-$ in $G_-$,
\item {(iv)} if $A=\Psi(f)$ and $M$ is as in Theorem \refhc{}, then the normalized
eigenfunction $m$ in Theorem \refgy{} of $A$ is $M_{-\infty}^{-1}M$.

\proof The first part (i) is a consequence of Theorem \refgy{}. Both (iii) and (iv)
follow from (ii). To prove (ii), recall if
$$\eqalign{f^{-1}e_{a,1}(x)&=E(x)M^{-1}(x)\in G_+\times D_-, \cr
g^{-1}e_{a,1}(x)&= E(x)N^{-1}(x)\in G_+\times D_-.\cr}$$ Then
$$M(x)N^{-1}(x)= e_{a,1}(x)^{-1} fg^{-1}e_{a,1}(x).$$ Suppose Im$\l >0$. Then the limit
of the right hand side is upper $a$-triangular when $x\to \infty$, and the limit is
lower
$a$-triangular when
$x\to -\infty$. So $MN^{-1}$ is both upper and lower $a$-triangular. Hence it is
$a$-diagonal, i.e.,  $MN^{-1}\in H_-$. So $fg^{-1}\in H_-$.

Conversely, if $g=hf$ for some $h\in H_-$ and
$f^{-1}e_{a,1}(x)=E(x)M(x)^{-1}\in G_+\times D_-$, then 
$$\eqalign{g^{-1}e_{a,1}(x)&= f^{-1}h^{-1}e_{a,1}(x)= f^{-1}e_{a,1}(x) h^{-1}\cr & =
E(x) (M(x)^{-1}h^{-1})\in G_+\times D_-.\cr}$$ So $\Psi(f)=\Psi(g)$.  \qed

In summary, we have shown that given $f\in D_-$, we can construct an $A\in
\cs_{1,a}$ such that $A=\Psi(f)$ by using various Birkhoff decomposition theorems
repeatedly. 

\refclaim[gq] Theorem. The natural right action of $D_-$ on the space $D_-/H_-$ of left
cosets induces a natural action $\ast$ of $D_-$ on $\cs^0_{1,a}$ via the
isomorphism $\bar\Psi$ from $D_-/H_-$ to $ \cs^0_{1,a}$. Equivalently, if
$A=\Psi(f)$ and
$g\in D_-$ then
$g\ast A= \Psi(fg^{-1})$. Moreover: 
\item {(i)} Let $g\in D_-$, $A\in \cs_{1,a}^0$, and $E$ the trivialization of $A$ normalized at
$x=0$, then we can factor $$gE(x)= \tilde E(x)\tilde g(x)\in G_+ \times D_-,$$ and $g\ast
A= \tilde E^{-1}\tilde E_x$.
\item {(ii)} If $g\in G_-^m$, then
$g\ast A=g\sharp A$, where $\sharp$ is the action of $G_-^m$ on $\cs_{1,a}$ defined in
Theorem \refmj{}. In other words, if $A=\Psi(f)$ and $g\in G_-^m$, then $g\sharp A = \Psi(fg)$. Or
equivalently, if $H_-f$ is the scattering coset of $A$ then $H_-fg$ is the scattering coset of
$g\sharp f$. 

\proof Given $f,g\in D_-$, we factor 
$$f^{-1}e_{a,1}(x)= E(x)M^{-1}(x), \qquad (fg)^{-1}e_{a,1}(x) = \tilde E(x) \tilde
M^{-1}(x).$$ Then
$g^{-1}E(x)= \tilde E(x)(\tilde M^{-1}(x) M(x))\in G_+\times D_-$. This defines the
action of $D_-$ on $\cs^0_{1,a}$, and it extends the action of $G_-^m$ on $\cs_{1,a}$
defined in Theorem
\refmj{}.\qed

\ss
\refpar[cn] Remark. If the scattering data of $A$ has $k$ poles counted with multiplicity,
then
$g_{z,\pi}\sharp A$ typically has $k+1$ poles, but it may have  $k$ or $k-1$ poles
for special choices of $z$ and $\pi$.  To see this,  let
$z\in C\setminus R$, and $\pi$ a projection such that $\pi a\not= a\pi$. If $\bar z$ is
not a pole of the scattering data of $A$ then $g_{z,\pi}\sharp A$ add one pole $\bar
z$ to the scattering data.  Let
$A=g_{\bar z,\pi}\sharp A_0$, where $A_0$ is the vacuum solution. Then  the scattering data of
$g_{z,\pi_1}\sharp A$
\item {(i)} has no poles if $\pi_1=\pi$,
\item {(ii)} has only one pole $z$ if $\pi_1$ and $\pi$ commute and $\pi+\pi_1\not= I$.

\bs

\newsection Poisson structure for the positive flows.\par

Let $H_+$ denote the subgroup of $G_+$ generated by $$\{e^{p(\l)}\n p(\l) \,{\rm is\,
a
\, polynomial\,}\, p(\l)a=ap(\l)\}.$$ In this section, we prove that the right dressing
action of $H_+$ on $D_-$ induces a Poisson group action of $H_+$ on $\cs_{1,a}^0$
and show that it generates the positive flows defined in section 5. We also study
the induced symplectic structure on the space of discrete scattering data
$G_-^m/(G_-^m\cap H_-)$, and the space of the continuous scattering data
$D_-^c/(D_-^c\cap H_-)$.

Set $e_{a,j,b} (x,t)(\l)= e^{a\l x + b\l^j t}$ and recall $e_{b,j}(x)(\l)=e^{b\l^j x}$.

\refclaim[me] Theorem. Let $a, b\in u(n)$ such that $[a,b]=0$. Then we can factor 
$$f^{-1} e_{a,j,b}(x,t) = E(x,t) M^{-1}(x,t) \in G_+\times D_-.$$ Moreover, $E$ and $M$
satisfy the following conditions: 
\item {(i)}$E^{-1}E_x =A$ is a solution of the $j$-th flow
defined by $b$, where $A(x,\l)= a\l + u(x)$.
\item {(ii)} $E^{-1} E_t =B$, where $$B(\cdot,\l)= b\l^j + Q_{b,1}(u) \l^{j-1} + \cdots
+ Q_{b,j}(u) = (M^{-1}b\l^jM)_+.$$

\proof Since $[a,b]=0$, $\exp(a\l x+ b\l^j t)= e^{a\l x} e^{b\l^j t}$. Use Theorem
\refmc{} to factor 
$$f^{-1}e_{a,j,b}(x,t)= f^{-1}e_{a,1}(x) e_{b,j}(t) = E_0(x)M_0^{-1}(x) e_{b,j}(t).$$
Use Theorem \refmc{} again to factor 
$$M_0^{-1}(x) e_{b,j}(t)= E_1(x,t) M(x,t)\in G_+\times D_-.$$ 

The variational form of $f^{-1}e_{a,j,b}=EM^{-1}$ implies
$$\cases{E^{-1}E_x= a\l + u,&\cr E^{-1}E_t=b\l^j + q_1\l^{j-1} + \cdots + q_j.&\cr}$$
So 
\refeq[pc]$$\left[{\partial\over \partial x} + a\l + u, \,\, {\partial \over \partial t} + 
b\l^j + q_1\l^{j-1} + \cdots + q_j\right] = 0.$$  Compare coefficient of $\l^i$ in
equation \refpc{} to get 
$$\cases{(q_i)_x+[u, q_i]=[q_{i+1},a],& if $0\leq i<j$, \cr u_t=(q_j)_x+ [u,q_{j+1}].&\cr}$$
This is the same system as \refdb{} defining the $Q_{b,is}'s$. Hence  $q_i=Q_{b,i}$.
\qed

\refclaim[hg] Corollary. The dressing action $\natural$ of $H_+$ on $D_-$ on the
right is well-defined and $H_-$ is fixed under this action. Hence an action
$\natural$ of $H_+$ on $D_-/H_-$ is defined, which leads to an action on
$\cs_{1,a}^0$. In fact, this action is defined as follows: Write
$A=\Psi(f)$,
$f^{-1}e_{a,1}(x)=E(x)M(x)^{-1}$. For $h\in H_+$, we factor
$$hM(x)=\tilde M(x)\tilde h(x)\in D_-\times G_+$$ to get $$h\natural A= (e^{a\l
x}\tilde M)^{-1}(e^{a\l x}\tilde M)_x.$$ 

\refclaim[pd] Corollary. If $A_0=\Psi(f_0)$, then $A(t)=\Psi(e_{b,j}(t)\natural f_0)$
is the solution of the $j$-th flow on $\cs_{1,a}$ defined by $b$ with $A(0)=A_0$.

\refclaim[aw] Corollary. Let $a_1, \cdots, a_{n}$ be a basis of the space $\ct$ of
diagonal matrices in  $u(n)$, 
$f\in D_-$, and $e_{a_1,\cdots, a_n}(x_1, \cdots, x_n)(\l)=\exp(\sum_{j=1}^n
a_jx_j\l)$.  Factor
$$f^{-1}e_{a_1, \cdots, a_n}(x)= E(x)M(x)\in G_+\times D_-.$$ Then there exists
$v:R^n\to\ct^\perp$ such that
\item {(i)} $E^{-1}E_{x_j}=a_j \l + [a_j,v]$ for all $1\leq j\leq n$,
\item {(ii)} $v$ is a solution of equation 
\refeq[ee]$$\left[a_i,{\p v\over \p x_i}\right] -\left[a_j,{\p v\over \p x_j}\right]
=[[a_i,v],[a_j,v]].$$

\refpar[] Remark. Equation \refee{} is the $n$-dimensional system associated to
$U(n)$ constructed in the paper of the first author [Te2].

\refclaim[Ai] Theorem. The action $\natural$ of $H_+$ on $\cs^0_{1,a}$ is Poisson. Moreover, the
map $\mu:\cs^0_{1,a}\to H_-=H_+^*$ defined by $\mu(A)(\l)= M_{-\infty}^{-1}M_\infty$
 is a moment map, where $A=\Psi(f)$, $f^{-1}e_{a,1}(x)=E(x)M(x)^{-1}\in G_+\times
D_-$ and $M_{\pm\infty}(\l)=\lim_{x\to\pm\infty} M(x,\l)$.

\proof Suppose $A=\Psi(f)$, i.e., $$\cases{f^{-1}e_{a,1}(x) = E(x)M^{-1}(x)\in
G_+\times D_-,&\cr  A=E^{-1}E_x= (e^{a\l x}M)^{-1}(e^{a\l x}M)_x.&\cr}$$ The
second equation implies
\refeq[fc]$$M^{-1}M_x + \l M^{-1}aM = A.$$ 
Set $\eta= M^{-1}\d M$, $B=\d A$ and $\psi=e^{a\l x}M$. 
Compute the variation directly from  equation \reffc{} to derive 
$$\eqalign{&\eta_x+ [A,\eta]= B, \qquad \lim_{x\to -\infty}\eta = 0,\cr &\eta(x)=
\psi(x)^{-1}
\int_{-\infty}^x(\psi B\psi^{-1})dy\, \psi(x).\cr}$$ For $\xi_+\in \ch_+$, since
$[\xi_+,a]=0$, we have $M^{-1}\xi_+ M= \psi^{-1}\xi_+\psi$. 
$$\eqalign{&<d\mu_A(B)\mu(A)^{-1},\quad\xi_+>\cr &= \lim_{x\to \infty}<
M(x)\psi^{-1}(x)\int_{-\infty}^x (\psi B\psi^{-1})dy \psi(x) M^{-1}(x), \,\,
\xi_+>,\cr  &= \lim_{x\to \infty} < e^{-a\l x} \int_{-\infty}^x (\psi B\psi^{-1})dy
e^{a\l x}, \quad
\xi_+> \cr &= \lim_{x\to \infty} <  \int_{-\infty}^x (\psi B\psi^{-1})dy, \quad
 e^{a\l x}\xi_+e^{-a\l x}> \cr &= \lim_{x\to \infty} <  \int_{-\infty}^x (\psi
B\psi^{-1})dy, \quad \xi_+>= \lim_{x\to \infty} \int_{-\infty}^x < \psi B
\psi^{-1},\quad \xi_+>dy \cr &=\lim_{x\to \infty} \int_{-\infty}^x <  B,
\quad \psi^{-1}e^{-a\l y}\xi_+e^{a\l y}\psi> dy  \cr  &= << B, \quad 
(M^{-1}\xi_+M)_->>.\cr}$$ The rest of the proof goes exactly the same as for Theorem
\refBb{}. \qed 

\ms
\refpar[cu] Remark. Let $a=\diag(ia_1, \cdots, ia_n)$, and $a_1 < \cdots < a_n$. Then
$\cu_a$ is the set of all diagonal matrices in $u(n)$, $\cu_a^\perp$ is the set of all 
matrices
$u\in u(n)$ such that $u_{ii}=0$ for all $1\leq i\leq n$. So  $H_+$ is abelian and the action of
$H_+$ on $\cs_{1,a}$ is in fact symplectic.
\ms

The following theorem was proved by Flaschke, Newell and Ratiu [FNR1, 2] for $n=2$ and by
one of us [Te2] for general $n$:

\refclaim[co] Theorem ([Te2]). The  Hamiltonian function on $\cs_{1,a}$
corresponding to the $j$-th flow defined by $a$ is:
\refeq[ge]$$F_{a,j}(u) = -{1\over j+1} \int_{-\infty}^\infty (Q_{a,j+2}, a) dx,$$ 
i.e., $\K F_{b,j}(u)= Q_{b,j+1}^\perp(u)$.

\refpar[hi] Remark. Let $b, c\in \cu_a$, and $\xi_{b,j}$ and $ \xi_{c,k}$ denote infinitesimal
vector field for the $H_+$-action on $\cs_{1,a}^0$ corresponding to $b\l^j$ and $c\l^k$
respectively. Then the bracket $[\xi_{b,j},\xi_{c,k}]$ is equal to the infinitesimal
vector field corresponding to $[b,c]\l^{k+j}$. Unless $[b,c]=0$, these two flows do not
commute. 

\ms 

\refpar[ne] Remark. If we replace the group $SU(n)$ by a
simple compact Lie group, then what we have discussed still holds if appropriate algebraic
conditions are prescribed. 

\ms

In the end of this section,  we will study the pull back of the
symplectic structure $w$ on $\cs_{1,a}$ to $D_-/H_-$ via the isomorphism $\Psi$.  
Note that $\Psi(G_-^m)$ ($\Psi( D_-^c)$ resp.) is the space of $A$'s with only discrete
(continuous resp.) scattering data. 

We have been using the base point $x=0$, i.e.,
$f(\l)=M(0,\l)$. But there is nothing special about $x=0$. In the following, we choose a base
point $y$ and let $y\to -\infty$. The expression with base point
$0$, and with $y$, differ by a term which cancels out when we evaluate integrals at the end
points. 

We only deal with the symplectic structure on $\Psi(D_-^c)$. However, this set is
pretty large. For example, Beals and Coifman and later Zhou show the
following:

\refclaim[na] Theorem ([BDZ]). Let $B_1$ denote the unit ball in $\cs_{1,a}$ with respect to
the
$L^1$-norm, i.e., $B_1$ is  set of all $A=a\l+u\in \cs_{1,a}$ such that  $\int_{-\infty}^\infty
\N u\N dr <1$. Then $B_1\subset  \Psi(D_-^c)$. 

Let $S$ denote the scattering transform that maps $A\in \cs_{1,a}'$ to
its scattering data  $S$ (defined in section 7).   The restriction of the
symplectic form $w$ on $\cs_{1,a}$ to  $S(B_1)$ was computed by Beals and Sattinger [BS]. We
will compute the restriction of $w$ to $\Psi(D_-^c)$ in terms of variations in $D_-^c$ below.  
Let 
$$<f(r),g(r)>=\int_{-\infty}^\infty
\Im(\tr(f(r)g(r))dr.$$ By the same computation as in Theorem
\refBb{}, the Poisson bracket on
$B_1\subset \cs_{1,a}$ is
$$\eqalign{\{\d_1 A, \d_2 A\}&= \lim_{x\to \infty}  <(E^{-1}(x)f^{-1}\d_1 f E(x))_-,
E(x)^{-1}f^{-1}\d_2 f E(x)> \cr &\quad  - \lim_{y\to -\infty}  <(E^{-1}(y)f^{-1}\d_1 f
E(y))_-, E(y)^{-1}f^{-1}\d_2 f E(y)> \cr  &= \lim_{x\to \infty} <E^{-1}(x)\d_1 E(x),
M^{-1}\d_2 M(x)> \cr &\quad  -
\lim_{y\to -\infty} <E^{-1}(y)\d_1 E(y), M^{-1}\d_2 M(y)>.\cr}$$ 
Now, let the vacuum be based at $y$, i.e., factor
$$f^{-1}(\l)e^{a\l(x-y)}=E(x,y,\l)M^{-1}(x,y,\l).$$ Hence
$f(\l)=M(y,y,\l)$ and $\d E(y,y,\l)=0$. The
$y$-term in the above description is now zero, and we have
$$\eqalign{&=\lim_{x\to \infty, y\to -\infty} -<M(x)^{-1}e^{-a\l(x-y)}\d_1 f f^{-1}
e^{a\l(x-y)}M(x), \,\, M^{-1}(x)\d M(x)>\cr
&=\lim_{x\to \infty, y\to -\infty} -<e^{-a\l(x-y)}\d_1 f f^{-1} e^{a\l(x-y)}, \,\,\d_2 M(x)
M^{-1}(x)>\cr
&= \lim_{x\to \infty, y\to -\infty} -<e^{a\l y} \d_1 M(y) M^{-1}(y) e^{-a\l y}, \,\, e^{a\l x}
\d_2 M(x) M^{-1}(x) e^{-a\l x}>.\cr}$$

Now by Theorem \refnb{}, we get

\refclaim[nc] Theorem. The Poisson structure on the unit ball $B_1$ in $\cs_{1,a}$ with respect to
the $L^1$-norm, written in terms of variations in $D_-^c$, is 
$$\{\d_1 A, \d_2 A\} = \int_{-\infty}^\infty <\d_1 P
P^{-1}, \d_2 Q Q^{-1}>,$$ where $A=\Psi(f)$, $f_+=Pv=Q\tilde v$ is the factorization of $f_+$
into upper $a$-triangular and lower $a$-triangular times unitary and $P_d, Q_d$
are holomorphic as given in Proposition \refhe{}.

Next we study the pull back the symplectic form 
$$w(q_1,q_2)=\int_{-\infty}^\infty \tr(\ad(a)^{-1}(q_1)(q_2))dx$$  on $\cs_{1,a}$ to the space 
 $\Psi(G_-^m)$. This space has many complicated algebraic components. For example the space
 of all $A\in \cs_{1,a}$ whose scattering data have only one pole (or equivalently, $A=\Psi(g)$,
where $g$ is a simple element)  can be parametrized by 
$$ \bigcup_{k=1}^n C_+\times \{V\in Gr(k,C)\n a(V)\not\subset V\}.$$ However, the 
space of $A$ whose scattering data has only two poles immediately
becomes complicated as the factorization of $g\in G_-^m$ as product of simple
elements is not unique. The following Proposition gives the restriction of
$w$ to the simplest component of $\Psi(G_-^m)$. We  believe that
the restriction of
$w$ to each algebraic component should be symplectic, but we have not yet found an efficient
way to  compute the general case. 

\refclaim[nd] Proposition. Let $a=\diag(-i, i, \cdots, i)$. Then:
\item {(i)}  The space of all
$A=\Psi(g_{z,\pi})$, where $\pi$ is the projection onto a one dimensional subspace $Cv$, is
isomorphic to
$$N=C_+\times (C^{n-1}\setminus {0})=\{(z,v)\n z\in C\setminus R, v=(v_2, \cdots,
v_n)\not=0\}.$$ 
\item {(ii)} The pull back of the symplectic form
$w$ to $N$  is $$2\,{\rm Re\/}\left(d\bar z\wedge
\partial \log(\n v\n^2) +(z-\bar z) \partial
\bar\partial
\log (\n v\n^2)\right),$$ where $\n v\n = \sum_{j=2}^n \n v_j\n^2$. 

\proof Let $\pi$ denote the projection of $C^n$ onto the one dimensional subspace
spanned by $(1,v)$, where $v=(v_2, \cdots, v_n)\in C^{n-1}$.  By Theorem \refmi{} (vi) and
formula \refhr{},  
$g_{z,\pi}\sharp 0= (u_{ij})$, where $(u_{ij})\in u(n)$, $u_{ij}=0$ if $2\leq i, j\leq n$ and 
$$u_{1j}(x)= {2i(z-\bar z) \bar v_j e^{i(z+\bar z)x}\over e^{-i(z-\bar z) x} + e^{i(z-\bar
z)}(\n v_2\n^2 + \cdots + \n v_n\n^2)}.$$
Then the proposition follows from at least two separate computations, neither of which is
very illuminating. We  hope to provide the more general results in a future paper. 
\qed

\ss
\refpar[nf] Remark.  Fix $z\in C_+$. Then the restriction of the symplectic form to the
subset $\{(z,v)\n v\in C^{n-1}\setminus \{0\}, \N v\N = 1\}$ of $M$ in the theorem above gives
the standard symplectic structure of $CP^{n-2}$. 

\bs

\newsection Symplectic structures for the restricted case.\par

Most of the interesting applications in geometry come from restrictions of the full flow
equation to a smaller phase space satisfying additional algebraic conditions. This leads to
a serious problem, not with the flows and the scattering cosets, but with the symplectic
structure. Generally the original symplectic structure we have used to this point vanishes
on the restricted submanifolds.  In this section, we describe a typical restrictions and the
construction of the hierarchy of symplectic structures. We give an outline of the theory,
and explain how it can be applied.  The details of this construction for involutions
appear in [Te2]. 

Let $U$ be a simple Lie group, $< , >$ a non-degenerate, ad-invariant bilinear form on the Lie
algebra $\cu$, and $\s$ an order $m$ automorphism of $U$. For simplicity, we denote the Lie
algebra automorphism
$d\s_e$ on $\cu$ again by $\s$.  Fix a primitive
$m$-th root of unity $\a$. Suppose $\s$ has an eigendecomposition on $\cu$:
$$\cu=\sum_{j=0}^m \cu_j,$$ where $\cu_j$ is the eigenspace of $\s$ on $\cu$ with
eigenvalue
$\a^j$. Then $$\cases{[\cu_j,\cu_k]\subset\cu_{j+k}, & for all $j,k$,\cr <\cu_j , \cu_k>=0, & if
$j\not=k$.\cr} $$ Here we use the convention
$\cu_j=\cu_k$ if $j\equiv k$ (mod $m$). 

Fix an element $a\in \cu_1$. Let $\cu_a$ denote the centralizer of $a$, and
$\cu_a^\perp$ the orthogonal complement of $\cu_a$ in $\cu$. Consider a subspace
of $\cs_{1,a}$:
$$\cs_{1,a}^\s=\{A=a\l + u\in \cs_{1,a} \n u\in \cu_0\cap \cu_a^\perp\} \quad \subset\,\,
\cs_{1,a}.$$
Then  $A\in \cs^\s_{1,a}$ satisfies the reality condition
\refeq[gf]$$\s(A(\a^{-1}\l))=A(\l).$$ Hence the trivialization $E$ of $A\in \cs_{1,a}$
normalized at the origin satisfies the condition
$\s(E(\a^{-1} \l))= E(\l)$. Then  

\refclaim[nh] Proposition. $\cs_{1,a}^\s$ is invariant under the action $\sharp$ of
$G_-^{m,\s}$ on $\cs_{1,a}$ defined in section 6, where 
 $G_-^{m,\s}$ is the subgroup of $g\in G_-^m$ such that $\s(g(\a^{-1}\l))=g(\l)$.
 
The trivialization 
$M$ of $A\in \cs_{1,a}^\s$ normalized at infinity also satisfies the same reality
condition as $E$.  So for $b\in \cu_1\cap \cu_a$, we have
$$\s(M^{-1}(\a^{-1}\l) b M(\a^{-1}\l))=\a M^{-1}(\l)bM(\l).$$ 
 Since $$M^{-1}bM \sim \sum_{j=0}^\infty Q_{b,j}\l^{-j},$$ we get
$$Q_{b, j}\in \cu_{-j+1}, \qquad {\rm for\,\, all\,\,} j\geq 0.$$ In particular, $[Q_{j+1},a]\in
\cu_{-j+1}$. So $[Q_{b,j+1}(u),a]$ is normal to $\cs_{1,a}^\s$ if $j\not\equiv 1$ (mod $m$), and
is tangent to $\cs_{1,a}^\s$ if $j\equiv 1$ (mod $m$).  And we have 

\refclaim[] Proposition. The $j$-th flow preserves $\cs_{1,a}^\s$ for all $j$ and any
order $m$ automorphism $\s$. If
$j\not\equiv 1$ (mod $m$), then the flows are identically constant.

However, the symplectic form 
$$w(\d_1u,\d_2u)=\int_{-\infty}^\infty \tr(-\ad(a)^{-1}(\d_1 u), \d_2 u) dx$$
vanishes on $\cs_{1,a}^\s$.

The sequence of symplectic structures constructed by Terng can be described
using a sequence of coadjoint orbits, which arise using a shift in the bi-linear form
$< , >$ on the loop algebra $\cg$. 

For $k\leq -1$, let $M_k$ denote the coadjoint $C(R,G_-)$-orbit at
$\left({d\over dx} + a\l\right)\nu_k^{-1}$ in $C(R,\cg_+)$, where $\nu_k(\l)=\l^{k+1}$.  Set 
$$\cs_{1,a,k}=(M_k\nu_k)\cap \cs_{1,a}.$$ Then
$\d u$ lies in the tangent space of $\cs_{1,a,k}$ at ${d\over dx}+a\l + u$ if and only if 
\refeq[fv]$$\d u(x)= \left(\left[\xi_-(x), {d\over d x} +
A(x)\right]\nu_k^{-1}\right)_+\nu_k,$$ where formally
$\xi_-(x)\in\cg_-$. Here $(\cdot )_+$ is the
projection into $\cg_+$, and the construction is entirely algebraic. For $A=a\l + u$,  write 
$$\xi_-(\d u) = \xi_{-1}(\d u)\l^{-1} + \xi_{-2}(\d u) \l^{-2} + \cdots.$$ Then equation
\reffv{} gives
$$\eqalign{&[\xi_{-1},a] = \d u, \cr
&[\xi_j,a]=\left[{d\over dx} + u, \xi_{j+1}\right], \qquad k\leq j\leq -1.\cr}$$ This
gives a recipe to compute $\xi_-(\d u)$ explicitly via a mixed
integro-differential operation:
$$\eqalign{\xi_{-1}^\perp(\d u) &= J_a^{-1}(\d u),\cr
\xi_j^\perp(\d u) &=(J_a^{-1}P_u)^{-j-1}J_a^{-1}(\d u), \cr}$$ where
$$\eqalign{J_a(v)&=[v,a],\cr
P_u(v)&= v_x +
[u,v]^\perp -[u,\eta_u(v)], \cr
\eta_u(v)&=\int_{-\infty}^x[u(y),v(y)]^d dy.
\cr}$$
Here $v^\perp$ and $v^d$ denote the projection onto $a$-off diagonal ($\cu_a^\perp$) and
$a$-diagonal ($\cu_a$) respectively.  Set $$J_k=J_a(J_a^{-1}P_u)^{k+1}.$$  The natural
shifted symplectic structure is given by
$$\eqalign{w_k(\d_1 u, \d_2 u)& =\int_{-\infty}^\infty <{d\over dx}+A,
\nu_k^{-1}[\xi_-(\d_1u),\xi_-(\d_2u)]> dx\cr
&=\int_{-\infty}^\infty \tr\left(\left({d\over
dx}+a\l+u\right)([\xi_-(\d_1u),\xi_-(\d_2u)])\right)_k dx,\cr&=\int_{-\infty}^\infty
\tr\left(\left(\d_1u\right)\xi^\perp_k\left(\d_2u\right) \right)dx,\cr
&=\int_{-\infty}^\infty \tr\left(\left(\d_1u\right)J_k^{-1}\left(\d_2 u\right)\right)
dx,\cr}$$ where $(\cdot)_k$ denote the coefficient of $\l^k$ in $(\cdot)$. 
In particular, the first two in the series are:
$$\eqalign{w_{-1}(\d_1u,\d_2u)&=w(\d_1u,\d_2u)=\int_{-\infty}^\infty
\tr((-\ad(a)^{-1}(\d_1u))\d_2u)dx,\cr
w_{-2}(\d_1u,\d_2u) &=\int_{-\infty}^\infty
\tr((\d_1u)(J_{-2})^{-1}_u(\d_2u))dx\cr &=\int_{-\infty}^\infty \tr((\d_1
u)J_a^{-1}P_uJ_a^{-1}(\d_2u))dx.\cr}$$
The natural coadjoint orbits require the relevant terms of $\xi_-$ to lie in the
Schwartz class. So the tangent space of the smaller submanifold
$\cs_{1,a,k}=(M_k\nu_k)\cap \cs_{1,a}$ at 
$u$ is $$\left\{\d u\n \xi_{-j}(\d u)(\infty)=0, 1\leq j\leq -k\right\}.$$  Hence $\cs_{1,a,k}$ is
a finite codimension submanifold of $\cs_{1,a}$ and the formulas we write down for 
$w_k$ are skew symmetric.  

For $k\geq 0$, let $M_k$ denote the coadjoint $C(R,G_+)$-orbit at
$\left({d\over dx} + a\l\right)\nu_k^{-1}$ in $C(R,\cg_-)$, and 
$\cs_{1,a,k}=(M_k\nu_k)\cap \cs_{1,a}$, where $\nu_k(\l)=\l^{k+1}$. Then
$\d u$ lies in the tangent space of $\cs_{1,a,k}$ at ${d\over dx}+A$ if and only if 
\refeq[fw]$$\d u(x)= \left(\left[\xi_+(x), {d\over d x} +
A(x)\right]\nu_k^{-1}\right)_-\nu_k,$$ where formally
$\xi_+(x)\in\cg_+$. Here $(  )_-$ is the
projection into $\cg_-$. For $A=a\l + u$,  write 
$$\xi_+(\d u) = \xi_0(\d u) + \xi_1(\d u) \l + \cdots.$$ Then equation
\reffw{} gives
$$\eqalign{&\left[{d\over dx}+u, \xi_0\right] = \d u, \cr
&\left[{d\over dx}+ u, \xi_j\right]=\left[\xi_{j-1},a\right], \qquad 1\leq j\leq
k.\cr}$$ Hence
$$\xi_j^\perp(\d u)= (P_u^{-1}J_a)^jP_u^{-1}(\d u) = J_j^{-1}(\d u).$$
 The natural
shifted symplectic structure is given by
$$w_k(\d_1 u, \d_2 u)=\int_{-\infty}^\infty \tr((\d_1u)J_k^{-1}(\d_2 u)) dx.$$
In particular, 
$$w_0(\d_1 u, \d_2 u)=\int_{-\infty}^\infty \tr((\d_1
u)P_u^{-1}(\d_2u))dx.$$

If $a\in \cu_1$, then $J_a=-\ad(a)$ maps $\cu_j$ to $\cu_{j+1}$. This implies that
$J_k$ maps
$\cu_j$ to $\cu_{j-k}$. Thus we obtain:

\refclaim[fva] Proposition ([Te2]). $w_k$ is a symplectic structure on
$\cs_{1,a,k}$. Moreover, 
\item {(i)} $w_k=0$ on $\cs_{1,a,k}\cap \cs_{1,a}^\s$ if $k\not\equiv 0$ (mod $m$),
\item {(ii)} $w_k$ is non-degenerate on $\cs_{1,a,k}\cap \cs_{1,a}^\s$ if $k\equiv 0$ (mod
$m$).

Recall that $F_{b,j}$ defined by formula \refge{} is the Hamiltonian for the $j$-th flow on
$\cs_{1,a}$ defined by $b$ with respect to the symplectic form $w_{-1}$, and $\K
F_{b,j}=Q_{b,j+1}^\perp$.  Since
$P_u(Q_{b,j}^\perp)= [Q_{b,j+1},a]$,  we get 

\refclaim[nl] Theorem ([Te2]). If $a$ is regular, then
\item {(i)} $J_r (\K F_{b, j})= [Q_{b, j+r+2},a]$,
\item {(ii)}  the Hamiltonian flow corresponding to $F_{b,j}$ on
$(\cs_{1,a}, w_r)$ is the $(j+r+1)$-th flow defined by $b$.

\refpar[ct] Examples.  \ss

{\bf Example 1}. \hskip 4pt Let $\s$ denote the involution
$\s(y)=- y^t$ of $SU(n)$, and $a=\diag(i,-i, \cdots, -i)$. Then $\cs^\s_{1,a,0}$ is the set of all
$A=a\l+u$ with $$u=\pmatrix{0&v\cr -v^t&0},$$ where $v:R\to \cm_{1\times
(n-1)}$ is a decay map from $R$ to the space $\cm_{1\times (n-1)}$ of real $1\times (n-1)$
matrices.   The even flows vanishes on
$\cs^\s_{1,a,0}$, and the odd flows are extensions of the usual hierarchy of flows for
the modified KdV. The third flow written in terms of $v:R\to \cm_{1\times (n-1)}$  is
the matrix modified KdV equation:
$$v_t=-{1\over 4}(v_{xxx} + 3(v_xv^tv+vv^tv_x).$$
(When $n=2$, $v$ is a scalar function and the above equation is the classic modified
KdV equation.) The $2$-form $w_0$ gives the appropriate non-degenerate symplectic
structure for the  matrix modified KdV equation and the hierarchy of odd flows. \ss

{\bf Example 2}. \hskip 4pt  It seems appropriate to mention the relation of the
restriction to the sine-Gordon equation. The sine-Gordon equation is 
written in space time coordinates $(\tau,y)$ as 
$${\p^2 q\over \p \tau^2}-{\p^2 q\over \p y^2} + \sin q = 0,$$ or $$2q_{xt}=\sin q$$ in
characteristic coordinates. The Lax pair is best written in characteristic coordinates:
$$\left[{\p\over \p x} +a\l + u, {\p\over \p t}+ \l^{-1} B\right]=0,$$
where $$a=\pmatrix{i&0\cr 0 & -i\cr}, \quad u=\pmatrix{0&q_x\cr -q_x&0}, \quad 
B=\pmatrix{\cos q& \sin q\cr -\sin q & \cos q\cr}.$$
 The restriction is the
same as for the modified KdV. The natural Cauchy problem is in space time coordinates
$(\tau,y)$, but the scattering theory has been developed for characteristic coordinates.
However, the classical B\"acklund transformations work well with either choice of
coordinates , and preserve whatever decay conditions have been described in
either coordinate systems. \ss

 {\bf Example 3}. \hskip 4pt We obtain the Kupershmidt and Wilson equation ([KW]) in
terms of a restriction by an order $n$ automorphism of $sl(n)$.  Let $\a=e^{2\pi i/m}$, and
$p\in SL(n)$  the matrix representing the cyclic permutation $(12 \cdots n)$, i.e.,
$p(e_i)=e_{i+1}$ (here we use the convention that $e_i=e_j$ if $i\equiv j$ (mod
$n$)).  Let $\s:sl(n)\to sl(n)$ be the order $n$ automorphism defined  by
$\s(y)=p^{-1}yp$. Then $X\in \cu_j$ if and only if
$\s(X) =\a^j X$.  Let $$a=\diag(1, \a, \cdots, \a^{n-1}) \in \cu_1.$$ 
Note that $\cu_a^\perp$ is the space of all matrices $X\in sl(n)$ such that $X_{ii}=0$ for all
$i=1, \cdots, n$. So
$$\cs_{1,a}^\s=\cases{\left\{\pmatrix{1&0\cr 0&-1\cr} \l + \pmatrix{0&v\cr v&0\cr}\,\bigg|
\, v\in
\cs(R, C)\right\},& if $n=2$,\cr
\left\{\pmatrix{1&0&0\cr 0&\a&0\cr
0&0&\a^2}\l +\pmatrix{0&v_1&v_2\cr v_2&0&v_1\cr v_1&v_2&0}\bigg| v_1, v_2\in
\cs(R,C)\right\},&if $n=3$.\cr}$$ In general, $A=a\l + u\in\cs_{1,a}^\s$ is determined by
$(n-1)$ functions (the first row of $u$). 
By Propositions \reffva{},  $\{w_{rn}\n r\in Z\}$ is a sequence of 
symplectic forms on $\cs_{1,a}^\s$. The $(n+1)$-th flow is the Kupershmidt-Wilson
equation. By Theorem \refnl{} it satisfies the Lenard relation:
$$u_t = [Q_{a,n+2},a] = J_0(\K F_{b,n})= J_n(\K F_{b,0}).$$
When $n=2$,  the third flow on $\cs_{1,a}^\s$ defined by $a$ gives the modified KdV
equation:
\refeq[ea]$$v_t={1\over 4}(v_{xxx}-6v^2v_x),$$ and all the odd flows are the hierarchy of
commuting flows of the modified KdV equation.  For $n>2$, this gives another
generalization of modified KdV equation. 

\bs

\newsection B\"acklund transformations for $j$-th flows.\par

This section contains a brief outline of ideas and results in a forthcoming paper
[TU1]. The classical B\"acklund transformations are originally geometric
constructions by which a two parameters family of constant  Gaussian curvature
$-1$ surfaces is obtained from a single surface of Gaussian curvature $-1$. This is
accomplished by solving two ordinary differential equations with a parameter $s$.
The second parameter is the initial data. Since surfaces of Gaussian curvature $-1$
are classically known to be equivalent to local solutions of the sine-Gordon
equation ([Da1], [Ei])
$$q_{xt}=\sin q$$ this provides a method of deriving new solutions of a partial differential
equation from a given solution via the solution of ordinary differential equations. Most of
the known ``integrable systems'' possess transformations of this type, which are
sometimes called {\it Darboux transformations\/}.  Ribaucour and Lie
transformations are other classical transformations that generate new solutions
from a given one. 

The action of the rational loop group we constructed in section 6 
can be extended to an
action which transforms solutions of the $j$-th flow equation. In this section we
describe very briefly the results of a future paper [TU1], which will construct an
action of the semi-direct product of
$R^*\sdp G_-^m$ on the solution space of the $j$-th flow. The construction of this
loop group action is motivated by the construction given by the second author in
[U1] for harmonic maps.  We will see
\item {(1)} the action of a simple element $g_{z,\pi}$ corresponds to a B\"acklund
transformation,
\item {(2)} the action of $R^*$ corresponds to the Lie transformations,
\item {(3)} the Bianchi permutability formula arises from the various ways of factoring
quadratic elements in the rational loop group into simple elements,
\item {(4)} the B\"acklund transformations arise from ordinary differential equations if one
solution is known,
\item {(5)} once given the trivialization of the Lax pair corresponding to a given
solution, the B\"acklund transformations become algebraic. 
\ms

Since the sine-Gordon equation arises as part of the algebraic structure (the $-1$-flow for
$su(2)$ with an involution constraint), we can check that we are generalizing the classical
theory. The choice of group structure depends on the choice of the base point (just as the
scattering theory depends on the choice of a vacuum, or the choice of $0\in R$). Hence the
group structure was not apparent to the classical geometers. 

One of the most interesting observations is that appropriate choices of poles for the
rational loop yield time periodic solutions. This yields an interesting insight into
the construction of time-periodic solutions (or the classical breathers) to the
sine-Gordon equation as explained in Darboux ([Da1]). For recent developments
concerning breathers of the sine-Gordon equation see [BMW], [De], [SS]. There are
no simple factors in the rational loop group corresponding to the placement of
poles for time periodic solutions. However, there are quadratic elements, whose
simple factors do not satisfy the algebraic constraints to preserve sine-Gordon,
but which nevertheless generate the well-known breathers (one way to think of
them is as the product of two complex conjugate B\"acklund transformations). The
product of these quadratic factors generate arbitrarily complicated time periodic
solutions. 

The classical theory of B\"acklund transformations is based on ordinary differential
equations.

\refclaim[ja] Theorem ([Ei]).  Suppose $q$ is a solution of the sine-Gordon equation,
and $s\not=0$ is a real number. Then the following first order system is solvable
for $q^*$:
\refeq[jb]$$\cases {(q^*-q)_x= 4s \sin ({q^*+q\over 2})\cr 
(q^*+q)_t = {1\over s} \sin ({q^*-q\over 2}).\cr}$$
Moreover,  $q^*$ is again a solution of the sine-Gordon equation.

\refpar[jc] Definition. If $q$ is a solution of the sine-Gordon equation, then given any $c_o\in
R$ there is a  unique solution $q^*$ for  equation \refjb{}
such that $q^*(0,0)= c_o$.  Then $B_{s,c_o}(q)=q^*$ is a transformation on the space of solutions
of the sine-Gordon equation, which will be called a {\it B\"acklund transformation\/}
for the sine-Gordon equation.  

\refclaim[jd] Proposition ([Ei]). Define $L_s(q)(x,t)= q(sx,s^{-1}t)$. Then 
$q$ is a solution of the sine-Gordon equation if and only if $L_s(q)$ is a solution of
the sine-Gordon equation.  ($L_s$ is called a {\it Lie transformation\/}). 

\refclaim[ko] Proposition ([Ei]). B\"acklund transformations and Lie transformations of the
sine-Gordon equation are related by the following formula:
$$B_{s,c_o}= L_s^{-1}B_{1,c_o}L_s.$$

There is also a Bianchi permutability theorem for surfaces with Gaussian curvature $-1$ in
$R^3$, which gives the following analytical formula for the sine-Gordon equation:

\refclaim[je] Theorem ([Ei]). Suppose $q_0$ is a solution of the sine-Gordon
equation,  $s_1^2\not= s_2^2$, and $s_1s_2\not=0$. Let
$q_i=B_{s_i,c_i}(q_0)$ for $i=1,2$. Then there exist $d_1,d_2\in R$,
which can be constructed algebraically,  such that 
\item {(1)} $B_{s_1,d_1}B_{s_2,c_2}= B_{s_2, d_2}B_{s_1,c_1}$,
\item {(2)} let $q_3= B_{s_1,d_1}B_{s_2,c_2}(q_0)$, then 
\refeq[jf]$$\tan {q_3-q_0\over 4}={s_1+s_2\over s_1-s_2}\tan {q_1-q_2\over
4}.$$  This is called the {\it Bianchi permutability formula\/} for the sine-Gordon
equation.   

\medskip

Next we  describe the action of $G^m_-$ on the spaces of solutions of the $j$-th
flow ($j\geq -1$). This construction is again using dressing action as in section  6
for the action of $G_-^m$ on $\cs_{1,a}$.  First we make some
definitions:
\ss

\refpar[kg] Definition. Let $\cm(j,a,b)$ denote the space of all solutions of the $j$-th
flow  (equation \refbi{}) on $\cs_{1,a}$ defined by $b$  with $[a,b]=0$ for $j=-1$ and $j\geq 1$
respectively.
 
\ms

Assume $j\geq 1$.  Let $A=a\l + u\in \cm(j,a,b)$, and $E(x,t,\l)$ the
trivialization of
$A$ normalized at $(x,t)=0$, i.e., 
$$\cases{E^{-1}E_x= a\l+u,&\cr E^{-1}E_t= b\l^j + Q_{b,1}\l^{j-1} + \cdots + Q_{b,j},&\cr E(0,0,
\l)= I.&\cr}$$ 
Given $g\in G_-^m$, by exactly the same method as in
section 6, we can factor 
$$g(\l)E(x,t,\l)=\tilde E(x,t,\l)\tilde g(x,t,\l),$$ such that $E(x,t,\cdot)\in G_+$ and
$\tilde g(x,t,\cdot)\in G_-^m$. Define $$\cases{g\bu E=\tilde E,&\cr g\bu A = \tilde
E^{-1}\tilde E_x&\cr}.$$ Then $g\bu A\in \cm(j,a,b)$, and $\bu$ defines an action of $G_-^m$ on
$\cm(j,a,b)$. 

Recall that the $g\in G_-^m$ can be generated by simple elements $g_{z,\pi}\in
G_-^m$, which are rational functions of degree $1$.
Choose a pole $z\in C\setminus R$, and a subspace $V\subset  C^n$, which is
identified with the Hermitian projection
$$\pi:C^n\to V.$$ Write 
$$g_{z,\pi}(\l)= \pi + {\l-z\over \l-\bar z} \pi^\perp$$ as in Proposition \refmg{}.

\refpar[rc] Definition.  $A\mapsto g_{z,\pi}\bu A$ is a B\"acklund transformation
for the $j$-th flow. 

\ms
Compute the action of $g_{z,\pi}$ explicitly as in section 6 to get:   

\refclaim[hm] Theorem. Let  $g_{z,\pi}$ be a generator in $G_-^m$, where $\pi$
is the projection of $C^n$ onto a
$k$-dimensional complex linear subspace
$V$. Let  $A= a\l + u\in \cm(j,a,b)$, and
$E(x,t,\l)$ the trivialization of $A$ normalized at $(x,t)=0$. Set $\tilde V(x,t)=
E(x,t, z)^*(V)$, and let $\tilde
\pi(x,t)$ denote the projection of $C^n$ onto $\tilde V(x,t)$.   Set 
\refeq[hs]$$\eqalign{\tilde \pi(x,t)&= E^*(x,t,z)U(U^*E(x,t,z)E^*(x,t,z)U)^{-1}U^*
E(x,t,z),\cr g_{z,\tilde \pi(x,t)(\l)} &=\tilde \pi(x,t) + {\l-z\over \l-\bar z}\tilde
\pi(x,t)^\perp\cr}$$ where $U$ is a
$n\times k$ matrix whose columns form a basis for $V$.
Then 
\item {(i)} $g_{z,\pi}\bu E = g_{z,\pi}E g_{z,\tilde \pi}\,^{-1}$,
\item {(ii)} $g_{z,\pi}\bu A= A + (z-\bar z)[\tilde \pi, a]$.

\refclaim[ht] Theorem. The $\tilde \pi$ constructed in Theorem \refhm{} is
the solution of the following compatible first order system:
\refeq[hq]$$\cases{(\tilde \pi)_x + [az+u,\tilde \pi] = (\bar z-z)[\tilde\pi,
a]\tilde\pi, &\cr (\tilde \pi)_t =(z-\bar z)\sum_{k=0}^j(-1)^k[\tilde \pi,
Q_{b,j-k}(u)](-z+(z-\bar z)\tilde\pi)^k,&\cr
\tilde\pi^*=\tilde\pi, \quad \tilde\pi^2=\tilde\pi, \quad \tilde\pi(0,0)=\pi.&\cr}$$
Moreover,
\item {(i)} equation \refhq{} is solvable for $\tilde \pi$ if and only if $A=a\l+u$ is a
solution of the $j$-th flow on $\cs_{1,a}$ defined by $b$,
\item {(ii)} if $A=a\l+ u$ is a solution of the $j$-th flow and $\tilde\pi$ is a solution
of equation \refhq{}, then $\tilde A= A+(z-\bar z)[\tilde\pi,a]$ is again a solution
of the $j$-th flow. \ei

\refpar[kp] Definition.
Let $R^*=\{r\in R\n r\not=0\}$ denote the multiplicative group, and
$R^\ast\sdp G_-^m$  the semi-direct product of $R^*$ and $G_-^m$ defined by the
homomorphism
$$\rho:R^\ast\to {\rm Aut\/}(G_-^m), \qquad \rho(r)(g)(\l)= g(r\l),$$
i.e., the multiplication in $R^\ast\sdp \W(G)$ is defined by
$$(r_1,h_1)\cdot(r_2,h_2)= (r_1r_2, h_1 (\rho(r_1)(h_2))).$$

\refclaim[kd] Theorem. The action $\bu$ of $G_-^m$ extends to an action of $R^*\sdp G_-^*$ on
the space $\cm(j,a,b)$ of solutions of the $j$-th flow on $\cs_{1,a}$ defined by $b$. In fact, if
$A=a\l + u\in \cm(j,a,b)$ and
$E$ is the trivialization of $A$ normalized at
$(x,t)=0$, then 
$$\cases{(r\bu E)(x,t,\l)= E(r^{-1}x, r^{-j}t, r\l),&\cr (r\bu A)(x,t,\l) = a\l + r^{-1}u(r^{-1}x,
r^{-j} t).\cr}$$ 

Since $(r^{-1},1)(1, g_{e^{i\a},\pi})(r, 1)= (1, g_{re^{i\a},\pi})$, we have

\refclaim[kh] Corollary. If $A\in \cm(j,a,b)$, then  $$r^{-1}\bu (g_{e^{i\a},\pi} \bu (r \bu A))=
g_{re^{i\a},\pi}\bu A.$$ 

Next we state an analogue of the Bianchi Permutability Theorem for the positive
flows:

\refclaim[jp] Theorem.  Let $z_1=r_1+is_1, z_2=r_2+is_2\in C\setminus R$ such that $r_1\not=
r_2$ or $s_1^2\not= s_2^2$, and
$\pi_1,\pi_2$ projections of $C^n$. Let $A_0\in \cm(j,a,b)$, and  $A_i=g_{z_i,\pi_i}\bu A_0$ for
$i=1,2$. Set 
\refeq[kx]$$\eqalign{\xi_i &=\bigl(-(z_1-z_2)I
+2i(s_1\pi_1-s_2\pi_2)\bigl)\pi_i\bigl((-(z_1-z_2)I +2i(s_1\pi_1-s_2\pi_2)\bigl)^{-1},\cr
\tilde \xi_i &=(-(z_1-z_2)I+ 2i(s_1\tilde \pi_1-s_2\tilde \pi_2))\tilde
\pi_i(-(z_1-z_2)I+ 2i(s_1\tilde \pi_1-s_2\tilde \pi_2))^{-1},\cr}$$ for
$i=1,2$, where $\tilde \pi_i$ is as in Theorem \refhm{} and $A_i= A_0+ 2is [\tilde \pi_i, a]$. 
Then 
\item {(i)} $g_{z_2,\xi_2}g_{z_1,\pi_1} = g_{z_1,\xi_1}g_{z_2,\pi_2}$,  
\item {(ii)}\refeq[jl]$$\eqalign{A_3 &= (g_{z_2,\xi_2}g_{z_1,\pi_1})\bu A_0= A_0+ 2i[s_1\tilde
\pi_1 + s_2\tilde\xi_2,a] \cr &= (g_{z_1,\xi_1}g_{z_2,\pi_2})\bu A_0 = A_0+
2i[s_1\tilde \xi_1 +s_2\tilde \pi_2, a].\cr}$$

Note that if $A$ lies in the orbit  of $R^*\sdp G_-^m$ through the vacuum solution
$A_0= a\l$ of the $j$-th flow then $A$ has no continuous scattering data.

\refpar[ke] Definition. An element $A$ in the orbit  of $R^*\sdp G_-^m$ through the vacuum
solution
$A_0= a\l$ of the
$j$-th flow is called a {\it  k-soliton\/} if the normalized eigenfunction has $k$
poles counted with multiplicities.  

\refpar[ky] Remark. The trivialization of the vacuum solution
$A_0=a\l$ of the $j$-th flow on $\cs_{1,a}$ defined by $b$ is
$E=\exp(a\l x + b\l^j t)$. So $1$-soliton is   
$$g_{z,\pi}\bu A_0= a\l + (z-\bar z)\left[e^{-a\bar z x - b \bar z^j t}U(U^*e^{a (z-\bar z) x+
b(z^j-\bar z^j)t}U)^{-1} U^*e^{azx+bz^j t}, \,\, a\right],$$
where $U$ is a matrix whose columns form a basis of Im$(\pi)$.  If 
$g=\prod_{i=1}^k g_{z_i,\pi_i}$, then $g\sharp
A_0$ can be written in terms of $1$-solitons
$g_{z_1,\pi_1}\sharp A, \cdots, g_{z_k,
\pi_k}\sharp A$ algebraically by applying the permutability formula \refjl{} repeatedly.

\refpar[kz] Remark. If $A$ is a soliton solution of the $j$-th
flow on $\cs_{1,a}$ defined by $b$, then by Corollary \refpd{} there exists $g\in G_-^m$ such that
  $$A(t)=\Psi(e_{b,j}(t)\natural g).$$ 
As noted in Remark \refcn{}, given $g_{z,\pi}\not\in H_-$, 
$g_{z,\pi}\bu A$ can be a $m+1, m$ or $m-1$-soliton if $A$ is a $m$-soliton solution.

\refclaim[kq] Proposition. The set of all soliton solutions of the $j$-th flow on $\cs_{1,a}$
defined by $b$ is isomorphic to the space $G_-^m/H_-^m$ of left cosets, where
$H_-^m=G_-^m\cap H_-$.

\refpar[maa] Remark. Although the space $G_-^m/H_-^m$ of all multi-solitons does not have
a manifold structure, the set  $\cb_m$  all $m$-solitons is the union of complex
algebraic varieties. For example, $$\eqalign{&\cb_1=\bigcup_{k=1}^{n-1}
(C_+\times X_k),\qquad {\rm where\/}\cr  &X_k=Gr(k,C^n)\setminus \{V\in Gr(k,
C^n)\n a(V)=V\}.\cr}$$  But
$\cb_m$ with $m\geq 2$ is much more complicated because of the non-uniqueness of the
factorization and the fact that generators of $G_-^m$  have complicated relations
such as the permutability formula  given in Theorem \refjp{} (i).  

Next we apply the action of $G_-^m$ on the first flows to get actions of $G_-^m$ on the space
$\cm$ of solutions of the $n$-dimension system \refee{} associated to $U(n)$. 
Given $v\in \cm$, the trivialization $E$ of $v$ normalized at the origin is the solution
$$\cases{E^{-1}E_{x_j}=a_j\l + [a_j,v],& $1\leq j\leq n$\cr E(0,\l)=I.&\cr}$$

\refclaim[av] Theorem. The group $R^*\times G_-^m$ acts on $\cm$, and 
the action $\bu$ is
constructed in the same manner as on the spaces of solutions of the first flow.  In
fact, given $g_{z,\pi}\in G_-^m$,  the
following initial value problem is solvable for $\tilde \pi$ and has a unique solution:
$$\cases{(\tilde \pi)_{x_j} +[a_jz+ [a_j,v], \tilde \pi] = (\bar z-z)[\tilde \pi,
a_j]\tilde
\pi,&\cr  \tilde \pi^*=\tilde \pi, \quad \tilde \pi^2=\tilde \pi, \quad \tilde
\pi(0)=\pi.&\cr}$$  Moreover, 
\item {(i)} $g_{z,\pi}\bu v = v - ((z-\bar z) \tilde \pi)^\perp$, where $y^\perp$
denote the projection of $y\in \cu$ onto $\ct^\perp$,
\item {(ii)} the trivialization of $g_{z,\pi}\bu v$ is $g_{z,\pi}E g_{z,\tilde \pi} ^{-1}$,
\item {(iii)} $\tilde \pi$ is the projection onto the linear subspace $E_z^*(V)$, where $V$ is
the image of the projection $\pi$,
\item {(iv)} $(r\bu v)(x)=r^{-1}v(r^{-1}x)$ for $r\in R^*$.\ei

\refpar[ar] Remark. The permutability Theorem \refjp{} holds for system \refee{} with the
same formula.
\ss

There are also analogous results for the $-1$-flow:

\refclaim[ki] Theorem. The group $R^*\sdp G_-^m$ acts on the space $\cm(-1,a,b)$ of solutions of
the $-1$-flow on $\cs_{1,a}$ defined by $b$:
\refeq[kr]$$\cases{u_t=[a, g^{-1}bg],&\cr g_t= g u, & $\lim_{x\to -\infty} g(x)=
I$.\cr}$$ Moreover, let $A\in \cm(-1,a, b)$, $E$ the
trivialization of $A$ normalized at $(x,t)=0$,  and
$g_{z,\pi}$ a simple element of $G_-^m$, then 
\item {(i)} Theorems \refhm{}, \refjp{} and Corollary \refkh{} hold with the same
formulas,
\item {(ii)} $\tilde \pi$ is the solution to
\refeq[kk]$$\cases{(\tilde \pi)_x + [az+u,\tilde \pi] = (\bar z-z)[\tilde\pi, a]\tilde\pi, &\cr
(\tilde \pi)_t ={1\over \n z\n^2} \left( (z-\bar z)\tilde \pi g^{-1}bg\tilde \pi - zg^{-1}bg\tilde \pi
+ \bar z \tilde
\pi g^{-1}bg\right),&\cr
\tilde\pi^*=\tilde\pi, \quad \tilde\pi^2=\tilde\pi, \quad \tilde\pi(0,0)=\pi,&\cr}$$
\item {(iii)} for $r\in R^*$ we have 
$$\cases{(r\bu E)(x,t,\l)= E(r^{-1}x, rt, r\l),&\cr (r\bu A)(x,t,\l) = a\l + r^{-1}u(r^{-1}x,
r t).&\cr}$$

\refpar[mba] Remark. Let $A\in \cm(-1,a,b)$, and $E$ its trivialization normalized at the origin.
Then $s(x,t)=E(x,t,-1)E(x,t,1)^{-1}$ is a harmonic map from $R^{1,1}$
with the metric
$2dxdt$ to
$U(n)$ and $s^{-1}s_x$ is conjugate to $a$ and $s^{-1}s_t$ is
conjugate to $b$. In particular, this says that $\cm(-1,a,b)$ is a subset of the space $\cf$ of
harmonic maps from
$R^{1,1,}$ to $U(n)$. The action of $G_-^m$  on $\cf$ constructed in [U1] leaves $\cm(-1,a,b)$
invariant, and agrees with the action $\bu$ we constructed here. 
\ms

Recall that given an involution $\s$ of $su(n)$, 
$G_-^{m,\s}$ is the subgroup of $g\in G_-^m$ such that $\s(g(-\l))=g(\l)$ and 
$\cs_{1,a}^\s=\{A\in \cs_{1,a}\n \s(A(-\l))=A(\l)\}$. Since the action
of $G_-^{m,\s}$ leaves $\cs_{1,a}^\s$  invariant,  we have

\refclaim[cq] Theorem.  The space $\cm^\s(2j-1, a,b)$ of  solutions of the $(2j-1)$-th flow associated to $U/K$
defined by $b$ is a subset of $\cm(2j-1,a,b)$ for $j\geq 0$. Moreover, the action $\bu$ of
$R^*\sdp G_-^{m,\s}$ leaves $\cm^\s(2j-1,a,b)$ invariant, and Theorems \refhm{}, \refht{},
\refkd{} and Corollary \refkh{} hold for $\cm^\s(2j-1,a,b)$ and $R^*\sdp G_-^{m,\s}$. 

\refclaim[mea] Theorem. The space 
$\cm^\s(-1,a,b)$  of solutions of the $-1$-th flow equation on $\cs^\s_{1,a}$
defined by $b$ is a subset of $\cm(-1,a,b)$.  Moreover, the action of $R^*\sdp G_-^{m,\s}$ leaves
$\cm^\s(-1,a,b)$ invariant and Theorem
\refki{} holds for $\cm^\s_{1,a}$ and $R^*\sdp G_-^{m,\s}$.

\refclaim[ku] Proposition. Let $\s$ denote the involution on $SU(n)$ defined by $\s(y)=\bar y$.
Then
\item {(i)} $g_{z,\pi}\in G_-^{m,\s}$ if and only if $z=-\bar z$ and $\bar \pi=\pi$,
\item {(ii)} if $\bar \pi= \pi$, then $g_{z,\pi}g_{-\bar z,\pi}\in G_-^{m,\s}$.

Next we explain the relation between the classical B\"acklund transformations and
the action of $G^{m,\s}_-$ on the space of solutions of the sine-Gordon equation (or the space
$\cm^\s(-1,a,b)$ with $\s$, $a, b$ defined as in Example \refct{}).  If $s\in R$, $\tilde\pi
^*=\tilde \pi =(\tilde
\pi)^t$ and $(\tilde \pi)^*=\tilde\pi$, then by Proposition \refku{}, $g_{is,\tilde \pi}\in
G_-^{m,\s}$. Hence 
$$\tilde\pi=
\pmatrix{\cos^2{ f\over 2} &\sin { f\over 2}\cos { f\over 2}\cr\sin { f\over 2}\cos { f\over 2}&
\sin^2 { f\over 2}\cr}$$ for some function $f$, i.e., $\tilde\pi$ is the projection onto
$\pmatrix{\cos{f\over 2}\cr \sin{f\over 2}\cr}$.   So the first order system \refkk{}  for
$\tilde\pi$  becomes 
\refeq[hta]$$\cases{f_x= {q_x\over 2} + 2s \sin f, &\cr f_t= {1\over 2s} \sin
(f-q).&\cr}$$ Write
$$\tilde u= g_{is,\b}\bu u= \pmatrix{0& \tilde q_x/2\cr -\tilde q_x/2&0\cr}.$$ But
$\tilde u= u + 2is[\tilde\pi,a]$, hence we have
$\tilde q= 2f-q$.  Writing equation \refhta{} in terms of $\tilde q$, we get 
$$\cases {(q^*-q)_x= 4s \sin ({q^*+q\over 2})\cr 
(q^*+q)_t = {1\over s} \sin ({q^*-q\over 2}),\cr}$$
which is  the classical B\"acklund transformation  for  the sine-Gordon equation. So we have:

\refclaim[jg] Proposition. Let $q$ be a solution of the sine-Gordon equation, and $0<c_0<\pi$.
Set 
$$\cases{A=a\l + \pmatrix{0&{q_x\over 2}\cr -{q_x\over 2} & 0},&\cr  f_o=1/2( q(0,0)+
c_o)&\cr
\pi=\pmatrix{\cos^2{ f_0\over 2} &\sin { f_0\over 2}\cos { f_0\over 2}\cr\sin { f_0\over 2}\cos
{ f_0\over 2}& \sin^2 { f_0\over 2}\cr}.&\cr}$$  Then 
$$B_{s,c_o}(q)= g_{is,\pi}\bu A.$$ (We will still use $\bu$ to denote the induced action of
$G_-^{m,\s}$ on the set of solutions of the sine-Gordon equation). 

\refclaim[kn] Proposition. Let  $q$ is a solution of the sine-Gordon equation. Then:
\item {(i)}  $s\bu q $ is the Lie
transformation $L_s(q)$. 
\item {(ii)} Proposition \refko{} is a consequence of the following equality
in the group $R^*\sdp G_-^m$:
$$(s^{-1},I)(1, g_{e^{i\a},\b})(s,I)= g_{se^{i\a},\b}.$$
\item {(iii)} The Permutability Formula \refjf{} in Theorem \refje{} is the same
formula \refjl{} given in Theorem \refjp{}, which follows from the following relation
of the generators of $G_-^m$:
$$g_{z_1,\xi_1}g_{z_2,\pi_2}= g_{z_2,\xi_2}g_{z_1,\pi_1},$$ where $\pi_i$ and $\xi_i$ are
related by formula \refkx{}. 

\refpar[fn] Remark.  It is noted by Xi Du [Du] that a classical Ribaucour
transformations as defined in Darboux ([Da1]) for surfaces of $K=-1$ in $R^3$
correspond to the action of an element
$g\in G_-^{m,\s}$, which is the product of two simple elements as in Proposition
\refku{} (ii). 

\ms 

Using the action of $G_-^m$, we obtain many solutions of the $j$-th flow that are periodic in
time.  This is an algebraic calculation, which shows that when the poles are
properly placed, the solutions are periodic.  Multi-solitons will be time periodic if
the periods of the component solitons are rationally related. 

\refclaim[df] Theorem.  Let $j>1$ be an integer, $a=\diag(ia_1,
\ldots, ia_n)$, and
$b=\diag(ib_1, \ldots, ib_n)$.  If $b_1, \ldots, b_n$ are rational numbers. Then the
$j$-th flow equation  defined by $a,b$: $$u_t=[Q_{b,j+1}(u),a]$$
has infinitely many $m$-soliton solutions that are periodic in $t$. 

The trivialization of the vacuum solution for the $-1$-flow defined by $a=\diag(i,\cdots, i, -i,
\cdots, -i)$ is $E(\l,x,t)=\exp(a(\l x+ \l^{-1}t))$. By Theorem \refki{}, the $1$-soliton
$g_{e^{i\o},\pi}\bu 0$ is a function of 
$$\exp(i(\cos\o (x+t)-i\sin(x-t)))=\exp(i\tau\cos \o + y \sin\o),$$ where $y=x-t$ and
$\tau=x+t$  are the space-time coordinates. This gives 

\refclaim[pa] Theorem. If $z=e^{i\o}$ and  $a=\diag(i,\cdots, i, -i,
\cdots, -i)$, then the $1$-soliton $g_{z,\pi}\bu 0$ for the
$-1$-flow (harmonic maps from $R^{1,1}$ to $SU(n)$)) is periodic in time with period
${2\pi\over\cos\o}$. A multiple soliton generated by a rational loop with poles at
$z_1=e^{i\o_1}, \cdots, z_r=e^{i\o_r}$ will be periodic with period $\tau$ if there
exists integers $k_1, \cdots, k_r$ such that 
$$\tau={2\pi k_j\over \cos \o_j} \qquad \forall \,\, 1\leq j\leq r.$$

The multi-solitons above satisfy the sine-Gordon equation if the rational loop
satisfies
$\overline{f(-\bar\l)}=f(\l)$. This means the poles occur in pairs $(e^{i\o_j}, -e^{-i\o_j})$ and
the projection matrices $\pi_j$ must be real. 

\refclaim[pb] Corollary. Multiple-breather solutions exists for the
sine-Gordon equation. 

\refpar[kv] Example.  If $\pi$ is real, then
$$(g_{e^{i\o},\pi}g_{-e^{-i\o}, \pi})\bu 0= 4 \tan^{-1} \left( {\sin \o \sin ( (x+t)\cos\o)\over
\cos\o \cosh ((x-t)\sin \o)}\right)$$ is the classical breather solution for the sine-Gordon
equation. Theorems \refmea{} and \refpa{} give $m$-breather solutions
explicitly, although the computations are quite long.  

\refclaim[dma] Corollary. There are infinitely many
harmonic maps from $R^{1,1}$ to a symmetric space that are periodic in time.

\bs

\newsection Geometric non-linear Schr\"odinger equation.\par

Consider the evolution of curves in $R^3$
\refeq[ew]$$\g_t= \g_x\times \g_{xx},$$ where $\times$ denote the cross-product
in $R^3$. This equation is known as the {\it vortex filament equation \/}, and has a
long and interesting history (cf. [Ri]). It is easy to see that $\N
\g_x\N^2$ is preserved under the evolution.  It follows that  if
$\g(\cdot, 0)$ is parametrized by its arc length, then so are all $\g(\cdot, t)$ for all
$t$. So  equation \refew{} can also be viewed as the
evolution of a curve that moves along the direction of binormal with the curvature
as its speed.  Let
$k(\cdot,t)$ and
$\tau(\cdot,t)$ be the curvature and torsion of the curve $\g(\cdot, t)$. Then there
exists a unique $\o(x,t)$ such that $\o_x=\tau$ and 
$$q(x,t)= k(x,t)e^{-i\int^x \tau(s,t)ds}$$ is a solution of the non-linear Schr\"odinger equation:
$$q_t={i\over 2} \left(q_{xx} + 2\n q\n^2 q\right).$$ 
There is another interesting evolution of curves in $S^2$ that is also associated to
the non-linear Schr\"odinger equation:

\refclaim[cj] Proposition. $\g(x,t)$ is a solution of equation \refew{} with $x$
as the arc length parameter if and only if $\phi(x,t)= \g_x(x,t):R^2\to
S^2$ satisfies the equation 
\refeq[ex]$$J(\phi_t)=\K_{\phi_x}\phi_x,$$ where $\K$ is the
Levi-Civita connection and $J$ is the complex structure of the standard
two sphere $S^2$.

Equation \refex{} is the geometric non-linear Schr\"odinger equation (GNLS) on $S^2$.  Such
equation can be defined on any complex Hermitian manifold $(M,g,J)$. Consider the
Schr\"odinger flow on the space
$\cs(R,M)$ of Schwartz maps, i.e., the equation for maps $\phi:R\times R\to M$:
$$J({\partial \phi\over \partial t})= \D \phi,$$ where $\D \phi=
\K_{\phi_x}\phi_x$ is the gradient of the energy functional on $\cs(R,M)$, or the
accelleration.

In this section, we  give a brief outline of ideas and results in a forthcoming paper
[TU3].  There is a Hasimoto type transformation that transforms the GNLS equation
associated to  $Gr(k,C^n)$ to the the second flow on $\cs_{1,a}$ defined by  
\refeq[rd]$$a=\left(\matrix {iI_k&0\cr 0&-iI_{n-k}\cr}\right)\in u(k)\times
u(n-k).$$  We have seen in Example \refca{} that identifying
$\cs_{1,a}$ as the space
$\cm_{k\times (n-k)}$ of $k\times (n-k)$ matrices, the second flow defined by $a$
is the matrix non-linear Schr\"odinger equation:
\refeq[mja]$$B_t = {i\over 2} (B_{xx} + 2 BB^\ast B).$$ 
Applying our theory to equation \refmja{}, we obtain many beautiful properties for
the GNLS associated to $Gr(k,C^n)$. For example, we have
\item {(i)}  a Hamiltonian formulation, 
\item {(ii)} long time existence for the Cauchy problem, 
\item {(iii)} a sequence of commuting Hamiltonian flows, 
\item {(iv)} explicit soliton solutions, 
\item {(v)} a non-abelian Poisson group action on the space of solutions of the GNLS,
\item {(vi)} a sequence of compatible symplectic  structures on the space $\cs(R,Gr(k,C^n))$
in which the GNLS is Hamiltonian and has a Lenard relation. 

\ms

Let $U(n)$ be equipped with a bi-invariant metric. It is well-known that $Gr(k,C^n)$ can be
naturally embedded as a totally geodesic submanifold $M$ of 
$U(n)$. In fact, $M$ is the set  of all $X\in U(n)$ such that $X$ is conjugate to $a$ as
described by formula \refrd{}. 
 The invariant complex structure on $M$ is given by
$$J_h(v)=[v,h].$$ Consider the following equation for maps $\phi:R^2\to M$:
\refeq[dj]$$J_\phi(\phi_t)= \K_{\phi_x}(\phi_x),\qquad \phi:R^2\to M$$ where $\K$
is the Levi-Civita connection of the standard Kahler metric on $M$.  A direct
computation gives
$$\K_{\phi_x}(\phi_x)= \phi (\phi^{-1}\phi_x)_x.$$
So equation \refdj{} becomes
$$\phi_t= -{1\over 2} (\phi^{-1}\phi_x)_x. \eqno(GNLS)$$

Next we want to associate to each solution of equation \refmja{} a solution of the GNLS. This is
a generalization of the Hasimoto transformation of the vortex filament equation to non-linear
Schr\"odinger equation.  As noted in Example
\refca{},  
$A=a\l+u\in$ is a solution of \refmja{} if and only if 
$$\o_\l= (a\l+ u)dx + (a\l^2+ u\l + Q_{a,2}(u))dt $$ is flat for all $\l$, where 
$$u=\pmatrix{0&B\cr -B^*&0\cr}, \qquad Q_{a,2}=\pmatrix{{1\over 2i} BB^*& {i\over 2}
B_x\cr {i\over 2} B_x^*& -{1\over 2i} B^*B\cr}.$$  In particular,
$\o_0=udx+Q_{a,2}(u)dt$ is flat. Let $g$ be the trivialization of $\o_0$, i.e., 
$$\cases{g^{-1}g_x= u,&\cr g^{-1}g_t= Q_{a,2}(u).&\cr}$$
Set $\phi=gag^{-1}$. Changing the gauge of $\o_\l$ by $g$ gives
$$\eqalign{\tau_\l &= g\o_\l g^{-1}-dg g^{-1}= (gag^{-1}\l)dx +(gag^{-1} \l^2 +
g_xg^{-1}\l) dt\cr &= \phi \l dx + (\phi \l^2 + Q_{a,2}(u)\l)dt.\cr}$$  Since
$\tau_\l$ is flat for all $\l$,  we get
\refeq[do]$$\cases{\phi_t=(Q_{a,2}(u))_x,&\cr \phi_x= -[\phi,Q_{a,2}(u)].&\cr}$$
But for $u\in \cu_a^\perp$, we have $a^{-1}ua=-u$. Hence
\refeq[ev]$$\eqalign{\phi^{-1}\phi_x &=
ga^{-1}g^{-1}(g_xag^{-1}-gag^{-1}g_xg^{-1})
\cr &=ga^{-1}uag^{-1}-g_xg^{-1}=-gug^{-1}-g_xg^{-1}= -2 g_xg^{-1}= -2Q_{a,2}(u).\cr}$$ So
the first equation of \refdo{} implies that $\phi$ is a solution of the GNLS.

\ms
Conversely, suppose $\phi:R^2\to M$ is a solution of the GNLS. Then there exists
$g:R^2\to U(n)$ such that $\phi= gag^{-1}$ and  $g^{-1}g_x(x,t)\in \cu_a^\perp$ for
all $(x,t)$. Set $$u=g^{-1}g_x, \qquad f= -{1\over 2} (\phi^{-1}\phi_x).$$
Then the equation \refev{} implies that $\phi^{-1}\phi_x= g_xg^{-1}= f$.
Differentiate $\phi=gag^{-1}$ with respect to $x$ to get 
$$\phi_x=[g_xg^{-1},\phi] = [f,\phi].$$ But the GNLS gives $\phi_t=f_x$. So 
$$(\phi \l)dx + (\phi \l^2 + f\l)dt$$ is flat for all $\l$. Changing the gauge by $g^{-1}$,
we find that
$$\b_\l= (a\l + u)dx + (a\l^2 + u\l +h)dt$$ is flat, where $h=g^{-1}g_t$.  Flatness of $\b_\l$ on the
$(x,t)$-plane for all $\l$ implies that $h=Q_{a,2}(u)$. So this proves that
$u$ is a solution of the second flow equation \refmja{}. To summarize, 

\refclaim[mga] Proposition. If $B:R^2\to \cm(k\times (n-k))$ is a solution of  the
matrix NLS \refmja{}, then there is
$g:R^2\to U(n)$ such that $$g^{-1}g_x= \pmatrix{0&B\cr -B^*&0\cr}, \qquad g^{-1}g_t
=\pmatrix{{1\over 2i} BB^*& {i\over 2} B_x\cr {i\over 2} B_x^*& -{1\over 2i} B^*B\cr}$$ and
$\phi=gag^{-1}$ is a solution of the GNLS. Conversely, if $\phi$ is a solution of the
GNLS, then there is $g:R^2\to U(n)$ such that $$\phi=g\pmatrix{iI_k&0\cr
0&-I_{n-k}\cr}g^{-1}, \qquad g^{-1}g_x=\pmatrix{0&B\cr -B^*&0\cr},$$ and $B$ is a solution of
the matrix NLS equation \refmja{}. 

To end this section, we will translate properties for the second flow \refmja{} to properties of
the GNLS equation. 
 
When $a$ is singular, formula \refge{} implies that  the corresponding
Hamiltonians for the first three flows on $\cs_{1,a}$ defined by $a$ are 
$$\eqalign{F_1(B) &= {1\over 4}\int_{-\infty}^\infty
\tr(iB_xB^*-iBB_x)dx={1\over 4}<iB_x,B>,\cr F_2(B) &= {1\over 4}
\int_{-\infty}^\infty
\tr(-B_xB_x^\ast + B^\ast BB^\ast B)dx \cr & ={1\over 8}(-<B_x,B_x>+
<B^*B,B^*B>),\cr F_3(B) &= {i\over 16}\int_{-\infty}^\infty
\tr(-BB^*_{xxx}+B_{xxx}B^*)+ 3\tr  (-BB^*BB^*_x+B_xB^*BB^*)dt \cr &= -{1\over
16}(<iB_{xxx},B> + 3<iB_x,BB^*B>).\cr}$$

\refpar[mha] Remark. If
$b\in \cu_a$, i.e., $[a,b]=0$, and $b\not=a$, then 
\item {(i)} $Q_{b,j}(u)$ in general is not a local operator in $u$,
\item {(ii)} the flow generated by $b\l^n\in \ch_+$  commutes with the flow
\refmja{},
\item {(iii)} given $b_1,b_2\in \cu_a$, the flow generated by $b_1\l^k$ and $b_2\l^j$ need not
commute and $$[\xi_{b_1,k}, \xi_{b_2,j}]=\xi_{[b_1,b_2],k+j},$$ where $\xi_{b,m}$ denotes
the 
 infinitesimal vector fields corresponding to $b\l^m$,
\item {(iv)} the action of $H_+=\{g\in G_+\n ga=ag\}$ on $\cs_{1,a}$ is Poisson. 

\ms

It follows from the discussion in section 7 that the Cauchy problem for equation
\refmja{} with initial condition  $A_0=a\l+u_0\in \cs'_{1,a}$ can be solved by using
 factorizations.  First  we use the direct scattering of  $A_0=a\l+u_0$ on
the line, i.e., solve equation \refgx{}. 
Set $f(\l)=\psi(0,\l)$. Then $f\in D_-$.  
Decompose $$f(\l)e^{a\l x+b\l^j t}= E(x,t,\l)M(x,t,\l)^{-1}$$ as in section 7 by applying Birkhoff
decompositions repeatedly. Then
$A=E^{-1}E_x$ is the solution of the Cauchy problem. 

The rational group $G_-^m$ acts on the space of solutions of the GNLS, and soliton
solutions can be calculated explicitly using the formulas in Theorem \refhm{}. 

\refclaim[mka] Proposition. Let $a=\diag(-i, i, \cdots, i)$, and choose a pole
$z=r+is\in C\setminus R$. Let
$\pi$ be the projection on the subspace spanned by 
$(1, v)^t=(1, v_2,
\cdots, v_n)^t$. Then the one-solitons for the $j$-th flow defined by $a$ generated
by B\"acklund transformations from the vacuum solution
$A_0(x,t,\l)=a\l$ are of the form
$$\eqalign{A(x,t,\l)&= a\l + u(x,t),\cr u(x,t) &=
\pmatrix{0&B(x,t)\cr-B^\ast(x,t)&0\cr},\cr}$$ where
$$B(x,t)= {4s e^{-2i(rx+{\rm Re\/}(z^j)t)}\over e^{-2(sx + \Im(z^j) t)}+e^{2(sx +
\Im(z^j) t)}\N v\N^2}\,\,\bar v.$$

\proof We use Theorem \refhm{} to make our computations.
We start with $A_0=a\l$. According to Theorem \refhm{}, 
$$A_0\mapsto a\l + (z-\bar z) [\tilde \pi, a],$$ where $a=\diag(-i, i, \cdots, i)$.
Here $\tilde \pi$ is the projection on $(1,\tilde v)^t$, where 
$$\hat v= (1,\tilde v) = (1, e^{2i(\bar z x + \bar z^j t)} v).$$
Let
$z=r+is$. Then
$$\tilde \pi(x,t)={1\over e^{-2(sx +\Im (z^j) t)} + e^{2(sx +\Im (z^j) t)} \N
v\N^2}\,\,  \hat v ^\ast \hat v.$$ The formula for $B$ follows. \qed

\bs

\newsection First flows and flat metrics.\par

	  The integrable equations of evolution we have been describing up to this point
have at most two independent variables. The flow of the first variable, regarded as
a spacial variable, is used to construct the initial Cauchy data from the scattering
coset (hence the ``first flow'' terminology). The second variable is considered to be
the time variable, and the flow in this variable is the evolution. Many authors
consider a commuting heirarchy of flows to generate functions of an infinite
sequence of time variables. However, the physical and geometric applications do
not require this consideration. 

We turn our attention to a family of geometric problems in $n$ spacial variables,
which we shall call $n$-{\it dimensional systems\/} or $n$-{\it dimensional flows\/}.
In the applications, the $n$ variables are on an equal footing, and the flows in each
variable is a first flow. The flows commute, and hence the resulting geometric
object is always a flat connection on a region of $R^n$ with  special properties. From
our viewpoint, the natural parameter (moduli) space of solutions is a coset space
of the sort we have just described. In many cases, we have obtained global results
on connections in $R^n$ via the decay theorems in section 7.

The $n$ commuting first flows associated to a rank $n$ symmetric space have been
discussed in a paper by the first author ([Te2]). We outline the general theory and
give some of the basic examples.  The results on coset spaces and B\"acklund
transformations apply naturally to these systems. 
	
\refpar[al] Definition ([Te2]). Let $U$ be a rank $n$ Lie group,  $\ct$ a maximal abelian
subalgebra of the Lie algebra $\cu$, and $a_1, \cdots, a_n$ a basis of $\ct$.  {\it The
$n$-dimensional system associated to $U$\/} is the following first order system:
\refeq[eea]$$[a_i,v_{x_j}]-[a_j, v_{x_i}] = [[a_i,v], [a_j,v]], \qquad v:R^n\to \ct^\perp.$$

\refpar[am] Definition ([Te2]). Let $U/K$ be a rank $n$ symmetric space, $\s:\cu\to \cu$ the
corresponding involution, $\cu=\ck+\cp$ the Cartan decomposition, $\ca$ a maximal
abelian subalgebra in $\cp$, and $a_1, \cdots, a_n$ a basis of $\ca$. {\it The $n$-dimensional
system associated to $U/K$ \/} is the first order system:
\refeq[ef]$$[a_i,v_{x_j}]-[a_j, v_{x_i}] = [[a_i,v], [a_j,v]], \qquad v:R^n\to
\cp\cap\ca^\perp.$$

\refclaim[an] Theorem ([Te2]). The following conditions are equivalent: 
\item {(i)} $v$ is a solution of equation \refeea{} (or \refef{}) 
\item {(ii)} the connection $1$-form $\o=\sum_{j=1}^n (a_j \l + [a_j,
v])dx_j$ is left flat, i.e., $d\o = -\o \wedge \o$,
\item {(iii)} $[{\partial \over \partial x_i} +(a_i\l + [a_i,v]),
{\partial \over \partial x_j} +(a_j\l + [a_j,v])]=0$ for all $ i\not=j$.\ei

\refclaim[cd] Theorem ([Te2]). Let $\ct$ be a maximal abelian subalgebra of $\cu$, $\ct^\perp$
the orthogonal complement of $\ct$ in $\cu$, and let
$\cs(R^n,\ct^\perp)$ denote the space of Schwartz maps from $R^n$ to $\ct^\perp$. 
Let $a_1, \cdots, a_n$ be regular elements and form a basis of
$\ct$. Then there exists a dense open subset $\cs_0$ of $\cs(R^n,\ct^\perp)$ such that given
$v_0\in \cs_0$, the following Cauchy problem for equation
\refeea{} has a unique solution:
$$\cases{[a_i,v_{x_j}]-[a_j,v_{x_i}]=[[a_i,v],[a_j,v]],&if $i\not=j$,\cr 
v(t,0,\cdots, 0)=v_0(t).&\cr}$$

At this point, it is important to give some explaination and application. Because
these flows all commute, a change of basis in the abelian subalgebra $\ct$ or $\ca$
can be represented by composition with an element of $GL(n,R)$. So we might as
well assume that $a_1, \cdots, a_n$ are generic or regular (have distinct
eigenvalues). Starting at $0\in R^n$, given an element in the coset space
$D_-/H_-$, we can solve for $E(x_1, 0, \cdots, 0, \l)$ and find ${d\over dx} + \l a_1 + u(x_1, 0,\cdots,0)$ as if
we were solving for initial Cauchy data. Instead of going to one of the
heirarchy of flows, we solve for the entire family of first flows in variables $(x_1,
\cdots, x_n)$. This gives us a map
$$D_-/H_- \,\mapsto \{\,{\rm flat\, connections \, on\, a \, region\, of\,} R^n\}.$$
In the case that our flows can be embedded in the unitary flows, 
$$D_-/H_- \,\mapsto \{\,{\rm flat\, connections \, on\, } R^n \, {\rm decaying \, at\, }
\infty\}$$ (by Theorem \refhc{}). This gives a proof of Theorem \refcd{}.

\refpar[rf] Remark. Note that we are constructing a more rigid structure then a
flat connection. We are actually constructing special connection one-forms, and
we do not allow arbitray coordinate change or gauge changes in the theory.  

By expressing the parameter space in this form, we have made a beginning towards
thinking about the natural symplectic structure on this solution space. In the case
of one-dimensional Cauchy data, the symplectic structures were averaged out
over the one-parameter flow. Here we need to average them out over an
$n$-dimensional flow to obtain a natural structure. 

The canonical examples of these flat connections are quite easy to describe.

\refpar[cv] Examples.
\ss
{\bf Example (i)\/}\hskip 4pt Let $U=GL(n,R)$, $\ct$ the maximal abelian subalgebra of
diagonal matrices, and $\{e_{11}, \cdots, e_{nn}\}$ a basis of $\ct$, where $e_{ij}$ denote the
matrix in
$gl(n)$ all whose entries are zero except that the $ij$-th entry is equal to $1$. Then
the
$n$-dimensional system associated to $GL(n,R)$ is the system for 
$$f=(f_{ij}):R^n\to gl(n,R), \quad f_{ii}=0,  \,\, 1\leq i\leq n$$
\refeq[ema]$$\cases{ (f_{ij})_{x_i}
+ (f_{ij})_{x_j} +
\sum_k f_{ik}f_{kj}=0, & if
$i\not=j$,\cr  (f_{ij})_{x_k}= f_{ik}f_{kj}, & if $i, j, k$ are distinct.\cr}$$

\ss
{\bf Example (ii)\/}\hskip 4pt Let $U/K= GL(n,R)/SO(n)$, and $\cu=\ck+\cp$ the corresponding
Cartan decomposition. Then 
$\cp$ is the set of all real symmetric $n\times n$ matrices, the space $\ca$ of all diagonal
matrices is a maximal abelian subalgebra in $\cp$, $e_{11}, \cdots, e_{nn}$ form a
basis of $\ca$, and $\cp\cap \ca^\perp$ is the space $gl_s(n)$ of all symmetric
$n\times n$ matrices whose diagonal entries are zero. The $n$ dimensional system \refef{}
associated to $GL(n,R)/SO(n)$ is the system for $$F=(f_{ij}):R^n\to gl(n, R), \quad f_{ij}=f_{ji},
\quad f_{ii}=0\,\,\, {\rm if \/}\,\, 1\leq i\leq n$$
\refeq[em]$$\cases{(f_{ij})_{x_i} + (f_{ij})_{x_j} + \sum_k f_{ik}f_{kj}=0, & if $i\not=j$,\cr 
(f_{ij})_{x_k}= f_{ik}f_{kj}, & if $i, j, k$ are distinct.\cr}$$
Note that system \refem{} is the system \refema{} restricted to maps $f=(f_{ij})$ that
are symmetric. 

\ss
{\bf Example (iii)\/}\hskip 4pt Let $U/K=U(n)/SO(n)$, and $u(n)=so(n) + \cp$ a Cartan
decomposition. Then $\cp$ is the set of all symmetric pure imaginary $n\times n$ matrices and the
space $\ca$ of all diagonal matrices in $\cp$ is a maximal abelian algebra. Let $ia_1, \cdots, ia_n$ be
a basis of $\ca$. Write $v:R^n\to \cp\cap \ca^\perp$ as $v=-i F$, where $F$ is a real
$n\times n$ symmetric matrix. Then equation \refef{} for $v$ is the equation \refem{} for $F$.
This is a special case of the general fact that the $n$-dimensional system associated to a
compact symmetric space is the same as that associated to its non-compact dual. 
\ss

{\bf Example (iv)\/}\hskip 4pt Let $U/K=SO(2n)/S(O(n)\times O(n))$, and 
$\cu=\ck+\cp$ the corresponding Cartan decomposition. Then 
$$\eqalign{\ck &= so(n)\times so(n)= \bigg\{\pmatrix{B&0\cr 0&C}\bigg|\,\, B, C\in
so(n)\bigg\},\cr
\cp &= \bigg\{\pmatrix{0&F\cr -F^t&0}\bigg|\,\, F\in gl(n)\bigg\},\cr}$$ and 
$$\ca =\bigg\{\pmatrix{0&-D\cr D&0\cr}\bigg|\,\, D\,\, {\rm is\,\,
diagonal\,}\,\bigg\}$$ is a maximal abelian subalgebra of $\cp$. Let
$a_i=\pmatrix{0&-e_{ii}\cr e_{ii}&0\cr}$. Then $a_1, \cdots, a_n$ form a basis of $\ca$, and
$\cp\cap\ca^\perp$ is the set of matrices of the form $\pmatrix{0&X\cr -X^t&0\cr}$ such
that $X=(x_{ij})$ is $n\times n$ matrix with $x_{ii}=0$ for all $i$.  Then
 equation \refef{} for $v=\pmatrix{0&F\cr -F^t&0\cr}$, with $F=(f_{ij}):R^n\to
gl(n,R)$ and $f_{ii}=0$ for all $1\leq i\leq n$, is 
\refeq[el]$$\cases{(f_{ij})_{x_i}+ (f_{ji})_{x_j} +\sum_k f_{ki}f_{kj} =0,& if $i\not=j$ \cr
(f_{ij})_{x_j}+ (f_{ji})_{x_i} + \sum_k f_{ik}f_{jk}=0,& if $i\not=j$\cr 
(f_{ij})_{x_k} = f_{ik}f_{kj}, & if $i, j,k$ are distinct.\cr}$$

Up to now, the flat connections were not constructed to relate to Riemannian
geometry. To explain the relation of these flat connections to geometry, we need
to set up some notations.  A {\it diagonal metric\/} is a metric of the form
$$ds^2=\sum_j b_j(x)^2dx_j^2.$$ 
If this diagonal metric is flat, then $(x_1, \cdots, x_n)$ is an orthogonal coordinate system on
$R^n$ in the sense of Darboux ([Da2]).  These examples arise in the study of
isometric immersions of constant sectional curvature $n$-manifolds into
Euclidean space. 

On the other hand, to study Lagrangian flat submanifolds in $C^n$ or Frobenius
manifolds (used in quantum cohomology), we consider {\it Egoroff metrics\/}.
These are metrics of the form 
$$ds^2=\sum_j\phi_{x_j} dx_j^2$$ for some function
$\phi$.  

The Levi-Civita connection $1$-form for the diagonal metric 
$\sum_{i=1}^n b_i(x)^2
dx_i^2$ is $$w=(w_{ij})=(-f_{ij}dx_i + f_{ji}dx_j), \qquad f_{ij}={(b_i)_{x_j}\over b_j}.$$ Or
equivalently, 
$$w= -\d F+ F^t \d, \qquad {\rm where\,\,} \d=\diag(dx_1, \cdots, dx_n).$$
Hence we are looking for flat connections o this special form. 

 The Levi-Civita connection of a Egoroff metric is  $w= [F,\d]$
with $F=F^t$.  It is easy to see that a diagonal metric is {\it Egoroff\/} if and only if
$f_{ij}=f_{ji}$. 
\ms

\refpar[fo] Definition.  A {\it Darboux connection\/} is a connection
of the form $-\d F+ F^t \d$, and an {\it Egoroff connection\/} is a connection of the form $[F,\d]$
with $F$ symmetric, where $\d=\diag(dx_1, \cdots, dx_n)$.
\ms

By definition of flatness, we get

\refclaim[fp] Proposition.  A Darboux connection $-\d F+ F^t\d$ is
flat if and only if $F=(f_{ij})$ satisfies
\refeq[en]$$\cases{(f_{ij})_{x_j}+ (f_{ji})_{x_i} + \sum_k f_{ik}f_{jk}=0,& if $i\not=j$\cr 
(f_{ij})_{x_k} = f_{ik}f_{kj}, & if $i, j,k$ are distinct.\cr}$$
 An Egoroff connection $[F,\d]$ (with $F=F^t$) is flat if and only if $F=(f_{ij})$ is a
solution of equation \refem{}, the $n$-dimensional system associated to the symmetric space
$GL(n,R)/SO(n)$. 

Let $w=-\d F+ F^t\d$ be  a flat Darboux connection. Then a metric $ds^2=\sum_j
b_j^2(x)dx_j^2$ has $w$ as its Levi-Civita connection if and only if $(b_1, \cdots, b_n)$
is a solution of  
\refeq[era]$${(b_i)_{x_j}\over b_j}= f_{ij}, \quad i\not=j.$$
In general, given a solution $F=(f_{ij})$ of equation \refen{},  equation \refera{} has
infinitely many local solutions parametrized by $n$ functions $b_i$ defined on the
line
$x_j=0$ for $j\not=i$. These are used as the initial conditions for the ordinary
differential equations \refera{}. 
\ms

Next we will explain the relation between the space of solutions of equation \refen{} and  the
set of flat $n$-submanifolds in $R^{2n}$ with flat normal bundle and maximal rank. First we
need the following definition:

\refpar[bh] Definition.  The rank of a submanifold $M^n$ of $R^m$ at $x\in M$ is the dimension
of the space of shape operators at $x\in M$. $M$ is said to have constant rank $k$ if the
rank of $M$ at $x$ is equal to $k$ for all $x\in M$. In general, the rank $k$ of $M$ at $x$ is less
than or equal to the codimension of $M$ in $R^m$. 

\ms

Using the local theory of submanifolds, it is easy to see that (cf. [Te2]) if $M^n$ is a
flat submanifold of
$R^{2n}$ with flat normal bundle and constant rank $n$, then locally there exist a coordinate
system $x:\co\to M\subset R^{2n}$,
 $A=(a_{ij}):\co\to O(n)$,
$b_i:\co\to R$ and parallel normal frame $\{e_{n+1},\cdots, e_{2n}\}$ such that the two
fundamental forms are:
$$\cases{I=\sum_i b_i(x)^2 dx_i^2,&\cr II=\sum_{i,j=1}^n a_{ji}b_i dx_i^2 e_{n+j}.&\cr}$$
The coordinate system $x$ is unique up to permutation and changing
$x_i$ to $-x_i$ (i.e., the action of the Weyl group $B_n$). Such coordinates will be
called {\it principal curvature coordinates\/}, and
$(b,A)$ will be called the {\it fundamental data\/} of $M$. 

\refclaim[bl] Theorem ([Te2]). Suppose $M^n$ is flat submanifold of $R^{2n}$ with flat normal
bundle and constant rank $n$, $x$ is a principal curvature coordinate system, and  $(b,A)$ is
the fundamental data of $M$. Set $$f_{ij}=\cases{(b_i)_{x_j}/b_j, & if $i\not=j$,\cr 0,& if
$i=j$.\cr}$$ Then $F=(f_{ij})$ is a solution of equation \refel{}, the system associated to the
rank $n$ symmetric space 
$SO(2n)/S(O(n)\times O(n))$.  Conversely, if
$F=(f_{ij})$ is a solution of equation \refel{}, then there exist an  open
subset $\co$ of $R^n$,  $b:\co\to R^n$, $A:\co\to O(n)$ and an immersion $X:\co\to R^{2n}$ 
such that 
\refeq[np]$$\cases{dA= A(-F\d + \d F^t),& where $\d=\diag(dx_1,\cdots, dx_n)$\cr 
(b_i)_{x_j}=f_{ij}b_j,& if $i\not= j$,\cr}$$ and
\item {(i)} the immersion $X$ is flat,  has flat normal bundle and constant rank $n$,
\item {(ii)}  $x$ is a principal curvature coordinate system for $X(\co)$ and $(b,A)$ is its
fundamental data,
\item {(iii)} given any constants $c_1, \cdots, c_n$ and set $b_i=\sum_j c_ja_{ji}$ for $ 1\leq
i\leq n$.  Then $(b_1, \cdots, b_n)$
 is a solution of the second equation of system \refnp{},
\item {(iv)} let $b=(a_{11}, \cdots, a_{1n})$, then $X(\co)\subset S^{2n-1}$.  \ei

\refpar[fq] Remark. If $F=F^t$  is a solution of equation \refem{}, then $F$ is a solution of
equation \refel{}. Let $A$ be as in Theorem \refbl{} and
$b=(a_{11}, \cdots, a_{1n})$, and $X:\co\to R^{2n}$ the corresponding immersion.
  It was observed by Dajczer and Tojeiro [DaR2] that the condition $F=F^t$
is equivalent to the condition that $X(\co)$ is a Lagrangian flat submanifold  of 
$R^{2n}=C^n$. 

\refpar[fr] Remark. Let $N^n(c)$ denote the $n$-dimensional space form of sectional curvature
$c$.  It was proved by Cartan ([Ca]) that if $c<c'$ then $N^n(c)$ can not be locally isometrically 
embedded in $N^m(c')$ when $m<2n-1$, but can  if $m\geq 2n-1$.   An analogue of B\"acklund's
theorem for immersions of
$N^n(c)$ into $N^{2n-1}(c')$ was constructed by Tenenblat and the first author [TT]
for
$c=-1$ and by Tenenblat [Ten] for
$c=0$. The corresponding Gauss-Codazzi equations for these immersions are called
the {\it generalized  sine-Gordon equation (GSGE)\/} and {\it generalized wave
equation GWE\/} respectively. GSGE and GWE arise as the $n$-dimensional system
associated to the symmetric spaces $SO(2n,1)/SO(n)\times SO(n,1)$ and
$SO(2n)/S(O(n)\times O(n))$ respectively (cf. [Te2]).  Du ([Du]) noted that the
equation for isometric immersions of
$N^k(c)$ in $N^m(c')$ is the $k$-dimensional system
associated to a suitable rank $k$ symmetric space. For example, the equation for
immersions of
$R^k$ into
$S^n$, $n\geq 2k-1$, is the $k$-dimensional system associated to $Gr(k,R^{n+1})$.
Du also proved that the B\"acklund transformations constructed in [TT], [Ten], and 
Ribaucour transformations constructed in [DaR1] are given by actions of certain
order two elements in $G_-^m$. 

\bs

Darboux' orthogonal coordinate systems arises naturally in the work of Dubrovin and Novikov
([DN1, 2]) and Tsarev ([Ts]) on Hamiltonian system of hydrodynamic type.  A
brief review of their results follows.   Given a smooth section $P$ of
$L(TR^n,TR^n)$, the following first order quasi-linear system for
$\g:R^2\to R^n$
\refeq[eua]$${\p \g\over \p t} = P(\g)({\p \g\over \p x})$$
is called a {\it hydrodynamic system\/}.  If $(u_1, \cdots, u_n)$ is a local coordinate system on
$R^n$, then $P(u)({\p\over \p u_i}) = \sum_{j=1}^n v_{ij}(u) {\p \over \p u_j}$ for some smooth map
$v=(v_{ij}):R^n\to gl(n)$.  System \refeua{} is said to be {\it diagonalizable\/} if given any point $x\in
R^n$ there is a local coordinate $u$ around $x$ such that the corresponding matrix map $v$ for the
smooth section $P$ is diagonal.  Let $ds_0^2$ denote the standard flat metric on $R^n$, and $\K$ its
Levi-Civita connection. Given two functionals $F$ and $G$ on $\cs(R,R^n)$, 
$$\{F, G\}(\g)=\int_{-\infty}^\infty (\d F(\g), \K_{\g_x}(\d G(\g))dx$$ defines a Poisson structure on
$\cs(R,R^n)$.  Dubrovin and Novikov ([DN1], [DN2]) proved that this is the only Poisson structure on
$\cs(R, R^n)$ that is given by a first order differential operator.  Given a zero order Lagrangian
$F$ with density $f:R^n\to R$, i.e.,
$$F(\g)=\int_{-\infty}^\infty f(\g(x))dx,$$ the Hamiltonian equation with respect
to the Poisson structure defined above is 
\refeq[eub]$${\p \g\over \p t} = \K_{\g_x} (\K f(\g)).$$ Such system is called {\it Hamiltonian
system of hydrodynamic type\/}. Novikov conjectured that if system \refeub{} is diagonalizable
then it is completely integrable. This conjecture is proved by Tsarev in [Ts]. In
these results, the boundary conditions for the Poisson bracket are ignored 
and equation
\refeub{} is defined on an open subset of $R^n$. In other words, this is the local
theory of Hamiltonian hydrodynamic systems.  Below we state some of Tsarev's
results:

\refclaim[bn] Theorem ([Ts]).  Suppose the Hamiltonian system for 
$$F(\g)=\int_{-\infty}^\infty
f(\g(x))dx$$ is diagonalizable with respect to a local coordinate system 
$(u_1, \cdots, u_n)$. Then $\K^2 f= \sum_i v_i(u) du_i\otimes
{\p \over \p u_i}$, and the Hamiltonian system \refeub{}  is 
\refeq[euc]$${\p u_i\over \p t} = v_i(u) {\p u_i\over \p x}.$$  Moreover:
\item {(i)} $(u_1, \cdots, u_n)$ is a local orthogonal coordinate system on $R^n$, i.e., the standard
metric $ds_0^2$ on $R^n$ is $\sum_{i=1}^n b_i(u)^2 du_i^2$ for some smooth functions $b_1, \cdots,
b_n$. Moreover, 
\refeq[eue]$${1\over v_i-v_j}{\p v_i\over \p u_j} = -{1\over b_j} 
{\p b_i\over \p u_j}.$$
\item {(ii)} If $(\tilde v_1, \cdots, \tilde v_n)$ is a solution of system \refeue{},
 then ${\p u_i\over
\p t}= \tilde v_i(u) {\p u_i\over \p x}$ is also a Hamiltonian system of hydrodynamic type. 
\item {(iii)} Suppose $\sum_{i=1}^n \tilde b_i^2(u) du_i^2$ and $\sum_{i=1}^n b_i(u)^2 du_i^2$ have
the same connection $1$-form, i.e., $(b_i)_{u_j}/b_j= (\tilde b_i)_{u_j}/\tilde b_j$ for all $i\not= j$.
Set $h_i=\tilde b_i/b_i$. Then 
\refeq[eud]$${\p u_i\over \p t} = h_i(u) {\p u_i\over \p x}$$ is a Hamiltonian system of
hydrodynamic type and commutes with system \refeuc{}. 
\item {(iv)} If $v_i, \cdots, v_n$  are distinct, then system \refeuc{} is completely integrable.

 To end this section, we review some of the elementary relations between 
Dubrovin's Frobenius manifolds ([Dub2]) and the
$n$-dimensional system \refem{} associated to the symmetric space $GL(n)/SO(n)$.  For more
deep and detailed results of Frobenius manifolds, we refer the reader to Dubrovin's article
[Dub2]. 

\refpar[az] Definition ([Dub2], [Hi2]). A {\it Frobenius manifold\/} of degree $m$ (not 
necessarily an integer) is a quintuple
$(R^n,x,g,\o, \xi)$, where $x$ is a coordinate system on $R^n$,
$g=\sum_j\phi_{x_j}dx_j^2$ a flat Egoroff metric,
$\o=\sum_j\phi_{ x_j}dx_j$ and $\xi=\sum_j\phi_{x_j}dx_j^3$ for some function
$\phi$ such that $\phi, g, \o, \xi$ satisfy the following conditions:
\item {(i)} $\K \o=0$, where $\K$ is the Levi-Civita connection of $g$,
\item {(ii)} $\K \xi$ is a symmetric $4$ tensor,
\item {(iii)} $\phi_{x_j}$ is homogeneous of degree $m$ for all $j$, i.e., 
$\phi_{x_j}(rx)=r^m \phi_{x_j}(x)$ for all $r\in
R^*$ and $x\in R^n$. 

\ni The coordinate system $x$ is called a {\it canonical coordinate system\/}. 
\ms

Each tangent plane of the Frobenius manifold has a natural multiplication defined as follows:
Set $v_i=\partial /\partial x_i$. Then 
$$v_iv_j=v_jv_i= \d_{ij} v_i, \quad \forall\,\, 1\leq i, j\leq n$$
defines a multiplication on the tangent plane of $R^n$. Moreover, $T(R^n)_x$ is
a commutative algebra and
\refeq[eh]$$\o(uv)=g(u,v), \qquad \xi(u,v,w)= g(uvw).$$
The dual of the $1$-form $\o$ is $e=\sum_j \partial/\partial x_j$, which is the identity, i.e.,
$ve=ev=v$ for all $v\in T(R^n)_x$.  

The following Proposition gives the relation between Frobenius manifolds and solutions of
the $n$-dimensional system associated to $GL(n)/SU(n)$. 

\refclaim[ia] Proposition ([Dub2], [Hi2]). Let $(R^n,x, g,\o,\xi)$ be a Frobenius
manifold of degree
$m$, 
$g=\sum_j b_j^2dx_j^2$, and  $w_{ij}=f_{ij}(-dx_i+dx_j)$ the Levi-Civita connection of $g$, i.e.,
$b_i^2=\phi_{x_i}$ and $f_{ij}=(b_i)_{x_j}/b_j$. Set
$$F=(f_{ij}), \qquad  S(x)=(S_{ij}(x))=(f_{ij}(x)(x_i-x_j)).$$ Then
\item {(i)} $F=(f_{ij})$ is a solution of equation \refem{}, the $n$-dimensional system
associated to $GL(n)/SO(n)$,
\item {(ii)} $F$ is invariant under the action of $R^*$ defined in section 9, i.e., $$r\cdot
F(x)=r^{-1} F(r^{-1}x), {\rm for\, all\,\, } r\not=0.$$
\item {(iii)} ${\p S\over \p x_i}= [[F,e_{ii}], S]$, where $e_{ii}$ is the diagonal
matrix with all entries zero except the $ii$-th entry is $1$,
\item {(iv)} $(b_1, \cdots, b_n)^t$ is an eigenvector of the matrix $(S_{ij})$ with eigenvalue
$m/2$. \ei

Since Frobenius manifolds are flat, there are also  coordinate systems such that all the
coordinate vector fields are covariant constant. A coordinate system 
$(t_1,\cdots, t_n)$ on a Frobenius manifold
$(R^n,x,g,\o,\xi)$ is called a {\it flat coordinate system\/} if 
$g$ has constant coefficients with respect to the $t$-coordinates.
Since $\K \o=0$ and $e$ is the dual of $\o$, we have $\K e=0$. So there exists a flat
coordinates $(t_1, \cdots, t_n)$ such that $e=\partial /\partial t_1$.
It follows from the condition that $\K \xi$ is symmetric,  there exists a function
$h(t)$ such that  
$$\cases{g=\sum_{jk} h_{t_1t_jt_k} dt_jdt_k,&\cr
\o= dt_1,&\cr \xi=\sum_{ijk} h_{t_it_jt_k}dt_idt_jdt_k.&\cr}$$  
Using condition \refeh{} and the fact that $\partial /\partial t_1$ is the identity,
 the multiplication can be written down explicitly in terms of $h$. In order for this
multiplication to be associative, $h$ has to satisfy a complicated non-linear equation, which
is the WDVV (Witten-Dijkgraaf-Verlinde-Verlinde) equation. The
WDVV equation arises in the study of Gromov-Witten invariants, and we refer
the readers to work of Dubrovin ([Dub2]) and Ruan and Tian ([RT]).

\vfil\eject

\Bibliography

\a //ABT//Ablowitz, M.J., Beals, R., Tenenblat, K.//On the solution of the
generalized wave and generalized sine-Gordon equations//Stud.
Appl. Math.//74//1987//177-203////

\b //AC//Ablowitz, M.J.,Clarkson, P.A.//Solitons,
non-linear evolution equations and inverse scattering//\break Cambridge Univ.
Press////1991//////

\a //AKNS1//Ablowitz, M.J., Kaup, D.j., Newell, A.C. and Segur, H.//Method for
solving the sine-Gordon equation//Phy. Rev. Lett.//30//1973//1262-1264////

\a //AKNS2//Ablowitz, M.J., Kaup, D.j., Newell, A.C. and Segur, H.//The inverse
scattering transform - Fourier analysis for nonlinear
problems//Stud. Appl. Math.//53//1974//249-315////

\b //AbM//Abraham, R., Marsden, J.E.//Foundations of
mechanics// Benjamin/Cummings////1978//////

\a //Ad//Adler, M.//On a trace functional for formal pseudo-differential
operators and the symplectic structure of the Korteweg-de Vries Type
Equations//Invent. Math.//50//1979//219-248////

\a //AdM//Adler, M., van Moerbeke, P.//Completely integrable systems,
Euclidean Lie algebras and curves//Adv. Math.//38//1980//267-317////

\b //Ar//Arnold, V.I.//Mathematical methods of classical mechanics//
Springer-Verlag////1978//////

\a //BC1//Beals, R., Coifman, R.R.//Scattering and inverse scattering for
first order systems//Commun. Pure Appl. Math.//37//1984//39-90////

\a //BC2//Beals, R., Coifman, R.R.//Inverse scattering and evolution
equations//Commun. Pure Appl. Math.//38//1985//29-42////

\a //BC3//Beals, R., Coifman, R.R.//Linear spectral problems, non-linear
equations and the $\bar\partial$-method//Inverse Problems//5//1989//87-130////

\a //BDZ//Beals, R., Deift, P., Zhou, X.//The inverse scattering transform on the
line//Important development in soliton theory////1993//7-32//Springer Ser.
Nonlinear Dynam.//

\a //BS//Beals, R., Sattinger. D.H.//On the complete integrability of
complete integrable systems//Commun. Math. Phys.//138//1991//409-436////

\a //BMW//Birnir, B., McKean, H., Weinstein, A.//The rigidity of
sine-Gordon breathers//Comm. Pure Appl. Math.//47//1994//1043-1051////

\a //Bo//Bobenko, A.I.//All constant mean curvature tori in $R^3, S^3, H^3$ in terms
of theta functions//Math. Ann.//290//1991//209-245////

\a //Bu//Budagov//Completely integrable model of classical field theory with
nontrivial particle interaction in two-dimensional space-time//Dokl. Akad. Nauk
SSSR//235//1977//24-56//(in Russian)//  

\b //BuC//Bullough, R.K., Caudrey, P.J.//Solitons//Topics in Current Physics, vol. 117, 
Springer-Verlag////1980////// 

\a //BFPP//Burstall, F.E., Ferus, D., Pedit, F., Pinkall, U.//Harmonic tori in
symmetric spaces and commuting Hamiltonian systems on loop
algebras//Annals of Math.//138//1993//173-212//// 

\p //BG//Burstall, F.E., Guest, M.A.//Harmonic two-spheres in compact symmetric
spaces//////////preprint//
 
\a //Ca//Cartan, E.//Sur les vari\'et\'es de courbure constante d'un espace
euclidien ou non-euclidien//Bull. Soc. Math. France//47//1919//132-208////

\a //Ch//Cherednik, I.V.//Definition of
$\tau$-functions for generalized affine Lie
algebras//Funct. Anal. Appl.//17//1983//243-244////

\b //CM//Chernoff, P., Marsden, J.//Properties of infinite dimensional Hamiltonian
system// Lecture Notes in Math., vol. 425,  Springer-Verlag, Berlin and New
York////1974////// 

\a //dR//Da Rios////Rend. Circ. Mat. Palermo//22//1906//117-135////

\p //DaR1//Dajczer,M., Tojeiro, R.//An extension of the classical Ribaucour
transformation//////////preprint//

\p //DaR2//Dajczer,M., Tojeiro, R.//The Ribaucour transformation for
flat totally real submanifolds of complex space forms//////////preprint//

\b //Da1//G. Darboux//Lecon sur la th\'eorie
g\'en\'erale des surfaces//Chelsea////1972////3rd edition//

\b //Da2//G. Darboux//Lecon sur les Syst\`emes
Orthogonaux//Gauthier-Villars////1910////2nd edition//

\a //De//Denzler, J.// Nonpersistence of breather families for the
perturbed sine-Gordon equation//Commun. Math. Phys.//158//1993//397-430////

\p //DPW//Dorfmeister, J., Pedit, F., Wu, H.//Weierstrass type representation of
harmonic maps into symmetric spaces//////////preprint//

\a //Dr//Drinfel'd, V.G.//Hamiltonian structures on
Lie groups, Lie bialgebras and the geometrical meaning of Yang-Baxter
equations//Sov. Math. Doklady//27//1983//68-71////

\a //DS//Drinfel'd, V.G., and Sokolov, V.V.//Equations of Korteweg-de Vries type
and simple Lie algebras//Dokl. Akad. Nauk SSSR//258//1981//11-16//(Trans. as
Soviet Math. Dokl. 23, 457-462//

\p //Du//Du, X.//Isometric immersions of space forms in space forms//////////thesis,
Northeastern University, 1997//

\a //Dub1//Dubrovin, B.A.//Differential geometry of strongly integrable
systems of hydrodynamic type//Funct. Anal. Appl.//24//1990//280-285////

\a //Dub2//Dubrovin, B.A.//Geometry of 2D topological field
theories//Lecture Notes in Mathematics, vol. 1620, Springer-Verlag////1996//////

\a //DN1//Dubrovin, B.A., Novikov, S.P.//Hamiltonian formalism of
one-dimensional systems of hydrodynamic type, and the Bogolyubov-Whitham
averaging method//Soviet Math. Dokl.//27//1983//781-785////

\a //DN2//Dubrovin, B.A., Novikov, S.P.//On Poisson brackets of hydrodynamic
type//Soviet Math. Dokl.//30//1984//294-297////

\b //Ei//Eisenhart, L.P.//A treatise on the
differential geometry of curves and surfaces//Ginn////1909//////

\b //FT//Faddeev, L.D., Takhtajan, L.A.//Hamiltonian Methods in the theory
of Solitons//Springer-Verlag////1987//////

\a //FPPS//Ferus, D., Pedit, F., Pinkall, U., Sterling,
I.//Minimal tori in $S^4$//J. reine angew. Math.//429//1992//1-47////

\a //FNR1//Flaschka, H., Newell, A.C., Ratiu, T.//Kac-Moody Lie algebras and
soliton equations, II. Lax equations associated with
$A^{(1)}_1$//Physica//9D//1983//303-323////

\a //FNR2//Flaschka, H., Newell, A.C., Ratiu, T.//Kac-Moody Lie algebras and
soliton equations, IV. Lax equations associated with
$A^{(1)}_1$//Physica//9D//1983//333-345////

\a //FK//Fordy, A.P., Kulish, P.P.//Nonlinear Schr\"odinger equations and simple
Lie algebra//Commun. Math. Phys.//89//1983//427-443////

\b //FW//Fordy, A.P., Wood, J.C.//Harmonic maps and
integrable systems//Vieweg////1994////editors//

\a //FRS//Frenkel, I.E., Reiman, A.G., Semenov-Tian-Shansky, M.A.//Graded Lie
algebras and completely integrable dynamical systems//Soviet
Math. Dokl.//20//1979//811-814////

\a //GGKM//Gardner, C.S., Greene, J.M., Kruskal, M.D., Miura, R.M.//Method for
solving the Korteweg-de Vries equation//Physics Rev. 
Lett.//19//1967//1095-1097////

\a //GDi//Gel'fand, I.M., Dikii, L.A.//Fractional Powers of Operators and
Hamiltonian Systems//Funkcional-\break 'nyi Analiz i ego Prilozenjia//10//1976//////

\a //GDo//Gel'fand, I.M., Dorfman, I. Ya//Hamiltonian operators and algebraic
structures related to them//Functional Anal. Appl.//13//1979//248-261////

\a //GZ//Gu, G.H., Zhou, Z.X.//On Darboux transformations for soliton
equations in high dimensional spacetime//Letters in Math.
Phys.//32//1994//1-10////

\b //Gu//Guest, M.//Integrable systems and harmonic maps//Cambridge University 
Press////////to appear//

\a //Ha//Hasimoto, H.//Motion of a vortex filament and its relation to
elastic//J. Phys. Soc. Jap.//31//1971//293-//// 

\b //He//Helgason//Differential Geometry and Symmetric Spaces//Academic
Press////1978//////

\a //Hi1//Hitchen, N.J.//Harmonic maps from a 2-torus to the
3-sphere//J. Differential Geom.//31//1990//627-710//// 

\p //Hi2//Hitchen, N.J.//Frobenius manifolds//////////preprint//

\p //Hi3//Hitchen, N.J.//Integrable systems in Riemannian
geometry//////////preprint//
   
\b //Ka//Kac, V.G.//Infinite Dimensional Lie Algebras//Cambridge University
Press////1985////// 

\b //Kon//Konopelchenko, B.G.//Nonlinear Integrable
equations//Lecture Notes in Physics, vol. 270//Springer-Verlag//1987////// 

\a //Kos//Kostant, B.//The solution to a generalized Toda lattice and
representation theory//Adv. Math.//34//1979//195-338////

\a //KdV//Koteweg, D.J., de Vries, G.//ON the change of form of long waves
advancing in a rectangular canal, and on a new type of long stationary
waves//Philos. Mag. Ser. 5//39//1895//422-443////

\a //Kr1//Krichever, I.M.//Integration of non-linear equations by methods of
algebraic geometry//Funct. Anal. Appl.//11//1977//12-26////

\a //Kr2//Krichever, I.M.//Methods of algebraic geometry in the theory of
non-linear equations//Russian Math. Surveys//32//1977//185-213////

\a //KW//Kuperschmidt, B.A., Wilson, G.//Modifying Lax equations and the second
Hamiltonian structure//Invent. Math.//62//1981//403-436////

\a //LP//Langer, J., Perline, R.//Local geometric invariants of integrable
evolution equations//J. Math. Phys.//35//1994//1732-1737////

\a //La//Lax, P.D.//Integrals of nonlinear equations of evolution adn solitary
waves//Comm. Pure. Appl. Math.//31//1968//467-490////

\a //Lu//Lu, J.-H.//Momentum mappings and reduction of Poisson
actions//Symplectic geometry, groupoid, and integrable
systems, Math. Sci. Res. Inst. Publ.//20//1991//209-226////

\a //LW//Lu, J.-H., Weinstein, A.//Poisson Lie groups, dressing actions, and Bruhat
decompositions//J. Differential Geom.//31//1990//501-526////

\a //M//Magri, F.//A simple model of the integrable Hamiltonian
equation//J. Math. Phys.//19//1978//1156-1162////

\a //Ma//Manakov, S.V.//An example of a completely integrable nonlinear wave
field with nontrivial dynamics (Lee model)//Teor. Mat. Phys.//28//1976//172-179////

\a //Mc//McIntosh, I.//Infinite dimensional Lie groups and the two dimensional
Toda lattice//Harmonic maps and integrable systems////1994//205-220//Fordy,
A.P., Wood, J.C. editors, Vieweg//

\b //Mi//Miura, R.M.//B\"acklund
transformations//Lecture Notes in Mathematics, vol.
515, Springer-Verlag////1976////editor//

\a //Mo//Moore, J.D.//Isometric immersions of space forms in space
forms//Pacific Jour. Math.//40//1972//157-166////

\b //N//Newell, A.C.//Solitons in Mathematics and Physics//SIAM,
CBMS-NSF vol. 48////1985//////

\b //NMPZ//Novikov, S., Manakov, S., Pitaevskii, L.B. Zakharov, V.E.//Theory
of Solitons//Plenum. New York////1984//////

\p //Pa//Palais, R.S.//The symmetries of solitons//////////Rudolph-Lipschitz
lectures, University of Bonn, 1997//

\b //PT//Palais R.S., and Terng C.L.//Critical Point Theory and Submanifold Geometry//
Lecture Notes in Math., vol. 1353,  Springer-Verlag, Berlin and New
York////1988////// 

\a //PiS//Pinkall, U., Sterling, I.//On the classification of constant mean curvature
tori//Ann. of Math.//130//1989//407-451////

\b //PrS//Pressley, A. and Segal, G. B.//Loop Groups//Oxford Science Publ., 
Clarendon Press, Oxford////1986//////

\a //R//Reyman, A.G.//Integrable Hamiltonian systems connected with graded
Lie algebras//J. Sov. Math.//19//1982//1507-1545////

\a //RS//Reyman, A.G., Semenov-Tian-Shansky//Current algebras and non-linear
partial differential equations//Sov. Math., Dokl.//21//1980//630-634////

\a //Ri//Ricca, R.L.//Rediscovery of Da Rios
equations//Nature//352//1991//561-562////

\a //RT//Ruan, Y., G. Tian//A mathematical theory
of quantum cohomology//J. Differential Geom.//42//1995//259-367////

\a //Sa//Sattinger, D.H.//Hamiltonian hierarchies on semi-simple Lie
algebras//Stud. Appl. Math.//72//1984//65-86////

\a //SZ1//Sattinger, D.H.,Zurkowski, V.D.//Gauge theory of B\"acklund
transformations.I//Dynamics of infinite dimensional
systems, Nato Sci. Inst. Ser. F. Comput.
Systems Sci.//37//1987//273-300//Springer-Verlag//

\a //SZ2//Sattinger, D.H.,Zurkowski, V.D.//Gauge theory of B\"acklund
transformations.II//Physica//26D//1987//225-250////

\a //SW//Segal, G., Wilson, G.//Loop groups and equations of KdV
type//Publ. Math. IHES//61//1985//5-65////

\a //Se//Semenov-Tian-Shansky, M.A.//Dressing transformations and Poisson
group actions//Publ. RIMS Kyoto Univ.//21//1985//1237-1260////

\a //Se//Semenov-Tian-Shansky, M.A.//Classical r-matrices, Lax equations, Poisson
Lie groups, and dressing transformations//Lecture Notes in
Physics, Springer-Verlag//280//1986//174-214////

\a //Sh//Shabat, A.B.//An inverse scattering problem//Diff.
Equ.//15//1980//1299-1307////

\a //SS//Shatah, J., Strauss, W.//Breathers as homoclinic geometric
wave maps//Physics D.//99//1996//113-133////

\a //Sy//Symes, W.W.//Systems of Toda type, Inverse spectral problems, and
representation theory//Inventiones Math.//59//1980//13-51////

\a //Ten//Tenenblat, K.//B\"acklund's theorem for submanifolds of space forms and
a generalized wave equation//Boll. Soc. Brasil. Mat.//16//1985//67-92////

\a //TT//Tenenblat, K., Terng, C.L.//B\"acklund's theorem for n-dimensional
submanifolds of $R^{2n-1}$//Ann. Math.//111//1980//477-490////

\a //Te1//Terng, C.L.//A higher dimensional generalization of the sine-Gordon
equation and its soliton theory//Ann. Math.//111//1980//491-510////

\a //Te2//Terng, C.L.//Soliton equations and differential
geometry//J. Differential Geom//45//1997//407-445////

\p //TU1//Terng, C.L., Uhlenbeck, K.//B\"acklund transformations
and loop group actions//////////preprint//

\p //TU2//Terng, C.L., Uhlenbeck, K.//Homoclinic wave maps into
compact symmetric spaces//////////in preparation//

\p //TU3//Terng, C.L., Uhlenbeck, K.//Hamiltonian theory of the
geometric non-linear Schr\"odinger equation on a Hermitian
symmetric space//////////in preparation//

\a //Ts//Tsar\"ev, S.P.//The geometry of Hamiltonian systems of
hydrodynamic type. The generalized Hodograph
method//Math. USSR Izvestiya//37//1991//397-419////

\a //U1//Uhlenbeck, K.//Harmonic maps into Lie group (classical solutions of the
Chiral model)//J. Differential Geometry//30//1989//1-50////

\a //U2//Uhlenbeck, K.//On the connection between harmonic maps and the self-dual
Yang-Mills and the sine-Gordon equations//Geometry \& Physics//2//1993//////

\a //Wa//Wadati, M.//The modified Korteweg-de Vries
equation//J. Phys. Soc. Japan//34//1973//380-384////

\a //Wi//Wilson, G.//The modified Lax equations and two dimensional Toda lattice
equations associated with simple Lie algebras//Ergodic Theory and
Dynamical Systems I//30//1981//361-380////

\a //ZF//Zakharov, V.E., Faddeev, L.D.//Korteweg-de Vries equation, a
completely integrable Hamiltonian system//Func. Anal.
Appl.//5//1971//280-287////

\a //ZMa1//Zakharov, V.E., Manakov, S.V.//On resonant interaction of wave packets
in non-linear media//JETP Letters//18//1973//243-247////

\a //ZMa2//Zakharov, V.E., Manakov, S.V.//The theory of resonant interaction of wave
packets in non-linear media//Sov. Phys. JETP//42//1974//842-850////

\a //ZMi1//Zakharov, V.E., Mikhailov, A.V.//Example of nontrivial interaction of
solitons in two-dimensional classical field theory//JETP
Letters//27//1978//42-46////

\a //ZMi2//Zakharov, V.E., Mikhailov, A.V.//Relativistically invariant
2-dimensional models of field theory which are integrable by means of the
inverse scattering problem method//Soviet Physics JETP//47//1978//1017-1027////

\a //ZS1//Zakharov, V.E., Shabat, A.B.//Exact theory of two-dimensional
self-focusing and one-dimensional of waves in nonlinear
media//Sov. Phys. JETP//34//1972//62-69////

\a //ZS2//Zakharov, V.E., Shabat, A.B.//Integration of non-linear equations of
mathematical physics by the inverse scattering method, II//Funct. Anal.
Appl.//13//1979//166-174////

\a //ZK//Zabusky, N.J., Kruskal, M.D.//Interaction of solitons in a
collisionless plasma and the recurrence of initial states//Physics Rev.
Lett.//15//1965//240-243////

\a //Zh1//Zhou, X.//Riemann Hilbert problem and inverse
scattering//SIAM J. Math. Anal.//20//1989//966-986////

\a //Zh2//Zhou, X.//Direct and inverse scattering transforms with arbitrary
spectral singularities//Comm. Pure. Appl. Math.//42//1989//895-938////

\def\uhlenaddr{{\hsize=2in \vsize=1in
\hbox{\vbox{
\leftline{Karen Uhlenbeck}
\leftline{Department of Mathematics}
\leftline{The University of Texas at Austin}
\leftline{RLM8.100 Austin, Texas 78712}
\leftline{email:uhlen@math.utexas.edu}
}}}}

\def\cltaddr{{\hsize=2in \vsize=1in
\hbox{\vbox{
\leftline{Chuu-lian Terng}
\leftline{Department of Mathematics}
\leftline{Northeastern University}
\leftline{Boston, MA  02115}
\leftline{email: terng@neu.edu}
}}}}

\bigskip\bigskip
\leftline{\cltaddr \hfill \uhlenaddr}

\end